\documentclass[twoside,a4paper,12pt,centertags,reqno]{amsart} % displays better with 13pt 
\usepackage{amsmath,amssymb,verbatim,vmargin}
\usepackage{color}
\usepackage{tikz}
\usepackage{pgfplots}
\usepackage{tikz-cd}
\usetikzlibrary{matrix,calc}
\usepackage{hyperref}
\hypersetup{colorlinks,bookmarks=true,linktocpage=true,citecolor=blue,linkcolor=magenta}

\usepackage{eulervm}   % formulas round   

\allowdisplaybreaks % align break
\usepackage{color}
\usepackage[all]{xy} % diagrams

\usepackage{mathtools}
\usepackage{MnSymbol} %wedge righthalfcup

\theoremstyle{plain}
\newtheorem{thm}{Theorem}[section]
\newtheorem{lem}[thm]{Lemma}

\newtheorem{prop}[thm]{Proposition}

\theoremstyle{definition}
\newtheorem{defn}[thm]{Definition}
\newtheorem{rem}[thm]{Remark}

\usepackage{enumerate}
\usepackage{enumitem}

\newcommand{\bRn}{\mathbb{R}^n}
\newcommand{\bRN}{\mathbb{R}^N}

\newcommand{\ep}{\varepsilon} 

\newcommand{\pd}{\partial}

\newcommand{\bR}{{\mathbb R}}

\newcommand{\bN}{{\mathbb N}}

\newcommand{\bL}{{\mathbb L}}

\newcommand{\cD}{{\mathcal D}}

\newcommand{\cK}{{\mathcal K}}
\newcommand{\cL}{{\mathcal L}}

\newcommand{\cQ}{{\mathcal Q}}
\newcommand{\cR}{{\mathcal R}}

\newcommand{\cV}{{\mathcal V}}

\newcommand{\EP}{\mathcal{E}} %\cE conflict

\newcommand{\fH}{{\mathbf H}} 
\newcommand{\fh}{{\mathbf h}}

\newcommand{\fk}{{\mathbf k}}

\newcommand{\fr}{{\mathbf r}}

\newcommand{\wt}{\widetilde}
\newcommand{\wh}{\widehat}

\newcommand{\supp}{\text{{\rm supp}}\,}
\newcommand{\dist}{\text{{\rm dist}}\,}

\newcommand{\divv}{{\text{{\rm div}}}} % \div already defined

\def\barint_#1{\mathchoice
            {\mathop{\vrule width 6pt
height 3 pt depth -2.5pt
                    \kern -9.5pt
\intop \kern -4pt}\nolimits_{#1}}%
            {\mathop{\vrule width 5pt height
3 pt depth -2.6pt
                    \kern -6.5pt
\intop \kern -4pt}\nolimits_{#1}}%
            {\mathop{\vrule width 5pt height
3 pt depth -2.6pt
                    \kern -6pt
\intop \kern -4pt}\nolimits_{#1}}%
            {\mathop{\vrule width 5pt height
3 pt depth -2.6pt
          \kern -6pt \intop \kern -4pt}\nolimits_{#1}}}
          
           \def\bariint_#1{\mathchoice
            {\mathop{\vrule width 15pt
height 3 pt depth -2.5pt
                    \kern -15.8pt
\intop \kern -8pt\intop \kern -4pt}\nolimits_{#1}}%
            {\mathop{\vrule width 9pt height
3 pt depth -2.6pt
                    \kern -10.5pt
\intop \kern -8pt\intop \kern -4pt}\nolimits_{#1}}%
            {\mathop{\vrule width 9pt height
3 pt depth -2.6pt
                    \kern -10pt
\intop \kern -8pt\intop \kern -4pt}\nolimits_{#1}}%
            {\mathop{\vrule width 9pt height
3 pt depth -2.6pt
          \kern -8pt \intop \kern -10pt\intop \kern -4pt}
      \nolimits_{  #1}}}

\def\barintlim_#1{\mathchoice
            {\mathop{\vrule width 6pt
height 3 pt depth -2.5pt
                    \kern -8.8pt
\intop \kern -4pt}\limits_{#1}}%
            {\mathop{\vrule width 5pt height
3 pt depth -2.6pt
                    \kern -6.5pt
\intop \kern -4pt}\limits_{#1}}%
            {\mathop{\vrule width 5pt height
3 pt depth -2.6pt
                    \kern -6pt
\intop \kern -4pt}\limits_{#1}}%
            {\mathop{\vrule width 5pt height
3 pt depth -2.6pt
          \kern -6pt \intop \kern -4pt}\limits_{#1}}}
          
           \def\bariintlim_#1{\mathchoice
            {\mathop{\vrule width 15pt
height 3 pt depth -2.5pt
                    \kern -15.8pt
\intop \kern -8pt\intop \kern -4pt}\limits_{#1}}%
            {\mathop{\vrule width 9pt height
3 pt depth -2.6pt
                    \kern -10.5pt
\intop \kern -8pt\intop \kern -4pt}\limits_{#1}}%
            {\mathop{\vrule width 9pt height
3 pt depth -2.6pt
                    \kern -10pt
\intop \kern -8pt\intop \kern -4pt}\limits_{#1}}%
            {\mathop{\vrule width 9pt height
3 pt depth -2.6pt
          \kern -8pt \intop \kern -10pt\intop \kern -4pt}
      \limits_{  #1}}}
          
\renewcommand{\iint}{\int \kern -3pt\int}       

\numberwithin{equation}{section}
\makeatletter
\@namedef{subjclassname@2020}{\textup{2020} Mathematics Subject Classification}
\makeatother

\usepackage[utf8]{inputenc}
\usepackage{chngcntr}
\usepackage{apptools}
\AtAppendix{\counterwithin{lemma}{section}}

\theoremstyle{plain} % slanted text
\newtheorem*{theorem*}{Theorem A}

\title[Gradient profile for vortex reconnection]
{Gradient profile for the reconnection of vortex lines with the boundary in type-II superconductors}  

\author{Yi C. Huang} 
\address{Universit\'e Sorbonne Paris Nord, Institut Galil\'ee, LAGA, CNRS (UMR 7539), F-93430 Villetaneuse, France}
\address{School of Mathematical Sciences, Nanjing Normal University, Nanjing 210023, People's Republic of China}
\email{Yi.Huang.Analysis@gmail.com}
\urladdr{https://orcid.org/0000-0002-1297-7674}
\author{Hatem Zaag} 
\address{Universit\'e Sorbonne Paris Nord, Institut Galil\'ee, LAGA, CNRS (UMR 7539), F-93430 Villetaneuse, France}
\email{Hatem.Zaag@univ-paris13.fr}
\urladdr{https://www.math.univ-paris13.fr/~zaag/}

\date{\today} 

\keywords{Vortex lines, type-II superconductors, reconnection with the boundary, finite time quenching/extinction, 
singularity formation, self-similar variables, shrinking set, mode dynamics, finite dimensional reduction, blowup/extinction profile, cusps, gradient profile.}
\subjclass[2020]{Primary 35K50, 35B40; Secondary 35K55, 35K57.}  
\thanks{The first author (YCH) is supported by the National NSF grant of China (no. 11801274). 
This paper is completed while YCH is on leave, funded by CSC Postdoctoral/Visiting Scholar Program (no. 202006865011), at Universit\'e Sorbonne Paris Nord.
YCH would like to thank Profs P. Auscher and H. Yin for early influences towards parabolic regularity and singularity formation of NLWs.}

\begin{document}

\begin{abstract} 
In a recent work, Duong, Ghoul and Zaag determined the gradient profile for blowup solutions of standard semilinear heat equation with power nonlinearities
in the (supposed to be) generic case.
Their method refines the constructive techniques introduced by Bricmont and Kupiainen and further developed by Merle and Zaag.
In this paper, we extend their refinement to the problem about the reconnection of vortex lines with the boundary in a type-II superconductor under planar approximation, 
a physical model derived by Chapman, Hunton and Ockendon featuring the finite time quenching for the nonlinear heat equation
$$\frac{\pd h}{\pd t}=\frac{\pd^2 h}{\pd x^2}+e^{-h}-\frac{1}{h^\beta},\quad\beta>0$$
subject to initial boundary value conditions 
$$h(\cdot,0)=h_0>0,\quad h(\pm1,t)=1.$$
We derive the intermediate extinction profile with refined asymptotics, and with extinction time $T$ and extinction point $0$, 
the gradient profile behaves as $x\rightarrow0$ like 
$$\lim_{t\rightarrow T}\,(\nabla h)(x,t)\quad\sim\quad\frac{1}{\sqrt{2\beta}}\frac{x}{|x|}\frac{1}{\sqrt{|\log|x||}}
\left[\frac{(\beta+1)^2}{8\beta}\frac{|x|^2}{|\log|x||}\right]^{\frac{1}{\beta+1}-\frac12},$$
agreeing with the gradient of the extinction profile previously derived by Merle and Zaag.
Our result holds with general boundary conditions and in higher dimensions.
\end{abstract}

\maketitle

\tableofcontents

\section{Introduction}

Let $\Omega$ be a bounded domain in $\bRN$ with $N\geq1$. Let $\beta>0$.
In this paper, we are interested in the quenching problem modelled on the nonlinear heat equation on $\Omega$
\begin{equation}\label{e:Quen}
\frac{\pd h}{\pd t}=\Delta h-F(h),
\end{equation}
where 
$$F=F_\beta=\frac{1}{h^{\beta}}+\wt F\in C^\infty(\bR_+), \quad\text{with\,\,} \bR_+=(0,\infty).$$ 
We shall assume that $\wt F$ satisfies
\begin{equation}\label{e:QuenF1}
\wt F(h)=o\left(\frac{1}{h^{\beta}}\right)\quad\text{and}\quad \wt F'(h)=o\left(\frac{1}{h^{\beta+1}}\right),\quad\text{as}\,\, h\rightarrow0,
\end{equation}
and $h$ is subject to initial data $h_0=h(\cdot,0)>0$ and the Dirichlet boundary condition $h\equiv1$ on $\pd\Omega$.
Finite time \textit{quenching} or \textit{extinction} of a solution $h$ to the Cauchy problem for \eqref{e:Quen} at $x_0\in\Omega$, 
means that for some $T\in\bR_+$, $h$ has limit value 0 at $(x_0,T)$.
We remark that the perturbative assumption \eqref{e:QuenF1} is typically satisfied by 
$$\wt F \equiv0, \quad\wt F(h)=e^{-h}, \quad or\quad \wt F(h)=\frac{1}{h^{\beta'}}\,\,\text{ for some}\,\, \beta'<\beta.$$
Moreover, since our main concern in this paper is in the quenching scenario $h\rightarrow0$,
apart from the smoothness $F\in C^\infty(\bR_+)$ there is no need (and no gain) to prescribe the behaviour of $F(h)$ for $h\geq1$.
Equation \eqref{e:Quen} can also be considered on the whole space $\bRN$, 
if $F$ further satisfies the decay assumption
\begin{equation} \label{e:QuenF2}
|F(h)|+|F'(h)|\leq Ce^{-h}\quad\text{as}\quad h\rightarrow+\infty, 
\end{equation}
and 
if the boundary condition is replaced by a growth condition, say,
\begin{equation} \label{e:Quen-DirGr}
h(x,t)\sim {a_1|x|}\quad\text{as}\quad |x|\rightarrow+\infty
\end{equation}
for some $a_1>0$.
We remark that \eqref{e:QuenF2}-\eqref{e:Quen-DirGr} are again not within our main concern.

Now, by introducing the following family of transformations
\begin{equation}\label{e:hu}
u(x,t)=u(\alpha,x,t)=\frac{\alpha^{\frac{\alpha}{\beta+1}}}{h(x,t)^\alpha} \quad (\alpha>0)
\end{equation}
for \eqref{e:Quen} we are then led to study the blowup problem of the following nonlinear heat equation with power-type nonlinearity and a particular perturbation term
\begin{equation}\label{e:Grad}
\begin{aligned}
&\frac{\pd u}{\pd t}=\Delta u-a\frac{|\nabla u|^2}{u}+f(u), \\
&f(u)=\alpha^{\frac{\beta}{\beta+1}}u^{1+\frac{1}{\alpha}}F(\alpha^{\frac{1}{\beta+1}}u^{-\frac{1}{\alpha}})=u^p+\wt f(u),
\end{aligned}
\end{equation}
where $u>0$, 
$u\equiv1$ on $\pd\Omega$, 
$(a,p)$ is computed from $(\alpha,\beta)$ by
\begin{equation}\label{e:apbeta}
1<a=1+\frac{1}{\alpha}<p=\frac{1+\alpha+\beta}{\alpha},
\end{equation}
and $\wt f\in C^\infty(\bR_+)$ satisfies
\begin{equation}\label{e:GradF1}
\wt f(u)=o(u^p)\quad\text{and}\quad \wt f'(u)=o(u^{p-1})\quad\text{as}\quad u\rightarrow+\infty.
\end{equation}
If $\Omega=\bRN$ we further assume that $f$ satisfies 
\begin{equation}\label{e:GradF2}
|f(u)|+|f'(u)|\leq Cu^{1+\frac{1}{\alpha}}\exp(-\alpha^{\frac{1}{\beta+1}}u^{-\frac{1}{\alpha}})\quad\text{as}\quad u\rightarrow0,
\end{equation}
and 
the solution $u$ is subject to the decay condition 
\begin{equation}\label{e:Grad-Dir}
u(x,t)\sim \frac{\wt a_1}{|x|}\quad\text{as}\quad |x|\rightarrow+\infty
\end{equation}
for some $\wt a_1>0$.
Conversely, given $1<a<p<\infty$ for \eqref{e:Grad}, 
we recover by $a=1+\frac{1}{\alpha}$ and $p=\frac{1+\alpha+\beta}{\alpha}$ the exponent $\beta$ for \eqref{e:Quen}
and the transformation index $\alpha$ in \eqref{e:hu}.
If we use, instead of $\beta$, the $p$ and $a$ notations to denote objects for $h$, we mean $\alpha=1$.

In the non-perturbed case of \eqref{e:Grad} ($a=\wt f=0$), namely, for the equation
\begin{equation}\label{e:Grad0}
\frac{\pd u}{\pd t}=\Delta u+|u|^{p-1}u,
\end{equation}
where $u(\cdot,t): \bRN\rightarrow\bR$, $p>1$ and $p<\frac{N+2}{N-2}$ if $N\geq3$,
Fujita \cite{Fuj66}, Levine \cite{Lev73}, Ball \cite{Bal77} and Kavian \cite{Kav87} obtained obstructions to global existence in time, using monotonicity properties and the maximum principle.
Another method has been followed by Merle and Zaag \cite{MerZaa97DMJ}
(see also the earlier papers by Giga and Kohn \cite{GigKoh85, GigKoh87, GigKoh89}, Berger and Kohn \cite{BerKoh88}, Herrero and Vel\'azquez \cite{HerVel92, HerVel93},
 Bricmont and Kupiainen \cite{BriKup94} and Zaag \cite{Zaa98}).
Once an asymptotic profile 
(that is a function that $u(t)$ approaches to, after a time-dependent rescaling,  as $t\rightarrow T<+\infty$) is derived formally,
the existence of a solution $u(t)$ which blows up at time $T$ with the suggested profile is then proved rigorously,
using analysis of \eqref{e:Grad0} near the given profile and reduction of the existence problem to a finite-dimensional one.
More precisely, in \cite{MerZaa97DMJ} the following asymptotic behaviour 
\begin{equation} \label{e:MZ97DMJ}
\left\|(T-t)^{\frac{1}{p-1}}u(\cdot,t)-\Phi_0\left(\frac{\cdot-x_0}{\sqrt{(T-t)|\log(T-t)|}}\right)\right\|_{L^\infty(\Omega)}\longrightarrow0\quad\text{as}\quad t\rightarrow T,
\end{equation}
is justified, where $x_0$ is the only blowup point, $T$ is the blowup time, and
\begin{equation} \label{e:Phi0}
\Phi_0(z)=\frac{1}{\left(p-1+\frac{(p-1)^2}{4p}|z|^2\right)^{\frac{1}{p-1}}}.
\end{equation}
Let $\kappa=\Phi_0(0)=(p-1)^{-\frac{1}{p-1}}$. The so-called final blowup profile for \eqref{e:Grad0}
$$\lim_{t\rightarrow T}u(x,t)\quad\sim\quad\left[\frac{(p-1)^2}{8p}\frac{|x-x_0|^2}{|\log|x-x_0||}\right]^{-\frac{1}{p-1}}\quad\text{as}\quad x\rightarrow x_0,$$
was later derived by Zaag in \cite{Zaa98}.

Note that the constructive method in \cite{MerZaa97DMJ} is robust 
and in fact has been extended by Zaag and collaborators to a considerably large class of parabolic problems such as:
Merle and Zaag \cite{MerZaa97NL} for the quenching problems 
(see also Duong and Zaag \cite{DuoZaa19} for MEMS devices and Duong, Kavallaris and Zaag \cite{DuoKavZaa21} for Gierer-Meinhardt system),
Masmoudi and Zaag \cite{MasZaa08}, Nouaili and Zaag \cite{NouZaa18}, Duong, Nouaili and Zaag \cite{DuoNouZaa19MAMS, DuoNouZaa20} for the complex Ginzburg-Landau equation,
Ebde and Zaag \cite{EbdZaa11}, Nguyen and Zaag \cite{NguZaa16}, Tayachi and Zaag \cite{TayZaa19}, Duong, Nguyen and Zaag \cite{DuoNguZaa19}, 
Ghoul, Nguyen and Zaag \cite{GhoNguZaa17} and Abdelhedi and Zaag \cite{AbdZaa21JDE, AbdZaa21DCDS} for perturbed nonlinear source terms (with or without gradient perturbation), 
Nouaili and Zaag \cite{NouZaa15} for a non-variational complex valued heat equation,
and Ghoul, Nguyen and Zaag \cite{GhoNguZaa18AIHP, GhoNguZaa18JDE} for semilinear parabolic systems where the nonlinearities have no gradient structure.

In a broader context, this method has also proved to be efficient for different PDEs of different type,
and no list can be exhaustive (see Merle, Rapha\"el and Szeftel \cite{MerRapSze20} for anisotropic type I blowup for heat equation,
del Pino, Musso and Wei \cite{DelPinMusWei20}, Harada \cite{Har20}, Collot, Merle and Rapha\"el \cite{ColMerRap20} for type II blowup for heat equation,
Rapha\"el and Schweyer \cite{RapSch14},
Collot, Ghoul, Masmoudi and Nguyen \cite{ColGhoMasNgu21} for Keller-Segel system,
Krieger, Schlag and Tataru \cite{KriSchTat09}, C\^ote and Zaag \cite{CotZaa13}, Donninger and Sch\"orkhuber \cite{DonSch14} for wave equation,
Martel \cite{Mar05} for generalized Korteweg-de Vries equation,
Elgindi, Ghoul and Masmoudi \cite{ElgGhoMas19} for incompressible Euler,
Buckmaster, Shkoller and Vicol \cite{BucShkVic22} for 2D isentropic compressible Euler,
Merle, Rapha\"el, Rodnianski and Szeftel \cite{MerRapRodSze21,MerRapRodSze19} for Schr\"odinger equation and 3D compressible fluids, etc.).

Back to Equation \eqref{e:Quen} (with $0<h<1$ and $F(h)=\frac{1}{h^\beta}$), few results are known except the case on $(-1,1)\subset\bR$ with Dirichlet boundary conditions,
where some criteria of quenching are obtained and, upon assuming that the initial data $h_0$ satisfies
$$\frac{\pd^2 h_0}{\pd x^2}-\frac{1}{h_0^\beta}\leq0,$$ 
the quenching rates and quenching profiles are obtained, 
see Levine \cite{Lev89}, Guo \cite{Guo91} (also Fila, Hulshof and Quittner \cite{FilHulQui92}), Filippas and Guo \cite{FliGuo93} and the references therein.
By adapting the techniques of \cite{BriKup94} and \cite{MerZaa97DMJ},
Merle and Zaag \cite{MerZaa97NL} succeeded in constructing a stable quenching solution.

Recently, Duong, Ghoul and Zaag \cite{DuoGhoZaa21} refined the constructive techniques of \cite{MerZaa97DMJ} and obtained the gradient profile of blowup solutions.   
Our aim here is to determine the gradient profile of the quenching solution constructed in \cite{MerZaa97NL}.
To recall the known results and also state ours, let us introduce some definitions.

\begin{defn}[Intermediate profile] \label{d:intfile}
Introduce the intermediate profile
\begin{equation}\label{e:Phihat}
\wh\Phi(z)=\left(\beta+1+\frac{(\beta+1)^2}{4\beta}|z|^2\right)^{\frac{1}{\beta+1}},\quad z\in\bRN,
\end{equation}
with the intermediate gradient profile being
\begin{equation}\label{e:Phihatn}
\nabla\wh\Phi(z)=\frac{\beta+1}{2\beta}z\left(\beta+1+\frac{(\beta+1)^2}{4\beta}|z|^2\right)^{-\frac{\beta}{\beta+1}},\quad z\in\bRN.
\end{equation}
Here $z$ will be a time-dependent rescaling of $x$ as in \eqref{e:MZ97DMJ}-\eqref{e:Phi0} above.
\end{defn}

Now we prescribe the final extinction profile (near extinction point) that is expected in our construction.
For later purpose we also prescribe its global behaviour.

\begin{defn}[Final profile] \label{d:H*}
Let $x_0\in\Omega\subset\bRN$.
Define $H^*_{x_0}\in C^\infty(\Omega\backslash\{x_0\})$ by 

\bigskip

(i) If $\Omega=\bRN$, then 
\begin{equation} \label{e:H*0}
H^*_{x_0}(x)=H^*(x-x_0),
\end{equation}
where $H^*\in C^\infty(\bRN\backslash\{0\})$ is defined as follows: 
\begin{equation}\label{e:H*}
H^*(x)=\left[\frac{(\beta+1)^2}{8\beta}\frac{|x|^2}{|\log|x||}\right]^{\frac{1}{\beta+1}},\quad 0<|x|\leq \min\left\{C(a_1,\beta),\frac12\right\},
\end{equation}
\begin{equation}\label{e:H*line}
H^*(x)=a_1|x|,\quad |x|\geq1,
\end{equation}
with $a_1$ given in \eqref{e:Quen-DirGr} and with the additional requirements that 
$$\forall\,\, x\neq 0, \quad H^*(x)>0\quad\text{and}\quad |\nabla H^*(x)|>0.$$

\bigskip

(ii) If $\Omega$ is bounded, then with $\varrho_0=\dist(x_0,\pd\Omega)$,
\begin{equation}\label{e:H*'}
H^*_{x_0}(x)=\left[\frac{(\beta+1)^2}{8\beta}\frac{|x-x_0|^2}{|\log|x-x_0||}\right]^{\frac{1}{\beta+1}}, \quad0<|x-x_0|\leq \min\left\{C(\beta), \frac{\varrho_0}{4}\right\},
\end{equation}
$$H_{x_0}^*(x)=1,\quad|x-x_0|\geq\frac{\varrho_0}{2},$$
with the additional requirements that 
$$\forall\,\, x\neq x_0, \quad H_{x_0}^*(x)>0\quad\text{and}\quad |\nabla H_{x_0}^*(x)|>0.$$
\end{defn}

\begin{rem}
Since our analysis will be local, we shall not precise
$$H^*(x)\quad\text{when}\quad\min\left\{C(a_1,\beta),\frac12\right\}\leq |x|\leq1;$$
$$H^*_{x_0}(x)\quad\text{when}\quad\min\left\{C(\beta), \frac{\varrho_0}{4}\right\}\leq|x-x_0|\leq \frac{\varrho_0}{2}.$$
Similarly, \eqref{e:H*line} is imposed only in accordance with \eqref{e:Quen-DirGr}, 
itself introduced in accordance with the vortex line model considered in \cite{ChaHunOck98}.
In view of this, it is more reasonable to assume $H^*(x)=a_\beta|x|^{\frac{2}{\beta+1}}$ for $|x|\geq1$ instead of \eqref{e:H*line}.
These different growth functions will not make any difference in the \textit{local Cauchy theory} addressed in Appendix \ref{a:C}.
Thus we relieve ourselves from further precision on the profile beyond \eqref{e:H*} and \eqref{e:H*'}.
For the same reason and when $\Omega$ is bounded, 
other types of boundary behaviours instead of $H_{x_0}^*(x)=1$ for $|x-x_0|\geq\frac{\varrho_0}{2}$ can be considered.
\end{rem}

\begin{rem}
For $\beta=1$, the (vortex line) profile $H^*(x)$ in \eqref{e:H*}-\eqref{e:H*line} is nearly straight (up to a logarithm $\log|x|$).
For $\beta>1$, the profile forms a cusp at $x=0$.
\end{rem}

\begin{rem}[Gradient profile]
The gradient of $H^*$ via \eqref{e:H*} satisfies that as $x\rightarrow0$, 
\begin{equation}\label{e:nablaH*}
\nabla H^*(x)\quad\sim\quad\frac{1}{\sqrt{2\beta}}\frac{x}{|x|}\frac{1}{\sqrt{|\log|x||}} 
\left[\frac{(\beta+1)^2}{8\beta}\frac{|x|^2}{|\log|x||}\right]^{\frac{1}{\beta+1}-\frac12}.
\end{equation}
When $|x|\rightarrow0$, $|\nabla H^*(x)|$ blows up if $\beta>1$ and extinguishes if $0<\beta\leq1$. 
\end{rem}

As the containing space of initial data $h_0$ and the evolution time-slices $h(\cdot,t)$ ($0\leq t<T$), 
$h$ being a solution to \eqref{e:Quen}, we shall need the following notion.

\begin{defn}[Trajectory space for quenching problem] \label{d:fH}
Define $\fH$ as follows:

\bigskip

(i) If $\Omega=\bRN$, then for some $x_0\in\bRN$, we set
\begin{equation} \label{e:HRN}
\fH=\fH_{\psi_{x_0}}=\left\{\fh: (\fh-\psi_{x_0})\in W^{1,\infty}(\bRN),\,\,\frac{1}{\fh}\in L^\infty(\bRN)\right\}, 
\end{equation}
where $\psi_{x_0}\in C^\infty(\bRN)$ satisfies 
$$\psi_{x_0}\equiv0,\,\,\text{if}\,\,|x-x_0|\leq1;$$
$$\psi_{x_0}(x)=a_1|x-x_0|,\,\,\text{if}\,\,|x-x_0|\geq2.$$ 
Here, $a_1$ is the constant given in \eqref{e:Quen-DirGr}.

\bigskip

(ii) If $\Omega$ is a bounded domain in $\bRN$, then we set
\begin{equation} \label{e:Hbdd}
\fH=\left\{\fh\in W^{1,\infty}(\Omega): \frac{1}{\fh}\in L^\infty(\Omega)\right\}.
\end{equation}
\end{defn}

When no ambiguity arises in the context,
we shall no longer precise the reference point $x_0$ in the definition of $\fH$.
The following was proved in \cite{MerZaa97NL}.

\bigskip 

\begin{theorem*}[Existence and stability of vortex reconnection with the boundary 
or finite time quenching obeying prescribed behaviours given in Definitions \ref{d:intfile} and \ref{d:H*}, Merle and Zaag] \label{t:MZ}
Assume that either $\Omega$ is bounded, or $\Omega=\bRN$ and \eqref{e:QuenF2}-\eqref{e:Quen-DirGr} hold.

\bigskip

(1) (Existence) For all $x_0\in\Omega$, there exists a positive $h_0\in\fH$ such that for a $T>0$,
equation \eqref{e:Quen} with initial data $h_0$ has a unique solution $h(\cdot,t)$ on $[0,T)$ satisfying 
$$\lim_{t\rightarrow T}h(x_0,t)=0.$$ 

Furthermore,

(i) (Intermediate extinction profile)
\begin{equation} \label{e:implicit}
\lim_{t\rightarrow T}\,\,\left\|\frac{(T-t)^{\frac{1}{\beta+1}}}{h(\cdot,t)}-\frac{1}{\wh\Phi(z_{x_0}(\cdot,t))}\right\|_{L^\infty(\Omega)}=0,
\end{equation}
where $\wh\Phi$ is given in \eqref{e:Phihat} and $$z_{x_0}(x,t)=\frac{x-x_0}{\sqrt{(T-t)|\log(T-t)|}}.$$ 

(ii) (Final extinction profile)
$h^*(x):=\lim_{t\rightarrow T}h(x,t)$ exists for all $x\in\Omega$ and
$$h^*(x)\sim H^*_{x_0}(x)\quad\text{as}\quad x\rightarrow x_0,$$
where $H^*_{x_0}$ is given precisely in \eqref{e:H*0}, \eqref{e:H*} and \eqref{e:H*'}.

\bigskip

(2) (Stability). For every $\ep>0$, there exists a neighbourhood $\cV_0$ of $h_0$ in $\fH$ with the following property:

\bigskip

For each $\wt h_0\in \cV_0$, there exists $\wt x_0\in\Omega$ and $\wt T>0$ satisfying
$$|T- \wt T|+|x_0-\wt x_0|\leq\ep$$
such that \eqref{e:Quen} with initial data $\wt h_0$ has a unique solution $\wt h(\cdot,t)$ on $[0, \wt T)$ satisfying 
$$\lim_{t\rightarrow  \wt T}\wt h(\wt x_0,t)=0.$$ 

\bigskip

Furthermore,

(i) (Intermediate extinction profile)
$$\lim_{t\rightarrow \wt T}\,\,\left\|\frac{(\wt T-t)^{\frac{1}{\beta+1}}}{\wt h(\cdot,t)}-\frac{1}{\wh\Phi(\wt z_{x_0}(\cdot,t))}\right\|_{L^\infty(\Omega)}=0,$$
where $\wh\Phi$ is given in \eqref{e:Phihat} and $$\wt z_{x_0}(x,t)=\frac{x-\wt x_0}{\sqrt{(\wt T-t)|\log(\wt T -t)|}}.$$

(ii) (Final extinction profile)
$\wt h^*(x):=\lim_{t\rightarrow  \wt T }\wt h(x,t)$ exists for all $x\in\Omega$ and
$$\wt h^*(x)\sim H^*_{\wt x_0}(x)\quad\text{as}\quad x\rightarrow \wt x_0,$$
where $H^*_{\wt x_0}$ is given precisely in \eqref{e:H*0}, \eqref{e:H*} and \eqref{e:H*'}.
\end{theorem*}

\bigskip

In this paper we shall prove that the final profile is stable under differentiation:
$$\begin{aligned}
 (\nabla h)^*(x)&\sim\nabla h^*(x)\\ 
\bigg(\lim_{t\rightarrow T} \nabla h(x,t)&\sim\nabla\lim_{t\rightarrow T}  h(x,t)\bigg)
\end{aligned}$$ 
as $x\rightarrow x_0$. More precisely, we shall prove the following result.

\begin{thm}[Construction of finite time extinction solutions with refined asymptotics and determination of the final extinction profile and the final gradient profile] \label{t:Quen}
Assume that either $\Omega$ is bounded, or $\Omega=\bRN$ and \eqref{e:QuenF2}-\eqref{e:Quen-DirGr} hold.

\bigskip

For all $x_0\in\Omega$, there exists a positive $h_0\in\fH$ such that for a $T=T(h_0)\in(0,e^{-1})$,
equation \eqref{e:Quen} with initial data $h_0$ has a unique solution $h(\cdot,t)$ on $[0,T)$ satisfying 
$$\lim_{t\rightarrow T}h(x_0,t)=0.$$ 

\bigskip

Furthermore, for this $h$ and for all $t\in[0,T)$, we have:

\bigskip

(i) (Intermediate extinction profile and refined asymptotics) 
\begin{equation} \label{e:refined1}
 \left\|\frac{(T-t)^{\frac{1}{\beta+1}}}{h(\cdot,t)}-\frac{1}{\wh\Phi(z_{x_0}(\cdot,t))}\right\|_{L^\infty(\Omega)}\leq C\frac{\log(|\log(T-t)|)}{|\log(T-t)|}
\end{equation}
and for each $K>0$ and $\Omega_{t,K}=\left\{x\in\Omega:|x-x_0|\leq K\sqrt{(T-t)|\log(T-t)|}\right\}$,
\begin{equation} \label{e:refined2} 
 \begin{aligned}
 \left\|(T-t)^{-\frac{1}{\beta+1}+\frac12}\nabla h(\cdot,t)-\frac{(\nabla\wh\Phi)(z_{x_0}(\cdot,t))}{|\log(T-t)|^{\frac12}}\right\|_{L^\infty(\Omega_{t,K})}
\leq C(K)\frac{\log(|\log(T-t)|)}{|\log(T-t)|}.
\end{aligned}
\end{equation}
Here, $\wh\Phi$ is given in \eqref{e:Phihat} and $z_{x_0}(x,t)$ is given in Theorem A.

\bigskip

(ii) (Final extinction profile and final gradient profile)

\bigskip

\indent\indent (ii-1) $h^*(x):=\lim_{t\rightarrow T} h(x,t)$ exists for all $x\in\Omega$, and 
\begin{equation} \label{e:profsame} 
h^*(x)\sim H^*_{x_0}(x)\quad\text{as}\quad x\rightarrow x_0,
\end{equation} 
where $H^*_{x_0}$ is given precisely in \eqref{e:H*0}, \eqref{e:H*} and \eqref{e:H*'}.

\bigskip

\indent\indent (ii-2) $(\nabla h)^*(x):=\lim_{t\rightarrow T} \nabla h(x,t)$ exists for all $x\in\Omega\backslash\{x_0\}$, and
\begin{equation} \label{e:profgrad} 
(\nabla h)^*(x)\sim \nabla H^*_{x_0}(x)\quad\text{as}\quad x\rightarrow x_0.
\end{equation}
\end{thm}

\begin{rem} \label{r:Quen}
The main contribution of the theorem is to discover precise error estimates \eqref{e:refined1}-\eqref{e:refined2} and their consequence \eqref{e:profgrad}.
Note that \eqref{e:implicit} suffices for the derivation of \eqref{e:profsame}, which is same as in Theorem A.
The error estimate in \eqref{e:refined1} refines $|\log(T-t)|^{-\frac12}$ which is implicit in \eqref{e:implicit} (see \cite{MerZaa97NL}).
As for \eqref{e:refined2}, the error term one obtains from the estimates of \cite{MerZaa97NL} would be $C|\log(T-t)|^{-\frac12}$,
which is of the same size as $\frac{\nabla\wh\Phi}{\sqrt{|\log(T-t)|}}$, preventing us from obtaining an equivalent for $\nabla h$.
Thus, the improved error estimate we get in \eqref{e:refined2} thanks to our strategy is crucial. 
\end{rem}

\begin{rem}
We obtained indeed an estimate stronger than \eqref{e:refined2}: $\forall\,t\in[0,T)$,
\begin{equation} \label{e:refined2'} 
 \begin{aligned}
& \left\|\frac{(T-t)^{\frac{1}{\beta+1}+\frac12}\nabla h(\cdot,t)}{h^2(\cdot,t)}-\frac{(\nabla\wh\Phi)(z_{x_0}(\cdot,t))}{|\log(T-t)|^{\frac12}\wh\Phi^2(z_{x_0}(\cdot,t))}\right\|_{L^\infty(\Omega)}\\
&\qquad\qquad\leq C\frac{\log(|\log(T-t)|)}{|\log(T-t)|}.
\end{aligned}
\end{equation}
See Proposition \ref{p:Grad}.
\end{rem}

\begin{rem} \label{r:Quen'}
In fact, $h$ extinguishes at time $T$ only at $x_0$, see Proposition \ref{p:profile}.
\end{rem}

\begin{rem}[Stability of the gradient profile] \label{r:stab}
We use the method which was first introduced by Bressan in \cite{Bre90, Bre92} for the heat equation with an exponential source,
then for \eqref{e:Grad0} in Bricmont and Kupiainen \cite{BriKup94} and Merle and Zaag \cite{MerZaa97DMJ} 
(also Fermanian, Merle and Zaag \cite{FerKamMerZaa00}).
This consists in a formal approach where the profile is obtained through an inner/outer expansion with matching asymptotics,
followed by a rigorous proof where the PDE is linearized around the profile candidate.
Then, the negative part of the spectrum (which is infinite dimensional) is controlled
thanks to the decaying properties of the Laplacian,
whereas the nonnegative part (which is finite dimensional) is controlled thanks to the degree theory.
Via this finite dimensional reduction and the interpretation of the parameters of the finite-dimensional problem
in terms of the blowup time and the blowup point,
the refined estimates \eqref{e:refined1}-\eqref{e:refined2} and the final profiles \eqref{e:profsame}-\eqref{e:profgrad} 
are stable under perturbations (in the same form) of initial data.
\end{rem}

\bigskip

\textbf{Structure}. 
In next section we recall the physical derivation of the vortex line model.
Then in Section \ref{s:form}, we review the dynamical system formulation of \eqref{e:Grad} 
together with some standard tools in the construction of blowup solutions with prescribed behavior,
and introduce the localization scheme of \eqref{e:Quen} that is particular for the quenching problems. 
We prepare the notion of shrinking set in Section \ref{s:ShriSet}. 
The proof of Theorem \ref{t:Quen} will be outlined in Section \ref{s:Map}, 
with technical ingredients proved in the seven sections together with the first two appendices that follow.
Appendix \ref{a:C} justifies the local well-posedness of the Cauchy problem for \eqref{e:Quen} in the space $\fH$.

\bigskip

\textbf{Notations}. Throughout this paper, $C$ is a constant which depends only on $N$, $\beta$ and $\alpha$,
and whose value may change from line to line in deriving various estimates.
If we need a constant depending on other parameters, we will specify its dependence.

\section{Reconnection of vortex lines with the boundary} 

Before going further we explain the motivation for studying this quenching problem.
Indeed, the nonlinear equation \eqref{e:Quen} appears in various physical contexts (combustion and population genetics for example). 
For $N=1$ and $\beta=1$, 
it is also the natural lab model employed by Chapman, Hunton and Ockendon to describe the magnetic field in a type-II superconductor in $\bR^3$,
which develops a particular type of line singularity called \textit{vortex}, 
see \cite[Section 3.2]{ChaHunOck98}.
In general, a vortex is not situated in a plane,
but under some reasonable physical conditions, 
the planar approximation is relevant.
In this case, a \textit{vortex line} at time $t\geq0$ can be viewed as 
$$L(t)=\{(x,y,z)=(x,0,h(x,t)): x\in\Omega\subset\bR\}.$$
The physical derivation (see \cite[(25)]{ChaHunOck98}) gives the following vortex line dynamics
\begin{equation}\label{e:Quen1}
\frac{\pd h}{\pd t}=\frac{1}{4\pi}\frac{\pd^2 h}{\pd x^2}+H_0e^{-h}-\frac{1}{4\pi}K_{1}(h),
\end{equation}
where $(H_0,0,0)$ is the applied magnetic field assumed to be constant,
and the smooth nonlinear function $K_1$ satisfies the following asymptotics 
$$K_1(h)\sim\frac{1}{h}\,\,\text{as}\,\,h\rightarrow0\quad\text{and}\quad K_1(h)\sim 2e^{-2h}\,\,\text{as}\,\,h\rightarrow\infty.$$
Thus, the boundary term $H_0e^{-h}$ dominates for $h$ large, while the image term $\frac{1}{4\pi}K_{1}(h)$ dominates for $h$ small.
Moreover, $H_0e^{-h}-\frac{1}{4\pi}K_{1}(h)$ in \eqref{e:Quen1} multiplied by the normalization constant $4\pi$ satisfies \eqref{e:QuenF1} for $\beta=1$.
Quenching of a solution $h$ to \eqref{e:Quen1} is then the \textit{vortex reconnection} with the boundary, namely, the $z=0$ plane in $\bR^3$.

Finally, we remark that the boundary growth \eqref{e:QuenF2} (hence \eqref{e:GradF2} for \eqref{e:Grad}) is no longer relevant if we consider solutions $h$ small.
The physical model becomes
$$\frac{\pd h}{\pd t}=\frac{1}{4\pi}\frac{\pd^2 h}{\pd x^2}-\frac{1}{4\pi}\frac1h.$$
See \cite[Section 3.2.2]{ChaHunOck98} for details.

\section{Dynamical system formulation and localization} \label{s:form} 

Following \cite{MerZaa97DMJ} and \cite{Zaa98}, we work with two rescaling transformations on \eqref{e:Quen} and \eqref{e:Grad}.
The first one reformulates the $t\rightarrow T$ asymptotic behaviours for solutions of \eqref{e:Grad} as a problem of long time dynamics,
while the second aims at deriving final profiles via localizing solutions of \eqref{e:Quen} away from the extinction point.

\bigskip

\subsection{Dynamical system formulation}

Note that for $F(h)=\frac{1}{h^\beta}$, 
$$h(t)=[(\beta+1)(T-t)]^{\frac{1}{\beta+1}}\quad\text{and}\quad u(t)=[(p-1)(T-t)]^{-\frac{1}{p-1}}$$
are respectively the ODE/flat solutions for \eqref{e:Quen} and \eqref{e:Grad}.
In fact, the former ODE on $h$ is solved subject to the natural extinction condition $h(T)=0$.

\bigskip

\underline{First rescaling transformation}. 
As in Giga and Kohn \cite{GigKoh85, GigKoh87, GigKoh89} and Merle and Zaag \cite{MerZaa97DMJ, MerZaa97NL},
we introduce the transform
\begin{equation} \label{e:utow}
w(y,s)=w_{T}(y,s)=(T-t)^{\frac{1}{p-1}}u(x,t),
\end{equation}
where $(y,s)$ are the so-called \textit{similarity variables} defined by
$$y=\frac{x}{\sqrt{T-t}}\quad\text{and}\quad s=-\log(T-t).$$
Here $s>0$ when $T<1$. The equation satisfied by $w$ is then
\begin{equation} \label{e:w}
\frac{\pd w}{\pd s}=\Delta w-\frac12y\cdot\nabla w-\frac{w}{p-1}-a\frac{|\nabla w|^2}{w}+w^p+e^{-\frac{ps}{p-1}}\wt f(e^{\frac{s}{p-1}}w),
\end{equation}
where $\wt f$ is in \eqref{e:Grad}.
Moreover, \eqref{e:implicit} for $x_0=0$ is then equivalent to
\begin{equation} \label{e:qtozero*}
\lim_{s\rightarrow\infty}\left\|w(\cdot,s)-\Phi\left(\frac{\cdot}{\sqrt s}\right)\right\|_{L^\infty(\bRN)}=0,
\end{equation}
where
\begin{equation}\label{e:Phi}
\Phi(z)=\Phi_a(z)=\left(p-1+\frac{(p-1)^2}{4(p-a)}|z|^2\right)^{-\frac{1}{p-1}}.
\end{equation}
See for example \cite[Remark 1.1]{MerZaa97DMJ} for the formal derivation of the profile $\Phi_0$.

\begin{rem} \label{r:exten}
We abused a little bit the notation $L^\infty(\bRN)$ in \eqref{e:qtozero*}, 
if $\Omega$ is bounded.
Since we are concerned with the $t\rightarrow T$ dynamics, 
it is however reasonable to interpret, via extension by zero of the function to be normed (e.g., $w(\cdot,s)-\Phi\left(\frac{\cdot}{\sqrt s}\right)$ in \eqref{e:qtozero*}), 
the $L^\infty$ norm with respect to the variable $y=\frac{x}{\sqrt{T-t}}$ over the $\frac{1}{\sqrt{(T-t)}}$-dilation of $\Omega$ as $L^\infty(\bRN)$ norm. 
In the sequel we shall denote these norms simply by $\|\cdot\|_{L^\infty}$.
\end{rem}

Now we can reformulate for \eqref{e:Grad} the results analogous to Theorem \ref{t:Quen}.

\begin{prop}[Construction of finite time blowup solutions with refined asymptotics 
and determination of the final blowup profile and the final gradient profile] \label{p:Grad}
There exists initial data $u_0\in W^{1,\infty}(\Omega)$ such that equation \eqref{e:Grad}, 
subject to \eqref{e:Grad-Dir} if $\Omega=\bRN$, has a unique solution $u$
which blows up in finite time $T=T(u_0)\in(0,e^{-1})$ only at the origin.
In particular, the following assertions hold:

\bigskip

(i) (Intermediate blowup profile with refined asymptotics) For all $t<T$,
\begin{equation} \label{e:u-int-p1}
\left\|(T-t)^{\frac{1}{p-1}}u(\cdot,t)-\Phi(z(\cdot,t))\right\|_{L^\infty(\Omega)}\leq C\frac{\log(|\log(T-t)|)}{|\log(T-t)|}
\end{equation}
and
\begin{equation} \label{e:u-int-p2}
\left\|(T-t)^{\frac{1}{p-1}+\frac12}\nabla u(\cdot,t)-\frac{(\nabla\Phi)(z(\cdot,t))}{\sqrt{|\log(T-t)|}}\right\|_{L^\infty(\Omega)}\leq C\frac{\log(|\log(T-t)|)}{|\log(T-t)|},
\end{equation}
where
\begin{equation}  \label{e:zPhi}
z(x,t)=\frac{x}{\sqrt{(T-t)|\log(T-t)|}},\quad \Phi(z)=\frac{1}{(p-1+b|z|^2)^{\frac{1}{p-1}}}
\end{equation}
and
$$b=b(a)=\frac{(p-1)^2}{4(p-a)}.$$

\bigskip

(ii) (Final blowup profile and final gradient profile) For all $x\neq0$,
$u(x,T):=\lim_{t\rightarrow T}u(x,t)$ and $(\nabla u)(x,T):=\lim_{t\rightarrow T}\nabla u(x,t)$ exist. Moreover,
$$u(x,T)\sim u^*(x),\quad x\rightarrow0,$$
$$(\nabla u)(x,T)\sim \nabla u^*(x),\quad x\rightarrow0,$$
where
$$u^*(x)=\left[\frac{b}{2}\frac{|x|^2}{|\log|x||}\right]^{-\frac{1}{p-1}},$$
$$ \nabla u^*(x)=-\frac{\sqrt{2b}}{p-1}\frac{x}{|x|}\frac{1}{\sqrt{|\log|x||}}\left[\frac{b}{2}\frac{|x|^2}{|\log|x||}\right]^{-\frac{p+1}{2(p-1)}}.$$
\end{prop}

\begin{rem}
To prove Theorem \ref{t:Quen}, one takes $\alpha=1$, $a=2$ and $p=\beta+2$ in Proposition \ref{p:Grad}.
For the reverse implication, one has to keep the parameter $\alpha>0$.
\end{rem}

\begin{rem}
Note that by \eqref{e:apbeta}, Proposition \ref{p:Grad} applies to the range $1<a<p<\infty$.
From the Remark after \cite[Proposition 1]{MerZaa97NL},
the case $a<1$ (in particular, $a=0$) for \eqref{e:Grad} with $\wt f=0$, can be reduced to Duong, Ghoul and Zaag \cite{DuoGhoZaa21} via
$$\pd_t v=\Delta v+v^{p'},$$ 
$$v=(1-a)^{\frac{1-a}{p-1}}u^{1-a},\quad p'=\frac{p-a}{1-a}>1.$$
The case $a=1$ reduces to the nonlinear heat equation with exponential nonlinearity:
$$\pd_t v=\Delta v+e^v,\quad v=(p-1)\log u,$$
and the question of determining the corresponding gradient blowup profile was raised in \cite[Remark 1.9]{GhoNguZaa17}. 
We shall address this issue in a forthcoming work.
\end{rem}

To better describe the asymptotic behaviours in \eqref{e:qtozero*}, we introduce 
\begin{equation} \label{e:qdefn}
q=w-\varphi,
\end{equation}
where 
$$\varphi(y,s)=\Phi\left(\frac{y}{\sqrt s}\right)+\frac{N\kappa}{2(p-a)s},$$
and one adds $\frac{N\kappa}{2(p-a)s}$ to simplify the calculations.
Observe that $q$ satisfies
\begin{equation} \label{e:qeqn}
\frac{\pd q}{\pd s}=\cL_V(q)+B(q)+T(q)+R+L(q),
\end{equation}
where
\begin{equation} \label{e:qeqn'}
\begin{aligned}
\cL_V&=\cL+V,\quad \cL=\Delta-\frac12y\cdot\nabla+1, \\
V&=p\left(\varphi^{p-1}-\kappa^{p-1}\right), \quad \kappa=\frac{1}{(p-1)^{\frac{1}{p-1}}},\\
B(q)&=(\varphi+q)^p-\varphi^p-p\varphi^{p-1}q, \\
T(q)&=-a\frac{|\nabla(\varphi+q)|^2}{\varphi+q}+a\frac{|\nabla\varphi |^2}{\varphi}, \\
R&=-\frac{\pd \varphi}{\pd s}+\Delta \varphi-\frac12y\cdot\nabla \varphi-\frac{\varphi}{p-1}-a\frac{|\nabla \varphi|^2}{\varphi}+\varphi^p, \\
L(q)&=e^{-\frac{ps}{p-1}}\wt f\left(e^{\frac{s}{p-1}}(\varphi+q)\right).
\end{aligned}
\end{equation}
Here $\wt f$ is given in \eqref{e:Grad}. 
In $B(q)$, $T(q)$ and $L(q)$, note that $\varphi+q=w>0$.
 
\bigskip

\textit{With this transformation, the aim is to find a solution $q$ to \eqref{e:qeqn} such that
\begin{equation} \label{e:qtozero}
\lim_{s\rightarrow\infty}\|q(s)\|_{L^\infty}=0,
\end{equation}
with the precise $\frac{\log s}{s}$-order vanishing rate as in \eqref{e:refined1}.}

\bigskip

To understand the involved dynamics of $q$, 
we first make some comments on the linear, nonlinear, remainder and lower order terms that appear in \eqref{e:qeqn'}.
For $B(q)$, $T(q)$ and $L(q)$, it is helpful to have in mind that via \eqref{e:qtozero} $q$ is assumed to be small.

\bigskip

$(\clubsuit1)$ Linear operator $\cL$.
Consider the Hilbert space with Gaussian measure  
$$L^2_\rho=L^2(\bRN,\rho dy),\quad\text{where\,\,} \rho(y)=\frac{e^{-\frac{|y|^2}{4}}}{(4\pi)^{\frac N2}},$$ 
with the inner product
$$(f,g)_\rho=\int_{\bRN}f(y)g(y)\rho(y)dy.$$
The adjoint operator $(\pd_i)^*$, with respect to $(\cdot,\cdot)_\rho$, of $\pd_i=\frac{\pd}{\pd y_i}$ is $\frac{y_i}{2}-\pd_i$, hence
$$\bL=\Delta-\frac12y\cdot\nabla=-\sum_{i=1}^N(\pd_i)^*\pd_i$$
is nonpositive, self-adjoint and generates the Ornstein-Uhlenbeck semigroup $\{e^{\theta\bL}\}_{\theta\geq0}$:
$$e^{\theta\bL}r(y)=\int_{\bRN}e^{\theta\bL}(y,x)r(x)dx,$$
$$e^{\theta\bL}(y,x)=\frac{1}{(4\pi(1-e^{-\theta}))^\frac{N}{2}}\exp\left[-\frac{|ye^{-\frac{\theta}{2}}-x|^2}{4(1-e^{-\theta})}\right].$$

Similarly, the operator $\cL=\Delta-\frac12y\cdot\nabla+1$ that appears in \eqref{e:qeqn} is also self-adjoint in $L^2_\rho$.
Moreover, it has explicit spectrum as follows
\begin{equation} \label{e:spectrum}
\text{Spec}(\cL)=\left\{\lambda_m=1-\frac m2\bigg|m\in\bN\right\}.
\end{equation}
Corresponding to the eigenvalue $\lambda_m$, we have the eigenspace $\EP_m$ given by
\begin{equation} \label{e:eigen}
\EP_m=\text{Span}\left\{ h_{m_1}(y_1)h_{m_2}(y_2)\cdots h_{m_N}(y_N)|m_1+m_2+\cdots+m_N=m\right\},
\end{equation}
where $h_\ell$ is the (rescaled) Hermite polynomial in one dimension, defined by
$$h_\ell(\xi)=\sum_{j=0}^{[\ell/2]}(-1)^j\frac{\ell!}{j!(\ell-2j)!}\xi^{\ell-2j},$$
and the first three such polynomials are $1$, $\xi$ and $\xi^2-2$.

\bigskip

$(\clubsuit2)$ Potential $V$.
By \eqref{e:V1}-\eqref{e:V2} in Lemma \ref{l:V}, 
$V$ is bounded and $V(\cdot,s)\rightarrow0$ in $L^2_\rho$ as $s\rightarrow\infty$, 
hence $\cL_V$ in $L^2_\rho$ can be regarded as a perturbation of $\cL$ by $V$.

Note that when restricted to the region $|y|\geq K_0\sqrt s$ or $|z|\geq K_0$, the difference
$$V(\cdot,s)-\left(-\frac{p}{p-1}\right)=\frac{p}{(p-1)+\frac{(p-1)^2}{4(p-a)}|z|^2}$$
is small when $K_0$ is large enough,
hence $\cL_V$ behaves like $\cL-\frac{p}{p-1}$.
But $\frac{p}{p-1}>1$ and the largest eigenvalue of $\cL$ is $1$,
thus $\cL-\frac{p}{p-1}$ has a purely negative spectrum and can be handled without difficulties in the region $|y|\geq K_0\sqrt s$.

\bigskip

The above analysis on $\cL$ and $V$ suggests the following two definitions.

\begin{defn} \label{d:inout}
We introduce the following decomposition, for any $r\in L^\infty(\bRN)$:
\begin{equation} \label{e:inout}
\begin{aligned}
r(y)&=\chi(y,s)r(y)+(1-\chi(y,s))r(y)\\
&=:r_b(y,s)+r_e(y,s),
\end{aligned}
\end{equation}
where
\begin{equation} \label{e:chi}
\chi(y,s)=\chi_0\left(\frac{|y|}{K_0\sqrt{s}}\right),
\end{equation}
$\chi_0$ being a one-dimensional cut-off satisfying
\begin{equation} \label{e:chi0}
\supp \chi_0\subset[0,2],\quad 0\leq\chi_0\leq1\quad\text{and}\quad \chi_0\equiv1\quad\text{on}\quad [0,1].
\end{equation}
\end{defn}

We refer to \eqref{e:inout} as \textit{inner-outer decomposition} since 
\begin{equation} \label{e:supprbe}
\supp r_b(s)\subset \{|y|\leq 2K_0\sqrt s\}\quad\text{and}\quad
\supp r_e(s)\subset \{|y|\geq K_0\sqrt s\}.
\end{equation}
At the blowup point $x=0$, $(0,t)$ for $t<T$ is always inside (resp. outside) the support of the first (resp. second) function.
The subscripts mean ``blowup" and ``exterior".

Next we note that the set of eigenfunctions of $\cL$ makes a basis of $L^2_\rho$. 

\begin{defn} \label{d:specde}
We write $r_b\in L^2_\rho$ into the following spectral decomposition
\begin{equation} \label{e:rbsum}
r_b(y,s)=r_0(s)+r_1(s)\cdot y+y^T\cdot r_2(s)\cdot y-2\text{Tr}(r_2(s))+r_-(y,s),
\end{equation}
where
\begin{equation} \label{e:r_m}
r_m(s)=\left\{P_\beta[r_b(s)]\right\}_{\beta\in\bN^N,\, |\beta|=m}, \,\,m\geq0,
\end{equation}
with $P_\beta[r_b]$ being the projection of $r_b$ on the eigenfunction $h_\beta$
\begin{equation} \label{e:Pvect}
P_\beta[r_b(s)]=\int_{\bRN}r_b(y,s)\frac{h_\beta}{\|h_\beta\|^2_{L^2_\rho}}\rho(y)dy,
\end{equation}
and
\begin{equation} \label{e:rnega}
r_-(y,s)=P_-[r_b(s)]=\sum_{\beta\in\bN^N,\,|\beta|\geq3}P_\beta[r_b(s)]h_\beta(y).
\end{equation}
Similarly, we write
\begin{equation} \label{e:rdecPerp}
r_\perp(y,s)=P_\perp[r_b(s)]=\sum_{\beta\in\bN^N,\,|\beta|\geq2}P_\beta[r_b(s)]h_\beta(y).
\end{equation}
\end{defn}

Throughout the paper we shall use the following expansion
\begin{equation} \label{e:rdec}
\begin{aligned}
r(y)&=r_0(s)+r_1(s)\cdot y+y^T\cdot r_2(s)\cdot y-2\text{Tr}(r_2(s))\\
&\qquad+r_-(y,s)+r_e(y,s),
\end{aligned}
\end{equation}
which is nothing but a combination of \eqref{e:inout} and \eqref{e:rbsum}. Note that $L^\infty(\bRN)\subset L^2_\rho$.
Moreover, the same convention as in Remark \ref{r:exten} applies if $\Omega$ is a bounded domain.

\begin{rem}
Note that the decompositions \eqref{e:inout} and \eqref{e:rdec} depend, through the cut-off function defined in \eqref{e:chi0}, on $K_0$.
Moreover, all the decompositions \eqref{e:inout}, \eqref{e:rbsum} and \eqref{e:rdec} depend on $s$,
and when we apply them to the slice $r(\cdot,s)$, there is readily a natural parameter $s$.
We shall specify it if we do not follow this rule.
\end{rem}

We continue with the analysis on the nonlinear part in \eqref{e:qeqn}.

\bigskip

$(\clubsuit3)$ Nonlinear term $B(q)$. It is superlinear, in the sense that it satisfies
$$|B(q)|\leq C|q|^{\overline p},$$
where $\overline p=\min(2,p)>1$.
In particular, $B(q)(s)$ decays to $0$ faster than $q(s)$ as $s\rightarrow\infty$.
See Lemma \ref{l:Bq} for relevant refined estimates on $B(q)$ that depend on $q$.

\bigskip

$(\clubsuit4)$ Gradient term $T(q)$. When $q(s)$ is small ($s$ is large enough), 
the gradient term $T(q)$ satisfies (see Lemma \ref{l:Tq} in Appendix \ref{a:A}) the following control in the support of $ \chi(\cdot,s)$ (i.e. the blowup region): 
$$|\chi(y,s)T(q)(s)|\leq C\chi(y,s)\left(\frac{|q|}{s}+\frac{|\nabla q|}{\sqrt s}+|\nabla q|^2\right).$$
This estimate is implicit in the proof of \cite[Lemma B.4]{MerZaa97NL} and follows from the local boundedness of the inverse of the profile, namely,
$$\chi(y,s)\frac{1}{\varphi(y,s)}\leq C,$$
together with Taylor expansion of $T(\theta q)$ for $\theta\in[0,1]$.
As one may see from the above estimate, a bound on $|\nabla q|$ may be necessary to control $T(q)$ in the blowup region.
However, we will not impose such a bound, since we can derive it from a bound on $q$ thanks to parabolic regularity. 
Now, regarding the control of $T(q)$ in the exterior region,
because of the spatial unboundedness of the inverse of the profile,
we can not proceed as in the blowup region.
That control requires a subtle analysis in the original variable $x$ (see the proof of Lemma \ref{l:Tq}),
which involves a delicate control of the solution in three different regions of the space, and not just in two (blowup/exterior regions).
Let us recall from \cite{MerZaa97NL} the following definition.

\bigskip

\begin{defn} \label{d:R123} 
For $T\in(0,1)$, $t\in(0,T)$ and $x\in \Omega$, let
$$s=-\log(T-t)>0,\quad y=\frac{x}{\sqrt{T-t}}\quad\text{and}\quad z=\frac{y}{\sqrt s}=\frac{x}{\sqrt{(T-t)|\log(T-t)|}}.$$
We define, for $K_0>0$, $\ep_0>0$ and $t\in(0,T)$ given, three regions of $x$ that cover $\Omega$:

\begin{equation} \label{e:R1t}
\begin{aligned}
\cR_1(K_0,\ep_0,t)&=\left\{x\in \Omega\bigg| |x|\leq K_0\sqrt{(T-t)|\log(T-t)|}\right\}\\
&=\{x\in \Omega||y|\leq K_0\sqrt{s}\}=\{x\in \Omega||z|\leq K_0\},
\end{aligned}
\end{equation}

\begin{equation} \label{e:R2t}
\begin{aligned}
\cR_2(K_0,\ep_0,t)&=\left\{x\in \Omega\bigg|\frac{K_0}{4}\sqrt{(T-t)|\log(T-t)|}\leq |x|\leq \ep_0\right\}\\
&=\left\{x\in \Omega\bigg| \frac{K_0}{4}\sqrt{s}\leq|y|\leq \ep_0e^{\frac s2}\right\}=\left\{x\in \Omega\bigg|\frac{K_0}{4}\leq |z|\leq \ep_0\frac{e^{\frac s2}}{\sqrt{s}}\right\},
\end{aligned}
\end{equation}

\begin{equation} \label{e:R3t}
\begin{aligned}
\cR_3(K_0,\ep_0,t)&=\left\{x\in \Omega\bigg| |x|\geq \frac{\ep_0}{4}\right\}\\
&=\left\{x\in \Omega\bigg| |y|\geq \frac{\ep_0}{4}e^{\frac s2}\right\}=\left\{x\in \Omega\bigg| |z|\geq \frac{\ep_0}{4}\frac{e^{\frac s2}}{\sqrt{s}}\right\}.
\end{aligned}
\end{equation}
We set $\cR_i(t)=\cR_i(K_0,\ep_0,t)$, and

\begin{equation} \label{e:Ri}
\cR_i=\cR_i(K_0,\ep_0)=\{(x,t)\in \Omega \times(0,T)|x\in \cR_i(K_0,\ep_0,t)\}.
\end{equation}
\end{defn}

\begin{rem}
Note that we slightly abused the notation $\cR_3(K_0,\ep_0,t)$, which is in fact independent of $t$,
but can be understood in the time-varying context as 
$$\left\{x\in \Omega\bigg| |x|\geq \frac{\ep_0}{4}\right\}\times\{t\}\subset \Omega \times(0,T).$$
Since we are concerned with the scenario as $0\leq t\rightarrow T$,
given $K_0>0$ and $\ep_0>0$, one can assume without loss of generality that $T$ is small, hence, $\cR_2(K_0,\ep_0,t)$ is non-empty for all $t<T$.
The same convention applies to $\cR_3(K_0,\ep_0,t)$ in the case where $\Omega$ is bounded, 
as one can work with $\ep_0>0$ which is small compared to $\dist(0,\pd\Omega)$.
\end{rem}

As we mentioned before Definition \ref{d:R123},
we shall combine estimates on $|\nabla q|$ thereby on $T(q)$ from the above three regions.
See Lemma \ref{l:Tq} for relevant precise estimates on $T(q)$ that depend on $q$ and $\nabla q$.

\bigskip

At last we comment on the remainder and the lower order terms in \eqref{e:qeqn}.

\bigskip

$(\clubsuit5)$ Remainder term $R$. It is small and decays to zero. In particular,
$$\|R(\cdot,s)\|_{L^\infty}\leq\frac{C}{s},\quad\forall\,s\geq1.$$
This is reasonable since $\varphi$ is an \textit{approximate solution} to \eqref{e:w}, 
in the sense that it approximates the solution $w$ via \eqref{e:qtozero}. 
See Lemma \ref{l:R} for further estimates on $R$.

\bigskip

$(\clubsuit6)$ Lower order term $L(q)$. It is bounded by $Ce^{-\frac{ps}{p-1}}$. See Lemma \ref{l:wtf}.

\bigskip

\subsection{Localization away from the singularity and the boundary}

In order to derive the final extinction profile and the final gradient profile,
following \cite{Zaa98} and \cite{MerZaa97NL} we shall localize the study of \eqref{e:Quen} away from the extinction point.

\bigskip

\underline{Second rescaling transformation}. To the aim just noted above, we introduce
\begin{equation} \label{e:kx}
k_x(\xi,\tau)=\frac{h\left(x+\sqrt{T-t(x)}\xi,t(x)+(T-t(x))\tau\right)}{(T-t(x))^{\frac{1}{\beta+1}}},
\end{equation}
where $x\in\Omega\backslash\{0\}$ and $t(x)<T$ is determined by the quasi-parabola  
\begin{equation} \label{e:tx}
\frac{K_0}{4}\sqrt{(T-t(x))|\log(T-t(x))|}=|x|,
\end{equation}
which is exactly the relation giving the inner boundary of the region \eqref{e:R2t}.
For the sake of convenience, we denote the radial function (abusing a little bit the notation)
\begin{equation} \label{e:thetax}
\theta(x)=T-t(x)=\theta(|x|),
\end{equation}
hence 
$$t(x)\rightarrow T, \quad\theta(x)\rightarrow0\quad\text{and}\quad|\log\theta(x)|\rightarrow+\infty,\quad \text{as\,\,} x\rightarrow0.$$
Moreover, for $t(x)\geq0$, by a simple calculation $T<e^{-1}$ implies that $|x|$ is bounded, 
in this sense our localization is also away from the boundary of $\Omega$.

It is an immediate consequence of \eqref{e:Quen} and \eqref{e:thetax} that $k=k_x$ satisfies the equation
\begin{equation} \label{e:keqn}
\frac{\pd k}{\pd \tau}=\Delta_\xi k-\theta(x)^{\frac{\beta}{\beta+1}}F\left(\theta(x)^{\frac{1}{\beta+1}}k\right),
\end{equation}
where 
$$N_\beta(k):=\theta(x)^{\frac{\beta}{\beta+1}}F\left(\theta(x)^{\frac{1}{\beta+1}}k\right)\sim\frac{1}{k^\beta} \quad\text{as}\quad\theta(x)^{\frac{1}{\beta+1}}k\rightarrow0.$$
We shall show that $k(\xi,\tau)$ behaves, when 
$$|\xi|\leq\alpha_0\sqrt{|\log(T-t(x))|} \quad\text{and}\quad \tau\in\left[\frac{t_0-t(x)}{T-t(x)},1\right)$$
for some $0\leq t_0<T$ and $0<\alpha_0\leq\frac{4}{K_0}$, like the function 
\begin{equation} \label{e:khatsol}
\wh k(\tau)=\left((\beta+1)(1-\tau)+\frac{(\beta+1)^2}{4\beta}\frac{K_0^2}{16}\right)^{\frac{1}{\beta+1}}.
\end{equation}
Note that $\wh k(\tau)$ is the solution of the following ODE
\begin{equation} \label{e:khat}
\frac{d\wh k}{d \tau}=-\frac{1}{\wh k^{\beta}}
\end{equation}
on $[0,\wh T)$ with $$\wh T=\wh T(\beta,K_0)=1+\frac{(\beta+1)K_0^2}{64\beta}>1,$$
subject to the initial value condition
\begin{equation} \label{e:khat0}
\wh k(0)=\wh\Phi\left(\frac{K_0}{4}\right).
\end{equation}
Note that at $\tau=0$, this \textit{flat} solution $\wh k$ has to agree with $k_x(\xi,0)$ for $\xi=0$.
Therefore, $\wh k(0)$ has to further agree with $\frac{h(x,t)}{(T-t)^{\frac{1}{\beta+1}}}$ for $t=t(x)$,
that is, in the boundary of $\cR_2$ and, more importantly, in the region $\cR_1$.
It is therefore natural to impose the initial value condition \eqref{e:khat0} in view of the profile approximation ansatz.

\begin{rem}
We give another \textit{a posteriori} explanation to motivate the initial value condition \eqref{e:khat0}. 
This will also be instrumental in later arguments. 
Analogous to \eqref{e:kx} we also transform the profile $\wh\Phi(z)$ in \eqref{e:Phihat} as follows
\begin{equation} \label{e:phix}
\wh\phi_x(\xi,\tau):=\frac{1}{(T-t(x))^{\frac{1}{\beta+1}}}\wh\Phi\left(\frac{x+\sqrt{T-t(x)}\xi}{\sqrt{(T-t(x))(1-\tau)|\log[(T-t(x))(1-\tau)]|}}\right).
\end{equation}
Since $\wh\Phi$ is radial and by \eqref{e:tx}, the condition \eqref{e:khat0} now means 
$$\wh k(0)=(T-t(x))^{\frac{1}{\beta+1}}\wh\phi_x(0,0)=\wh\Phi\left(\frac{K_0}{4}\right).$$
\end{rem}

\begin{lem} \label{l:Htheta}
For fixed $K_0$, we have the following asymptotics as $x\rightarrow0$:

(i) 
\begin{equation} \label{e:Hsim}
H^*(x)\sim \wh k(1)\theta(x)^{\frac{1}{\beta+1}},
\end{equation}

(ii) 
\begin{equation} \label{e:nabHsim}
\nabla H^*(x)\sim \frac{8}{(\beta+1)K_0}\frac{x}{|x|}\frac{\wh k(1)}{\sqrt{|\log\theta(x)|}}\theta(x)^{\frac{1}{\beta+1}-\frac12}.
\end{equation}
\end{lem}

This lemma can be found in \cite[Lemma 2.2]{MerZaa97NL}. Indeed, as $x\rightarrow0$, 
\begin{equation} \label{e:asymp}
\log \theta(x)\sim2\log|x|\quad\text{and}\quad\theta(x)\sim\frac{8}{K_0^2}\frac{|x|^2}{|\log|x||}.
\end{equation}
Therefore, to prove Lemma \ref{l:Htheta} it suffices to use \eqref{e:H*} and \eqref{e:nablaH*}.

\begin{rem}
In summary, we can perform the following chain of transformations
$$h\mapsto u\mapsto w\mapsto q$$ 
for $(x,t)$ in the entire region $\Omega\times (0,T)$, with $T<1$, 
and the second rescaling transform $h\mapsto k_x$ when $x$ is away from the singularity $x=0$ and the boundary of $\Omega$.
\end{rem}

\section{New shrinking set with refined rates and less constraints} \label{s:ShriSet} 

We aim to construct a solution $h$ to \eqref{e:Quen} which extinguishes at some point $x_0\in\Omega$ and at some finite time $T$.
Observe that the local Cauchy problem for \eqref{e:Quen} is well-posed in the function space $\fH$ that has been introduced in Definition \ref{d:fH}, 
see \cite[Lemma 2.1]{MerZaa97NL} or Appendix \ref{a:C} for details.
For some $t_0\in[0,T)$, we thus define a shrinking subset of $C([t_0,T);\fH)$ intended for the containing space of $h$ that solves \eqref{e:Quen}.
Without loss of generality, we can take $x_0=0\in\Omega\subset\bRN$ and assume $T<e^{-1}$.

Given the analysis on $(\cL+V, T(q))$ and the related Definitions \ref{d:inout}-\ref{d:specde} and \ref{d:R123},
we introduce the following notion which is the key ingredient in our arguments.

\begin{defn}[Shrinking set with refined rates] \label{d:shrink}
Let $T\in(0,e^{-1})$. Fix
$$\underline{\alpha}>3\quad\text{and}\quad \underline{\alpha}+1\leq\overline{\alpha}<\infty.$$ 
Let $K_0>0$, $\ep_0>0$, $A>0$, $\alpha_0>0$, $\delta_0>0$, $C_0>0$ and $\eta_0>0$.

\bigskip 

(I) For all $t_0\in[0,T)$ and for all $t\in[t_0,T)$, 
we define 
$$S^*(t_0,t)=S^*(t_0, K_0, \ep_0, A,\alpha_0,\delta_0, C_0,\eta_0,t)$$ 
as the set of functions $\fh\in \fH$, where $\fH$ is defined in \eqref{e:HRN}-\eqref{e:Hbdd}, satisfying

\bigskip 

(i) Estimates in $\cR_1(K_0,\ep_0,t)$: We require, with $s=-\log(T-t)>1$, that
$$q(\cdot,s)\in V_{K_0,A}(s).$$ 
Here, $q(\cdot,s)=w(\cdot,s)-\varphi(\cdot,s)$ is defined in \eqref{e:qdefn} 
through transforming $h(\cdot,t):=\fh(\cdot)$ first into $u(\cdot,t)$ via \eqref{e:hu} and then into $w(\cdot,s)$ via \eqref{e:utow}, 
and $V_{K_0,A}(s)$ is the set of functions $\fr\in L^\infty(\bRN)$ such that $r(\cdot,s):=\fr(\cdot)$ satisfies
\begin{equation} \label{e:rcomp}
\begin{aligned}
(r_0(s), r_1(s))&\in \left[-\frac{A}{s^2}, \frac{A}{s^2}\right]^{1+N}=:\cQ_A(s), \\
|r_2(s)|&\leq \frac{A^2\log s}{s^2}, \\ 
\left\|\frac{r_-(\cdot,s)}{1+|\cdot|^3}\right\|_{L^\infty}&\leq\frac{A^{\underline{\alpha}}\log s}{s^{\frac52}}, \\ 
\|r_e(\cdot,s)\|_{L^\infty}&\leq\frac{A^{\overline{\alpha}}\log s}{s},
\end{aligned}
\end{equation}
where $r_m(s)$ for $m\in\{0, 1, 2\}$, $r_-(\cdot,s)$ and $r_e(\cdot,s)$ are the components of $r(\cdot,s)$ in the decomposition \eqref{e:rdec}
associated to $s$ and the parameter $K_0$ in defining $\chi_0$ via \eqref{e:chi0}.

\bigskip 

(ii) Estimates in $\cR_2(K_0,\ep_0,t)$: For all 
$$t\in[t_0,T)\quad\text{and}\quad|x|\in \left[\frac{K_0}{4}\sqrt{(T-t)|\log(T-t)|},\ep_0\right],$$
and for all
$$|\xi|\leq\alpha_0\sqrt{|\log \theta (x)|}\quad\text{and}\quad\tau\in\left[\frac{t_0-t(x)}{T-t(x)},1\right),$$
we require
\begin{equation} \label{e:kcomp1}
\left|k_x(\xi,\tau)-\wh k(\tau)\right|\leq \delta_0,
\end{equation}
\begin{equation} \label{e:kcomp2}
|\nabla_\xi k_x(\xi,\tau)|\leq \frac{C_0}{\sqrt{|\log \theta (x)|}}.
\end{equation}
Here, $k_x$ is defined in \eqref{e:kx} via $h(\cdot,t):=\fh(\cdot)$,
$$\tau=\tau(x,t)=\frac{t-t(x)}{\theta(x)},\quad \theta(x)=T-t(x),$$
where $t(x)$ and $\wh k$ are defined respectively in \eqref{e:tx} and \eqref{e:khat}.

\bigskip 

(iii) Estimates in $\cR_3(K_0,\ep_0,t)$: 
For all $|x|\geq\frac{\ep_0}{4}$, we require on $h(\cdot,t):=\fh(\cdot)$ that
\begin{equation} \label{e:hcomp1}
\left|h(x,t)-h(x,t_0)\right|\leq\eta_0,
\end{equation}
\begin{equation} \label{e:hcomp2}
\left|\nabla h(x,t)-\nabla h(x,t_0)\right|\leq\eta_0.
\end{equation}

\bigskip

(II) For all $t_0\in[0,T)$, we define 
$$\begin{aligned}
S^*(t_0)&=S^*(t_0, K_0, \ep_0, A,\alpha_0,\delta_0, C_0,\eta_0)\\
&=\bigg\{h\in C([t_0,T);\fH): \\
&\qquad\forall t\in[t_0,T), h(\cdot,t)\in S^*(t_0, K_0, \ep_0, A,\alpha_0,\delta_0, C_0,\eta_0,t)\bigg\}.
\end{aligned}$$
\end{defn}

The involved parameters for $S^*(t_0, K_0, \ep_0, A,\alpha_0,\delta_0, C_0,\eta_0)$ 
are arranged in the order of their appearance in this definition.
It is clear that the set $V_{K_0,A}(s)$ is shrinking to 0 (respectively, $S^*(t_0,t)$ to the extinction profile) when $s\rightarrow\infty$ (respectively, $t\rightarrow T$).
Recall that our aim is to find a solution $q$ of \eqref{e:qeqn} with $\|q(s)\|_{L^\infty}\rightarrow0$.

\begin{rem} \label{r:defn}
Compared to the shrinking set used in \cite{MerZaa97NL}, 
our definition is new in the following aspects:
1) in Part (i), we allow a flexible dependence in the magnitude parameter $A$, remove the gradient estimate for $q(\cdot,s)$ and refine the shrinking rates in \eqref{e:rcomp} on $q_-$ and $q_e$;
2) we also remove, as in \cite{DuoZaa19}, the second order assumption $|\nabla_{\xi}^2 k_x|\leq C$ in Part (ii).
We keep however the same assumptions in Part (iii), see \cite{DuoZaa19} for a semigroup-adapted version of \eqref{e:hcomp2} when $\Omega$ is bounded.
\end{rem}

\begin{rem}[$L^\infty$-estimates] \label{r:bdd}
Note that \eqref{e:rcomp} from $q(\cdot,s)\in V_{K_0,A}(s)$ implies, for $A\geq1$ and some $C>0$ which is independent of $A$, 
the following desired asymptotics
\begin{equation} 
\|q(\cdot, s)\|_{L^\infty}\leq  C\frac{A^{\overline{\alpha}}\log s}{s}.
\end{equation}
See Lemma \ref{l:qnablaq} for the proof and for related estimates on $q$ via $h(t)\in S^*(t_0,t)$.
The shrinking rate $\frac{\log s}{s}$ ``matches" the $\frac1s$-order decay of the ``remainder" term $R$. 
Accordingly, the decay rates on $q_2$, $q_-$ and $q_e$ are designed under this consideration.
\end{rem}

\begin{rem} \label{r:soboSt}
It is important that one tailors the shrinking set by further imposing gradient estimates on $h$ outside the extinction region
(more precisely, via $\nabla_\xi k_x$ in $\cR_2$ and $\nabla h$ in $\cR_3$), 
see Lemma \ref{l:qnablaq}.
In particular, this allows us to control via $q$ and $\nabla q$ the gradient term $T(q)$ that appears in \eqref{e:qeqn'}, 
see \eqref{e:Tqe} in Lemma \ref{l:Tq}.
\end{rem}

\begin{rem} 
Via the transformations $h\mapsto u \mapsto w$, the constraint $q(s)\in V_{K_0,A}(s)$ in Part (i) describes $h$ mainly in $\cR_1$,
as all the estimates of \eqref{e:rcomp} excluding the one on $q_e$ are stated for the components of $q_b$,
which by definition agrees with $q$ in $\cR_1$. 
\end{rem}

\section{A roadmap to the main theorem} \label{s:Map} 

In this section we outline the main steps to prove Theorem \ref{t:Quen} and state seven propositions that are needed.
Their proofs are left to the seven sections that follow.

\bigskip

Let us consider initial data $h(t_0)$ in the following form: for $(d_0,d_1)\in\bR^{1+N}$, define
\begin{equation} \label{e:ht0}
\begin{aligned}
h(x,t_0&;d_0,d_1)\\
&:=(T-t_0)^{\frac{1}{\beta+1}}\alpha^{\frac{1}{\beta+1}}
\left[\Phi(z)+(d_0+d_1\cdot z)\chi_0\left(\frac{|z|}{K_0/16}\right)\right]^{-\frac{1}{\alpha}}\bigg|_{t=t_0}
\chi_1(x,t_0)\\&
\qquad+H^*(x)\left(1-\chi_1(x,t_0)\right),
\end{aligned}
\end{equation}
where
\begin{equation} \label{e:chi1}
\chi_1(x,t_0)=\chi_0\left(\frac{|x|}{(T-t_0)^{\frac12}|\log(T-t_0)|^{\frac{p}{2}}}\right),
\end{equation}
$H^*(x)$, $(z,\Phi(z))$ and $\chi_0$ are given respectively in Definition \ref{d:H*}, \eqref{e:zPhi} and \eqref{e:chi0}. 
In particular, for $z|_{t=t_0}$ large, the initial data $h(t_0)$ agrees with $H^*(x)$,
while on $\{\chi_1(\cdot,t_0)=1\}$, $h(t_0)$ and its transformation $u(t_0)$ are well-prepared around $\Phi(z)|_{t=t_0}$.

\bigskip

The roadmap which finally leads to the proof of Theorem \ref{t:Quen} consists of six parts.

In Parts I and II, we present the regularity estimates of the solution $h$ that evolves from data as in \eqref{e:ht0}
and is \textbf{assumed to be} trapped in some shrinking set.

The central issue is then to find $(d_0,d_1)$ in some region $\cD$ of $\bR^{1+N}$ such that
the solution $h_{d_0,d_1}=h(\cdot,\cdot,d_0,d_1)$ to \eqref{e:Quen}, 
generated by the initial data $h(\cdot,t_0;d_0,d_1)$ as in \eqref{e:ht0}, is \textbf{indeed} trapped in some shrinking set as given in Definition \ref{d:shrink}, 
namely, for suitable parameters
$$h_{d_0,d_1}\in S^*(t_0, K_0, \ep_0, A,\alpha_0,\delta_0, C_0,\eta_0).$$
Following \cite{BriKup94} and \cite{MerZaa97DMJ}, we deal with this existence issue in two steps:

(i) Finite dimensional reduction: 
Using the \textit{a priori} estimates on $q$, $k_x$ and $h$ in Part III, we find suitable parameters $t_0<T$,
$K_0$, $\ep_0$, $A$, $\alpha_0$, $\delta_0$, $C_0$ and $\eta_0$ so that the following property holds:
if for some $t_*\in[t_0,T)$, we have
$$h(t)\in S^*(t_0, K_0, \ep_0, A,\alpha_0,\delta_0, C_0,\eta_0,t)\,\,\text{for all}\,\,t\in[t_0,t_*]$$
and
 $$h(t_*)\in \pd S^*(t_0, K_0, \ep_0, A,\alpha_0,\delta_0, C_0,\eta_0,t_*),$$
then
$$(q_0(s_*),q_1(s_*))\in \pd\cQ_A(s_*)=\pd\left[-\frac{A}{s_*^2}, \frac{A}{s_*^2}\right]^{1+N},$$
where $s_*=-\log(T-t_*)$. This reduction will be justified in Part IV.

(ii) A topological argument is used to find a pair of parameters $(d^*_0,d^*_1)$ so that 
$$h_{d^*_0,d^*_1}\in S^*(t_0, K_0, \ep_0, A,\alpha_0,\delta_0, C_0,\eta_0).$$
For the sake of completeness, this argument will be formally presented in Part V.

In Part VI, we shall conclude the proof of Theorem \ref{t:Quen} (postponing all the technical arguments for the needed propositions).
This final part is made possible since on one hand, the intermediate/final profiles are consequences of the regularity estimates via the new shrinking set,
and on the other hand, the existence of initial data which generate solutions trapped in the shrinking set is justified in Part V.

\bigskip

\textbf{Part I. Initialization of the Evolution Problem.}

\bigskip

In the following proposition,
we identify a set for the parameters $(d_0,d_1)\in\bR^{1+N}$
so that the initial data $h(\cdot,t_0;d_0,d_1)$ constructed by \eqref{e:ht0} is already in $S^*(t_0,t_0)$.

\begin{prop}[Initialization of the evolution problem] \label{p:init}
There exists $K_{1}>0$ such that for each $K_0\geq K_{1}$ and $\delta_1>0$, there exist $\alpha_1=\alpha_1(K_0,\delta_1)>0$ and $C_0^*=C_0^*(K_0)>0$
such that for all $\alpha_0\leq\alpha_1$, there exists $\ep_1=\ep_1(K_0,\delta_1,\alpha_0)>0$, such that for all $\ep_0\leq \ep_1$, for all $A\geq1$,
there exists $t_1=t_1(K_0,\delta_1,\ep_0,A)<T$ such that
for all $t_0\in[t_1,T)$, there exists a rectangle $\cD(t_0, K_0,A)\subset\bR^{1+N}$ with the following properties:

\bigskip

If $h(t_0)=h(\cdot,t_0;d_0,d_1)$ is given by \eqref{e:ht0}, then

\bigskip 

(i) For all $(d_0,d_1)\in\cD(t_0, K_0,A)$, we have
$$h(t_0)\in S^*(t_0, K_0, \ep_0, A,2\alpha_0,\delta_1, C_0^*,0,t_0).$$ 
More precisely, 

-- (i0) we have
$$h(t_0)\in \fH.$$

-- (i1) by setting $s_0=-\log(T-t_0)$, then for $q$ defined via \eqref{e:qdefn} we have
$$(q_0(s_0), q_1(s_0))\in \cQ_A(s_0),$$ 
$$|q_2(s_0)|\leq \frac{\log s_0}{s_0^2},$$
$$\left\|\frac{q_-(\cdot,s_0)}{1+|\cdot|^3}\right\|_{L^\infty}\leq\frac{\log s_0}{s_0^{\frac52}},$$ 
$$\|q_e(\cdot,s_0)\|_{L^\infty}\leq\frac{\log s_0}{s_0}.$$
In addition,
$$\left\|\frac{(\nabla q)_\perp(\cdot,s_0)}{1+|\cdot|^3}\right\|_{L^\infty}\leq\frac{\log s_0}{s_0^{\frac52}}.$$ 

-- (i2) by setting 
$$\tau_0=\frac{t_0-t(x)}{\theta(x)},\quad \theta(x)=T-t(x),$$
where $t(x)$ is defined in  \eqref{e:tx},
then for $k_x(\xi,\tau)$ with $|\xi|\leq 2\alpha_0\sqrt{|\log \theta(x)|}$, we have
$$|k_x(\xi,\tau_0)-\wh k(\tau_0)|\leq \delta_1,$$ 
$$|\nabla_\xi k_x(\xi,\tau_0)|\leq \frac{C_0^*}{\sqrt{|\log \theta (x)|}},$$
where $k$ and $\wh k$ are defined respectively in \eqref{e:kx} and \eqref{e:khat}.

\bigskip

(ii) We have the following correspondences
$$(d_0,d_1)\in\cD(t_0, K_0,A)\Leftrightarrow (q_0(s_0),q_1(s_0))\in\cQ_A(s_0),$$
$$(d_0,d_1)\in\pd\cD(t_0, K_0,A)\Leftrightarrow (q_0(s_0),q_1(s_0))\in \pd\cQ_A(s_0),$$
and $(q_0(s_0),q_1(s_0))$ is an affine function of $(d_0,d_1)$ when $(d_0,d_1)\in\pd\cD(t_0, K_0,A)$.
\end{prop}

\begin{proof}
See Section \ref{s:init}.
\end{proof}

Note that the last estimate in $(i1)$ is not contained in $h(t_0)\in S^*(t_0,t_0)$.

\bigskip

\textbf{Part II. Parabolic Regularity under the Partial Trapping Assumption.}

\bigskip

Under the partial trapping assumption, the following gradient estimates are valid for the solution $q$ to \eqref{e:qeqn} that is generated by $q(s_0)$, $s_0=-\log(T-t_0)$.
Here $q(s_0)$ is obtained, via \eqref{e:hu}, \eqref{e:utow} and \eqref{e:qdefn}, from initial data $h(x,t_0;d_0,d_1)$ in \eqref{e:ht0}.

\begin{prop}[Parabolic regularity via partial trapping in the shrinking set] \label{p:Sobo} 
There exists $K_{1.5}>0$ such that for each $K_0\geq K_{1.5}$ and $\delta_{1.5}>0$, 
there exist $\alpha_{1.5}=\alpha_{1.5}(K_0,\delta_{1.5})>0$ and $C_0=C_0(K_0)>0$
such that for all $\alpha_0\leq\alpha_{1.5}$, there exists $\ep_{1.5}=\ep_{1.5}(K_0,\delta_{1.5},\alpha_0)>0$, such that for all $\ep_0\leq \ep_{1.5}$, for all $A\geq1$,
there exists $t_{1.5}=t_{1.5}(K_0,\delta_{1.5},\ep_0,A)<T$ and $\eta_{1.5}=\eta_{1.5}(\ep_0)$ 
such that for all $t_0\geq t_{1.5}$, $\delta_0\leq \delta_{1.5}$ and $\eta_0\leq\eta_{1.5}$,
the following property holds: 

\bigskip

Assume that $h$ is the solution to \eqref{e:Quen}
generated by the initial data $h(t_0;d_0,d_1)$ defined in \eqref{e:ht0} where $(d_0,d_1)$ is chosen so that $(q_0(s_0), q_1(s_0))\in \cQ_A(s_0)$ with $s_0=-\log(T-t_0)$.
Assume further that $h$ satisfies, for some $t_*\in[t_0,T)$,
$$\forall \,\, t\in[t_0,t_*],\quad h(t)\in S^*(t_0, K_0, \ep_0, A,\alpha_0,\delta_0, C_0,\eta_0,t).$$
Then, for some $C=C(K_0,C_0)>1$ and for $s\in[s_0,-\log(T-t_*)]$,
\begin{equation} \label{e:parreg1}
\|\nabla q(s)\|_{L^\infty}\leq  C\frac{A^{\overline{\alpha}}\log s}{s}
\end{equation}
and
\begin{equation} \label{e:parreg2'}
|(\nabla q)(y,s)|\leq C\frac{A^2\log s}{s^{2}}(1+|y|^2)+C\frac{A^{\overline{\alpha}}\log s}{s^{\frac52}}(1+|y|^3),
\end{equation}
where $s=-\log(T-t)$ and $q$ is in \eqref{e:qdefn}. 
In particular, for $s\in[s_0,-\log(T-t_*)]$,
\begin{equation} \label{e:parreg2}
\left\|\frac{(\nabla q)(\cdot,s)}{1+|\cdot|^3}\right\|_{L^\infty}\leq  C\frac{A^{\overline{\alpha}}\log s}{s^2}.
\end{equation}
\end{prop}

\begin{proof}
See Section \ref{s:Sobo}.
\end{proof}

\bigskip

\textbf{Remark.} This proposition only assumes partial trapping in $S^*(t_0,t)$ ($t\leq t_*<T$) 
and is weaker than the full trapping to be addressed later in Part IV.
It however serves for our purpose (to be used in conjunction with the propositions in next part).
This proposition is a consequence of Proposition \ref{p:init} and parabolic regularity. 

\bigskip

As consequences of \eqref{e:parreg1}, \eqref{e:parreg2'}, \eqref{e:qbEst'}, \eqref{e:Tqb} and \eqref{e:Tqe'}, we have
\begin{equation} \label{e:Tq*bd}
|T(q)(y,s)|\leq \frac{C}{s}
\end{equation}
and
\begin{equation} \label{e:Tq*}
|T(q)(y,s)|\leq C\frac{A^{2} \log s}{s^{\frac52}}(1+|y|^3)
\end{equation}
for $s\in[s_0,-\log(T-t_*)]$ where $s_0=s_0(A)$ is large enough. 

\bigskip

\textbf{Part III. \textit{A Priori} Estimates in $\cR_1$, $\cR_2$ and $\cR_3$.}

\bigskip

In this technical part we give \textit{a priori} estimates for the solutions $q$, $k_x$ and $h$.

\bigskip

We first present some self-improving properties in the non-expanding directions for $q(s)\in V_{K_0,A}(s)$.
This concerns the components $q_2$, $q_-$ and $q_e$ of $q$.

\begin{prop}[\textit{A priori} estimates in $\cR_1$] \label{p:impr1}
There exists $K_{2}>0$ such that for each $K_0\geq K_{2}$, there exists $A_2=A_2(K_0)>1$ such that for each $A\geq A_2$, $\ep_0>0$ and $C_0\leq A^3$,
there exists $\eta_2=\eta_2(\ep_0)>0$ and $t_2=t_2(K_0, \ep_0, A, C_0)\in(0, T)$ such that for each $t_0\in[t_2,T)$, $\delta_0\leq \frac{\wh k(1)}{2}$, $\alpha_0>0$
and $\eta_0\leq\eta_2$, we have the following properties:

\bigskip

If $h$ is the solution to \eqref{e:Quen}, 
generated by the initial data $h(t_0;d_0,d_1)$ defined in \eqref{e:ht0} where $(d_0,d_1)$ is chosen so that $(q_0(s_0), q_1(s_0))\in \cQ_A(s_0)$ with $s_0=-\log(T-t_0)$,
and satisfies for some $t_*\in[t_0,T)$ that
$$\forall \,\, t\in[t_0,t_*],\quad h(t)\in S^*(t_0, K_0, \ep_0, A,\alpha_0,\delta_0, C_0,\eta_0,t),$$
then for $q$ the solution of \eqref{e:qeqn} which is defined from $h$ via \eqref{e:qdefn}, we have
$$|q_2(s_*)|\leq \frac{A^2\log s_*}{s_*^2}-\frac{1}{s_*^3},$$ 
$$\left\|\frac{q_-(\cdot,s_*)}{1+|\cdot|^3}\right\|_{L^\infty}\leq\frac{A^{\underline{\alpha}}\log s_*}{2s_*^{\frac52}},$$\
$$\|q_e(\cdot,s_*)\|_{L^\infty}\leq\frac{A^{\overline{\alpha}}\log s_*}{2s_*},$$
where $s_*=-\log(T-t_*)$ and the notations $q_2$, $q_-$ and $q_e$ are given by \eqref{e:rdec}.
\end{prop}

\begin{proof}
See Section \ref{s:impr1}.
\end{proof}

\bigskip

The following result is stated for general functions $k=k(\xi,\tau)$ which are not necessarily being $k_x$ as transformed from $h_{|\cR_2}$ via \eqref{e:kx}.
In particular, it is independent of the shrinking set although the involved estimates are modelled on those satisfied by $k_x$.
Note that in this result, the $K_0$-dependence comes from $\wh k$ as defined in \eqref{e:khatsol}.

\begin{prop}[\textit{A priori} estimates in $\cR_2$, or rather, on $k$] \label{p:impr2} 
There exists $K_{3}>0$ such that for each $K_0\geq K_{3}$, $\delta_3\leq1$, $C_0'^*>0$, $C''^*_0>0$, 
we have the following property:

\bigskip

Assume that $k=k(\xi,\tau)$ is a solution to the nonlinear heat equation
\begin{equation} \label{e:knox}
\pd_\tau k=\Delta_\xi k-N_\beta(k)
\end{equation}
for $\tau\in[\tau_1,\tau_2)$ with $0\leq \tau_1\leq\tau_2\leq1$,
where $N_\beta$ is given in \eqref{e:keqn}.

Assume in addition that for all $|\xi_0|\geq1$, there holds: 
$\forall\,\xi\in B\left(0,\frac{7|\xi_0|}{4}\right)$, $\forall\,\tau\in[\tau_1,\tau_2)$,
\begin{equation} \label{e:R2-1}
k(\xi,\tau)\geq\frac{\wh k(\tau)}{2} \quad\text{and}\quad |\nabla k(\xi,\tau)|\leq\frac{C_0'^*}{|\xi_0|},
\end{equation}
and $\forall\,\xi\in B\left(0,2|\xi_0|\right)$, 
$$\left|k(\xi,\tau_1)-\wh k(\tau_1)\right|\leq\delta_3\quad\text{and}\quad|\nabla k(\xi,\tau_1)|\leq \frac{C_0''^*}{|\xi_0|}.$$

Then there exists $\xi_3=\xi_3(C_0'^*, C''^*_0)$ such that for all $\xi_0$ with $|\xi_0|\geq|\xi_3|$, there exists $\ep=\ep(K_0,C_0'^*,\delta_3,\xi_0)$
such that $\forall\,\xi\in B(0,|\xi_0|)$, $\forall\,\tau\in[\tau_1,\tau_2)$,
$$\left|k(\xi,\tau)-\wh k(\tau)\right|\leq\ep \quad\text{and}\quad|\nabla k(\xi,\tau)|\leq\frac{C_{3}C_0''^*}{|\xi_0|},$$
where $\ep\rightarrow0$ as $(\delta_3,|\xi_0|)\rightarrow(0,+\infty)$ and $C_{3}>1$ is a universal constant. 
\end{prop}

\begin{proof}
See Section \ref{s:impr2}.
\end{proof}

\bigskip

Finally, we also need a regularity result which is directly stated on $h$.

\begin{prop}[\textit{A priori} estimates in $\cR_3$] \label{p:impr3}  
For all $\ep>0$, $\ep_0>0$, $\sigma_0>0$ and $\sigma_1>0$,
there exists $t_4=t_4(\ep,\ep_0,\sigma_0,\sigma_1)<T$ such that we have the following property:

\bigskip

For all $t_0\in[t_4,T)$, assume that $h$ is a solution of \eqref{e:Quen} on $[t_0,t_*]$ for some $t_*\in[t_0,T)$ satisfying

\indent\indent (i) for $|x|\in[\frac{\ep_0}{6},\frac{\ep_0}{4}]$, $t\in[t_0,t_*]$,
\begin{equation} \label{e:R3-1}
\sigma_0\leq h(x,t)\leq\sigma_1\quad\text{and}\quad |\nabla h(x,t)|\leq \sigma_1,
\end{equation}

\indent\indent(ii) for $|x|\geq\frac{\ep_0}{6}$, 
\begin{equation} \label{e:R3-2}
h(x,t_0)=H^*(x),
\end{equation}
where $H^*$ is defined in Definition \ref{d:H*}. 

Then for $|x|\geq\frac{\ep_0}{4}$, $t\in[t_0,t_*]$,
\begin{equation} \label{e:hR3}
|h(x,t)-h(x,t_0)|+|\nabla h(x,t)-\nabla h(x,t_0)|\leq\ep.
\end{equation}
\end{prop}

\begin{proof}
See Section \ref{s:impr3}.
\end{proof}

Unlike Proposition \ref{p:impr1}, Propositions \ref{p:impr2}-\ref{p:impr3} are independent of the shrinking set.
We point out however that \textit{a priori} estimates for $k_x$ in $\cR_2$ via the shrinking set, in contrast with Proposition \ref{p:impr2}, 
will be considered in Lemma \ref{l:rigidk} of Section \ref{s:finredu}.

\bigskip

\textbf{Part IV. Finite Dimensional Reduction.} 

\bigskip

The following proposition shows that, roughly speaking, 
the leading size information of $S^*(t_0,t)$ is encoded in the $(1+N)$-dimensional vector $(q_0(s),q_1(s))$.

\begin{prop}[Finite dimensional reduction] \label{p:finredu}
We can choose in order parameters $K_0$, $\delta_0$, $A$, $C_0$, $\alpha_0$, $\ep_0$, $\eta_0$ and $t_0<T$
such that the following properties hold:

\bigskip

Assume that $h(t_0)$ is given by \eqref{e:ht0} with $(d_0,d_1)\in\cD(t_0,K_0,A)$. Then

(i) (Initialization)
$$h(t_0)\in S^*(t_0, K_0, \ep_0, A,\alpha_0,\delta_0, C_0,\eta_0,t_0).$$
Assume in addition that for some $t_*\in[t_0,T)$, we have for all $t\in[t_0,t_*]$
$$h(t)\in S^*(t_0, K_0, \ep_0, A,\alpha_0,\delta_0, C_0,\eta_0,t)$$
and
 $$h(t_*)\in \pd S^*(t_0, K_0, \ep_0, A,\alpha_0,\delta_0, C_0,\eta_0,t_*),$$
where $h(t)$ is the solution to \eqref{e:Quen} generated with initial data $h(t_0)$.
Then

(ii) (Finite dimensional reduction) there holds
$$(q_0(s_*),q_1(s_*))\in \pd\cQ_A(s_*),$$
where $s_*=-\log(T-t_*)$.

(iii) (Transversality) there exists $\nu_0>0$ such that for all $\nu\in(0,\nu_0)$
$$(q_0(s_*+\nu),q_1(s_*+\nu))\notin \cQ_A(s_*+\nu),$$
hence there exists $\wt \nu_0>0$ such that for all $\wt\nu\in(0,\wt\nu_0)$
$$h(t_*+\wt\nu)\notin S^*(t_0, K_0, \ep_0, A,\alpha_0,\delta_0, C_0,\eta_0,t_*+\wt\nu).$$
\end{prop}

\begin{proof}
See Section \ref{s:finredu}.
\end{proof}

In the above statement we specified ``in order", which is the essence in the proof.

\bigskip

\textbf{Part V. Contradiction via Topological Arguments.}

\bigskip

From Proposition \ref{p:finredu}, we claim that there exists $(d_0,d_1)\in \cD(t_0,K_0,A)$ such that
$$h_{d_0,d_1}=h(\cdot,\cdot\, ;d_0,d_1)\in S^*(t_0, K_0, \ep_0, A,\alpha_0,\delta_0, C_0,\eta_0),$$
where $h(\cdot,\cdot\, ;d_0,d_1)$ is the solution to \eqref{e:Quen} generated by $h(\cdot,t_0;d_0,d_1)$ as in \eqref{e:ht0}.

\bigskip

\textbf{Remark.} For later use, we denote this pair of parameters by $(d^*_0,d^*_1)$.

\bigskip

We proceed by contradiction 
and assume that for all $(d_0,d_1)\in \cD(t_0,K_0,A)$,
$$h(\cdot,\cdot\, ;d_0,d_1)\notin S^*(t_0, K_0, \ep_0, A,\alpha_0,\delta_0, C_0,\eta_0).$$
From (i) of Proposition \ref{p:finredu},
$$h(\cdot,t_0\, ;d_0,d_1)\in S^*(t_0, K_0, \ep_0, A,\alpha_0,\delta_0, C_0,\eta_0,t_0).$$
Therefore, for each $(d_0,d_1)\in \cD(t_0,K_0,A)$ we denote by $t_*(d_0,d_1)<T$ the infimum of all $t\in[t_0,T)$ such that
$$h(\cdot,t\, ;d_0,d_1)\notin S^*(t_0, K_0, \ep_0, A,\alpha_0,\delta_0, C_0,\eta_0,t).$$
By (ii) of Proposition \ref{p:finredu}, for $(q_0,q_1)$ obtained from $h_{d_0,d_1}$, we have
$$(q_0,q_1)(s_*(d_0,d_1))\in \pd\cQ_A(s_*(d_0,d_1)),$$
where $s_*(d_0,d_1)=-\log(T-t_*(d_0,d_1))$.
Hence we can define the following function
\begin{equation} \label{e:Psi}
\begin{aligned}\Psi:  \,\,\cD(t_0,K_0,A)&\rightarrow\pd \cQ^\sharp, \\
(d_0,d_1)&\mapsto \frac{[s_*(d_0,d_1)]^2}{A}(q_0,q_1)(s_*(d_0,d_1)),
\end{aligned}
\end{equation}
where $\cQ^\sharp=[-1,1]^{1+N}$ is the unit cube.
Then, for $\Psi$ defined in \eqref{e:Psi}, we claim that

(i) $\Psi$ is a continuous mapping from $\cD(t_0,K_0,A)$ to $\pd\cQ^\sharp$.

(ii) There exists a non-trivial affine function $G: \cD(t_0,K_0,A)\rightarrow \cQ^\sharp$ such that
$$\Psi\circ G^{-1}|_{\pd\cQ^\sharp}=Id_{\pd\cQ^\sharp}.$$
Part (i) of the claim follows from the continuity in $\fH$ of the solution $h(\cdot,t\, ;d_0,d_1)$ at a fixed time $t$ with respect to the initial data 
(thus $(q_0,q_1)$ is continuous with respect to $(d_0,d_1)$),
and the continuity of $s_*(d_0,d_1)$ with respect to $(d_0,d_1)$, 
which is a direct consequence of the transverse property from (iii) of Proposition \ref{p:finredu}. 
From (ii) of Proposition \ref{p:init}, we have
$$\forall\,\,(d_0,d_1)\in \pd\cD(t_0,K_0,A),\quad s_*(d_0,d_1)=s_0.$$
Then Part (ii) of the claim follows. 
Indeed, we can take $G=\frac{s_0^2}{A}g$ with $$g: (d_0,d_1)\mapsto (q_0(s_0),q_1(s_0))$$ being the affine function constructed in Section \ref{s:init}.
Now, we deduce that $\Psi$ has a non-zero degree on the boundary,
which is a contradiction (index theory). 
This proves the existence of $(d_0,d_1)\in \cD(t_0,K_0,A)$ such that $h(\cdot,\cdot\, ;d_0,d_1)\in S^*(t_0).$

\bigskip

\textbf{Part VI. Full Trapping Implies the Gradient Profile.}

\bigskip

Theorem \ref{t:Quen} is a consequence of the following proposition.

\begin{prop}[Derivation of the gradient profile via refined shrinking rates] \label{p:profile}
Consider $t_0<T$, $K_0$, $\ep_0$, $A$, $\alpha_0$, $\delta_0$, $C_0$ and $\eta_0$ 
such that 
\begin{equation} \label{e:del0eta0}
\delta_0\leq \frac{1}{2}\wh k(1)\quad\text{and}\quad \eta_0\leq \frac12\inf_{|x|\geq\frac{\ep_0}{4}}h(x,t_0),
\end{equation}
where $h(x,t_0)$ is in the form of \eqref{e:ht0} with $(d_0,d_1)=(d^*_0,d^*_1)$ as determined in Part V. 
Then $h(t)$ extinguishes in finite time $T$ only at the point $x=0$,
that is, 
\begin{equation} \label{e:0extinc}
\lim_{t\rightarrow T}h(0,t)=0 
\end{equation}
and $\forall \,x\neq0$, there exists $\omega(x)>0$ such that 
\begin{equation} \label{e:0extinc!}
\liminf_{t\rightarrow T}\inf_{|x'-x|\leq\omega(x)}h(x',t)>0.
\end{equation}
Moreover, the refined estimates \eqref{e:refined1}-\eqref{e:refined2} and the final profiles \eqref{e:profsame}-\eqref{e:profgrad} hold for the solution and its gradient.
\end{prop}

\begin{proof}
See Section \ref{s:profile}.
\end{proof}

\begin{rem}
It is harmless to assume $T<e^{-1}$ small enough so that 
$$-\log T\geq s_0.$$
This is reasonable since $T$ depends on the initial data.
In other words, the asymptotic estimates \eqref{e:refined1}-\eqref{e:refined2} hold for all $t\geq0$ instead of $t\geq t_0=T-e^{-s_0}$ 
with $t_0$ the initial time of trapping that is determined in Proposition \ref{p:finredu}.
Another option is to rewrite the shrinking set definition as $T$-dependent,
and instead of finding $s_0$- or $t_0$-threshold, one determines finally the threshold on $T$.
This is adopted for example in \cite{Duo19}.
\end{rem}

\section{Initialization of the evolution problem} \label{s:init}

In this section, we prove Proposition \ref{p:init}.

\bigskip

\textbf{Proof of $(i0)$}.
Observe that for some $C(t_0,d_0,d_1)>0$,
$$\inf_{x\in\Omega} h(x,t_0)\geq C(t_0,d_0,d_1).$$
By Definitions \ref{d:H*} and \ref{d:fH}, together with \eqref{e:Phi}, one verifies that 
$$h(t_0)\in\fH\quad ( \text{more precisely,\,\,}h(t_0)\in\fH_\psi \text{\,\,if\,\,} \Omega=\bRn),$$ 
where $$\psi(x)=H^*(x)\left(1-\chi_0(|x|)\right)(1-\chi_1(x,t_0)).$$

\bigskip

\textbf{Proof of $(i2)$}.
Since in this part we require in defining the shrinking set almost the same conditions as in \cite{MerZaa97NL}, 
the arguments are completely the same as the corresponding ones on \cite[Step 1, pages 1531-1535]{MerZaa97NL}.

\bigskip

\textbf{Proof of $(i1)$ and $(ii)$}.
This is also implicit in \cite{MerZaa97NL}. 
By \cite[(109)]{MerZaa97NL},
$$q_0(s_0)=d_0\int_{\bRN}\chi_0\left(\frac{|y|}{\sqrt{s_0}K_0/16}\right)d\rho(y)
-\frac{\kappa}{2(p-a)s_0}\int_{\bRN}\chi_0\left(\frac{|y|}{\sqrt{s_0}K_0}\right)d\rho(y),$$
$$q_1(s_0)=\frac{d_1}{\sqrt{s_0}}\int_{\bRN}\frac{y^2}{2}\chi_0\left(\frac{|y|}{\sqrt{s_0}K_0/16}\right)d\rho(y).$$
The map $g: (d_0,d_1)\mapsto (q_0,q_1)(s_0)$ is an affine function. Introduce
$$\cD(t_0, K_0,A)=g^{-1}(\cQ_A(s_0)),$$
which is then a cube in $\bR^{1+N}$.
If $(d_0,d_1)\in \cD(t_0, K_0,A)$, or equivalently 
$$|q_m(s_0)|\leq\frac{A}{s_0^2}\quad\text{for}\quad m\in\{0,1\},$$
then as shown in \cite{MerZaa97NL}
$$|d_0|\leq \frac{C}{s_0}\quad\text{and}\quad|d_1|\leq C\frac{A}{s_0^{\frac32}}.$$
However, as proved on \cite[pages 1535-1536]{MerZaa97NL}, for $K_0\geq20$ and $t_0\geq t_{1'}$, 
$$|q_-(y,s_0)|\leq\left(|d_0|+|d_1|+\frac{C}{s_0}\right)\frac{(1+|y|^3)}{s_0^{\frac32}},\quad |q_e(y,s_0)|\leq\frac{C}{s_0},$$
$$\left|\left(\nabla q\right)_\perp(y,s_0)\right|\leq C\frac{(|d_0|+|d_1|+1/s_0)}{\sqrt{s_0}}\frac{(1+|y|^3)}{s_0^{\frac32}}.$$
Here $C=C(K_0)>0$ and $$(\nabla q)_\perp=\left[\left(\frac{\pd q}{\pd y_i}\right)_\perp\right]_{i=1}^N,$$
where $(\cdot)_\perp$ is given in \eqref{e:rdecPerp}. 
Now using the new conditions required in Part (i) of the shrinking set,
for $t_0\geq t_{1''}(K_0,A,C)$ and $s_0=-\log(T-t_0)$ we obtain
$$|q_2(s_0)|\leq \frac{\log s_0}{s_0^2},$$
$$\left\|\frac{q_-(\cdot,s_0)}{1+|\cdot|^3}\right\|_{L^\infty}\leq\frac{\log s_0}{s_0^{\frac52}},\quad \|q_e(\cdot,s_0)\|_{L^\infty}\leq\frac{\log s_0}{s_0},$$
\begin{equation} \label{e:nabqperp}
\left\|\frac{(\nabla q)_\perp(\cdot,s_0)}{1+|\cdot|^3}\right\|_{L^\infty}\leq\frac{\log s_0}{s_0^{\frac72}}\leq\frac{\log s_0}{s_0^{\frac52}}.
\end{equation}
Here, the estimate on $q_2$ remains the same as in \cite{MerZaa97NL}.

The fact that (ii) is true follows from the construction of $g$ and $\cD(t_0, K_0,A)$.

\begin{rem} \label{r:nablaq0}
On \cite[Page 1536]{MerZaa97NL} it was also observed  that
\begin{equation} \label{e:nabqe}
\left|(\nabla q)_e(y,s_0)\right|\leq \frac{1}{s_0}.
\end{equation}
Via integration by parts 
(see for example \cite[Page 1545]{MerZaa97NL}), 
we have for $m=0,1$,
\begin{equation} \label{e:01by12}
(\nabla q)_m(s)=O(e^{-s})+(m+1)q_{m+1}(s). 
\end{equation}
Moreover, according to the first inequality in \eqref{e:nabqperp}, 
$$\left|(\nabla q)_\perp(y,s_0)\right|\leq\frac{\log s_0}{s_0^{\frac72}}(1+|y|^3)\leq C(K_0)\frac{\log s_0}{s_0^2}.$$
By \eqref{e:nabqe}, \eqref{e:01by12} and the bounds on $q_1(s_0)$ and $q_2(s_0)$, this implies that for $s_0$ large
$$\left\|(\nabla q)(\cdot,s_0)\right\|_{L^\infty}\leq \frac{C}{s_0}.$$
\end{rem}

\section{Parabolic regularity via partial trapping in the shrinking set} \label{s:Sobo} 

In this section we prove Proposition \ref{p:Sobo}.

\bigskip

We shall use the parabolic regularity techniques of \cite{EbdZaa11} which is later developed in \cite[Proposition 4.17]{TayZaa19}.
Namely, we use the gradient regularity of the ``semigroup'' $e^{t\cL}$ involved in the standard Duhamel integral representation of \eqref{e:qeqn}:
$$q(s)=e^{(s-\underline s)\cL}q(\underline s)+\int_{\underline s}^se^{(s-s')\cL}\,F(s')ds',$$
$$F=Vq+B(q)+T(q)+R+L(q),$$
and we take $\underline s=s_0$ if $s\leq s_0+1$ and $\underline s=s-1$ if $s> s_0+1$.
In particular, for $\theta>0$, 
$$\|\nabla(e^{\theta\cL}r)\|_{L^\infty}\leq Ce^{\frac{\theta}{2}}\|\nabla r\|_{L^\infty}$$
and
$$\|\nabla(e^{\theta\cL}r)\|_{L^\infty}\leq C\frac{e^{\frac{\theta}{2}}}{(1-e^{-\theta})^{\frac{1}{2}}}\|r\|_{L^\infty}.$$ 
These mapping properties follow from the integral representation
$$e^{\theta\cL}r(y)=\int_{\bRN}e^{\theta\cL}(y,x)r(x)dx$$
and Mehler's formula
$$e^{\theta\cL}(y,x)=e^{\theta}e^{\theta\bL}(y,x)=\frac{e^\theta}{(4\pi(1-e^{-\theta}))^\frac{N}{2}}\exp\left[-\frac{|ye^{-\frac{\theta}{2}}-x|^2}{4(1-e^{-\theta})}\right].$$
For other estimates concerning $e^{\theta \cL}$ and its perturbed version $e^{\theta \cL_V}$, see \cite[Appendix C]{NguZaa17}.
We shall also use Gr\"onwall's argument as in \cite{TayZaa19}.
Below we only indicate the main estimates that are needed in adapting \cite[Proposition 4.17]{TayZaa19}.

\bigskip

\textbf{Proof of \eqref{e:parreg1}}. By \eqref{e:qest}, \eqref{e:V1}, \eqref{e:Bq'}, \eqref{e:Tqe}, \eqref{e:R1} and \eqref{e:Lq}, we have
$$|\left(Vq+B(q)+(1-\chi)T(q)+R+L(q)\right)(s)|\leq C\frac{A^{\overline{\alpha}}\log s}{s},$$
and by \eqref{e:Tqb},
$$\|F(s)\|_{L^\infty}\leq C\frac{A^{\overline{\alpha}}\log s}{s}+C\left(\frac{\|\nabla q(s)\|_{L^\infty}}{\sqrt s}+\|\nabla q(s)\|_{L^\infty}^2\right).$$
The remaining arguments for \eqref{e:parreg1} are similar to Part (i) of \cite[Proposition 4.17]{TayZaa19}. 

\bigskip

\textbf{Proof of \eqref{e:parreg2'}}. Let us define the function
$$H(y,s)=\frac{A^2\log s}{s^{2}}(1+|y|^2)+\frac{A^{\overline{\alpha}}\log s}{s^{\frac52}}(1+|y|^3).$$
By \eqref{e:qbEst'} and \eqref{e:qeEst'}, for $s$ large enough it holds that
$$|q(y,s)|\leq CH(y,s).$$ 
By \eqref{e:V1}, \eqref{e:Bq'}, \eqref{e:Tqe'}, Lemma \ref{l:R} on $R$, and \eqref{e:Lq}, we have
$$|Vq+B(q)+(1-\chi)T(q)+R+L(q)|\leq CH.$$
Using \eqref{e:Tqb} 
and \eqref{e:parreg1} proved in Part (i), 
we further have
$$\begin{aligned}
|F(y,s)|&\leq CH(y,s)+C\left(\frac{|\nabla q|}{\sqrt s}+|\nabla q|^2\right)\\
&\leq CH(y,s)+C\left(\frac{1}{\sqrt s}+\frac{A^{\overline{\alpha}}\log s}{s}\right)|\nabla q|(y,s).
\end{aligned}$$
Note that by the proof of Proposition \ref{p:init} in Section \ref{s:init} (see also Remark \ref{r:nablaq0}), 
one infers indeed the following gradient estimate
$$|\nabla q(y,s_0)|\leq C\frac{\log s_0}{s_0^{\frac52}}(1+|y|^3).$$
For \eqref{e:parreg2'} it remains to adapt the proof of Part (ii) of \cite[Proposition 4.17]{TayZaa19}.

\bigskip

The proof of Proposition \ref{p:Sobo} is complete.

\section{A priori estimates in $\cR_1$} \label{s:impr1}

In this section we prove Proposition \ref{p:impr1}.

\bigskip

First of all we prove the following result which is the heart of our arguments.
It concerns the ``mode dynamics" of \eqref{e:qeqn} on $q$ through the solution $h$ to \eqref{e:Quen} that is trapped in the shrinking set.
The statement on $q_0$ and $q_1$ will not be used in this section, 
but we include it for later purpose (in the proof of Proposition \ref{p:finredu}).

\begin{prop}[Mode dynamics] \label{p:qdyn}
There exist $K_{5}>1$ and $A_5>1$ such that for each $K_0\geq K_{5}$, $A\geq A_5$, $\ell^*>0$, 
there exists $t_5=t_5(K_0, A,\ell^*)< T$ such that:

\bigskip

For all $t_0\in[t_5,T)$, $\ell\in[0,\ell^*]$, $\ep_0>0$, $\alpha_0>0$, $\delta_0\leq\frac{\wh k(1)}{2}$, $C_0>0$ and $\eta_0\leq \eta_5(\ep_0)$ for some $\eta_5(\ep_0)>0$, 
assume that

\indent $\bullet$ $h(t_0)=h(\cdot,t_0;d_0,d_1)$ is given by \eqref{e:ht0} with $(d_0,d_1)$ chosen so that 
$$(q_0(s_0),q_1(s_0))\in \cQ_A(s_0)$$
where $s_0=-\log(T-t_0)$.

\indent $\bullet$ for some $\sigma\geq s_0$ and $\forall \, t\in[T-e^{-\sigma},T-e^{-(\sigma+\ell)}]$, we have 
$$h(t)\in S^*(t_0, K_0, \ep_0, A,\alpha_0,\delta_0, C_0,\eta_0,t),$$
where $h$ is the solution of \eqref{e:Quen} generated by $h(t_0)$.

Then, for $q$ defined in \eqref{e:qdefn} (via $h$) and all $s=-\log(T-t)\in[\sigma, \sigma+\ell]$,
we have

(i) (Unstable modes $q_0$ and $q_1$): for some $\overline C>0$ and $m=0,1$,

\begin{equation} \label{e:01mod}
\left|q_m'(s)-\left(1-\frac{m}{2}\right)q_m(s)\right|\leq \frac{\overline C}{s^2}.
\end{equation}

(ii) (Neutral mode $q_2$): for some $ C>0$, 

\begin{equation} \label{e:2mod}
\left|q_2'(s)+\frac2s q_2(s)\right|\leq \frac{CA}{s^3}. 
\end{equation}

(iii) (Control of the negative spectrum part $q_-$ and the exterior part $q_e$):

\indent\indent (iii-1) If $\sigma>s_0$, then

\begin{equation} \label{e:q-1}
\begin{aligned}
\left\|\frac{q_-(\cdot,s)}{1+|\cdot|^3}\right\|_{L^\infty}&\leq\frac{C\log s}{s^{\frac52}}\\&\quad\times
\left({A^{\underline{\alpha}}}{e^{-\frac{s-\sigma}{2}}}+{A^{\overline{\alpha}}}{e^{-(s-\sigma)^2}}+A^2e^{s-\sigma}((s-\sigma)^2+1)+A^2(s-\sigma)\right),
\end{aligned}
\end{equation}

\begin{equation} \label{e:qe1}
\|q_e(\cdot,s)\|_{L^\infty}\leq\frac{C\log s}{s}\left({A^{\overline{\alpha}}}{e^{-\frac{s-\sigma}{p}}}+A^{\underline{\alpha}}e^{s-\sigma}+A^2(s-\sigma)e^{s-\sigma}\right).
\end{equation}

\indent\indent (iii-2) If $\sigma=s_0$, then

\begin{equation} \label{e:q-2}
\left\|\frac{q_-(\cdot,s)}{1+|\cdot|^3}\right\|_{L^\infty}\leq\frac{C\log s}{s^{\frac52}}\left(1+(s-s_0)^3e^{s-s_0}\right),
\end{equation}

\begin{equation} \label{e:qe2}
\|q_e(\cdot,s)\|_{L^\infty}\leq\frac{C\log s}{s}\left(1+(s-s_0)e^{s-s_0}\right).
\end{equation}
\end{prop}

\begin{proof}
We can prove Parts (i) and (ii) by adapting the corresponding arguments for  \cite[Lemma 4.1]{DuoNguZaa19} and \cite[Lemma 2.3, (64)-(65)]{MerZaa97NL}.

\bigskip

$\bullet$ For Part (iii), we follow the streamlined arguments in \cite{DuoGhoZaa21}. 
By Duhamel formula, we can rewrite \eqref{e:qeqn} in the following integral form
$$\begin{aligned}
q(s)&=\cK(s,\sigma)q(\sigma)+\int_\sigma^s\cK(s,s')(B(q)+T(q)+R+L(q))(s')ds'\\
&=:\vartheta(s)+\Theta(s),
\end{aligned}$$
where $\cK=\cK_V$ is the fundamental solution associated to $\pd_sq=\cL_Vq$.

The following fundamental result concerns the dynamics of the negative spectrum part $v_-$ and the exterior part $v_e$
for $v$ satisfying the linear equation $\pd_sv=\cL_Vv$.
It goes back to Bricmont and Kupiainen \cite[Lemma 1]{BriKup94} (see also Merle and Zaag \cite[Lemma 3.13]{MerZaa97DMJ}) and, 
in the current form, can be found in Nguyen and Zaag \cite[Lemma 2.9]{NguZaa17} 
(its proof uses in particular \eqref{e:V1} and \eqref{e:V4} in Lemma \ref{l:V}).

\begin{lem}[Linear dynamics of $v_-$ and $v_e$ for $v$ satisfying $\pd_sv=\cL_Vv$] \label{l:BKNZ}
The exists $K^*\geq1$ such that for all $K_0\geq K^*$ and $\ell^*>0$, 
there exists $s^*=s^*(K_0,\ell^*)$ such that for all $\sigma\geq s^*$ and $v(\sigma)\in L^\infty(\bRN)$ satisfying
$$\sum_{m=0}^2|v_m(\sigma)|+\left\|\frac{v_-(\cdot,\sigma)}{1+|\cdot|^3}\right\|_{L^\infty}+\|v_e(\sigma)\|_{L^\infty}<\infty,$$
it holds that for any $s\in[\sigma,\sigma+\ell^*]$,
$$\begin{aligned}
&\left\|\frac{\mathbb{V}_-(\cdot,s)}{1+|\cdot|^3}\right\|_{L^\infty}\leq C\frac{e^{s-\sigma}((s-\sigma)^2+1)}{s}\left(|v_0(\sigma)|+|v_1(\sigma)|+\sqrt s|v_2(\sigma)|\right)\\
&\qquad\qquad\qquad\qquad+Ce^{-\frac{s-\sigma}{2}}\left\|\frac{v_-(\cdot,\sigma)}{1+|\cdot|^3}\right\|_{L^\infty}+C\frac{e^{-(s-\sigma)^2}}{s^{\frac32}}\|v_e(\sigma)\|_{L^\infty},
\end{aligned}$$
$$\begin{aligned}
&\|{\mathbb{V}_e(s)}\|_{L^\infty}\leq Ce^{s-\sigma}\left(\sum_{m=0}^2s^{\frac m2}|v_m(\sigma)|+s^{\frac32}\left\|\frac{v_-(\cdot,\sigma)}{1+|\cdot|^3}\right\|_{L^\infty}\right)\\
&\qquad\qquad\qquad\qquad+Ce^{-\frac{s-\sigma}{p}}\|v_e(\sigma)\|_{L^\infty},
\end{aligned}$$
where $\mathbb{V}(s)=\cK(s,\sigma)v(\sigma)$.
\end{lem}

With this lemma we can then proceed as follows.
Take $K_0\geq K^*$, with $K^*$ given by Lemma \ref{l:BKNZ}.
For $\ell^*>0$, take $t_5=T-e^{-s_5}$ so that 
$$s_5=s_5(K_0,\ell^*)\geq \max\{s^*(K_0,\ell^*),\ell^*\},$$ 
where $s^*(K_0,\ell^*)$ is given by Lemma \ref{l:BKNZ}.
We further take $t_0\geq t_5$ and $$\sigma\geq s_0=-\log(T-t_0)\geq s_5.$$
Note that we have the following implication for $(s', s)$:
\begin{equation} \label{e:sigmas}
\sigma\leq s'\leq s\leq \sigma+\ell^*\Longrightarrow \frac1s\leq \frac{1}{s'}\leq \frac1\sigma\leq\frac2s.
\end{equation}
According to the expression $q(s)=\vartheta(s)+\Theta(s)$, we argue in two steps.

\bigskip 

\textbf{Step I: Free evolution term $\vartheta(s)$.} Using \eqref{e:sigmas}, together with \eqref{e:rcomp} via $q(\sigma)\in V_{K_0,A}(\sigma)$,
we have the following crucial component-wise estimates:
$$|q_0(\sigma)|\leq\frac{CA}{s^2},\quad |q_1(\sigma)|\leq\frac{CA}{s^2},\quad |q_2(\sigma)|\leq \frac{CA^2\log s}{s^2},$$
$$\left\|\frac{q_-(\cdot,\sigma)}{1+|\cdot|^3}\right\|_{L^\infty}\leq\frac{CA^{\underline{\alpha}}\log s}{s^{\frac52}},\quad \|q_e(\cdot,\sigma)\|_{L^\infty}\leq\frac{CA^{\overline{\alpha}}\log s}{s}.$$
Applying Lemma \ref{l:BKNZ}, we see that
$$\left\|\frac{\vartheta_-(\cdot,s)}{1+|\cdot|^3}\right\|_{L^\infty}
\leq C\frac{\log s}{s^{\frac52}}\left(A^{\underline{\alpha}}e^{-\frac{s-\sigma}{2}}+A^{\overline{\alpha}}e^{-(s-\sigma)^2}+A^2e^{s-\sigma}\left((s-\sigma)^2+1\right)\right),$$
$$\begin{aligned}
\|{\theta_e(s)}\|_{L^\infty}\leq C\frac{\log s}{s}\left(A^{\overline{\alpha}}e^{-\frac{s-\sigma}{p}}+A^{\underline{\alpha}}e^{s-\sigma}\right).
\end{aligned}$$
If $\sigma=s_0$, we can use Proposition \ref{p:init} and obtain that
$$\left\|\frac{\vartheta_-(\cdot,s)}{1+|\cdot|^3}\right\|_{L^\infty}\leq C\frac{\log s}{s^{\frac52}}\left(1+e^{s-s_0}\left((s-s_0)^2+1\right)\right),$$ 
$$\begin{aligned}
\|{\vartheta_e(s)}\|_{L^\infty}\leq C\frac{\log s}{s}\left(1+e^{s-s_0}\right).
\end{aligned}$$
It is worthy to mention that when using the initialization Proposition \ref{p:init} instead of shrinking set conditions, 
the involved estimates are $A$-independent.

\bigskip 

\textbf{Step II: Duhamel term $\Theta(s)$.} 
We only deal with the case $\sigma>s_0$ (if $\sigma=s_0$, we use Proposition \ref{p:init} again).
By the assumptions that $h(\cdot,t_0)$ is given by \eqref{e:ht0} and $q(s)\in V_{K_0,A}(s)$, 
we can use elementary estimates on $(q,\nabla q)$ and on each term in
$$v( s'):=B(q)(s')+T(q)(s')+R(s')+L(q)(s').$$ 
These estimates are consequences of the definition of the shrinking set and are recalled in Appendices \ref{a:A}-\ref{a:B} (Lemmata \ref{l:qnablaq}-\ref{l:Tq} and \ref{l:R}-\ref{l:wtf}). 
More precisely, for 
$$\sigma\geq s_0\geq s_5\geq \ell^*,\quad \ell\leq \ell^*,\quad s\in[\sigma, \sigma+\ell]\quad \text{and} \quad s'\in[\sigma,s],$$ 
we infer from Lemma \ref{l:Bq} that
\begin{equation} \label{e:BqEst}
|B(q)(y,s')|\leq  \frac{C(1+|y|^3)}{s^{\frac52}}\quad\text{and}\quad |B(q)(y,s')|\leq \frac{C}{s},
\end{equation}
and from \eqref{e:Tq*bd}-\eqref{e:Tq*} that 
\begin{equation} \label{e:TqEst}
|T(q)(y,s')|\leq  \frac{CA^2(1+|y|^3)\log s}{s^{\frac52}}\quad\text{and}\quad |T(q)(y,s')|\leq \frac{C}{s}. 
\end{equation}
Recall that we have the following implication: for $m$ which may depend on $s$,
$$|v(y)|\leq m(1+|y|^3)\Longrightarrow|v_-(y)|\leq Cm(1+|y|^3).$$
By \eqref{e:BqEst}, \eqref{e:TqEst} and definitions of $v_0$, $v_1$ and $v_2$,
together with the estimates on $R$ and $L(q)$ from Lemmata \ref{l:R}-\ref{l:wtf}, 
we have the following component-wise estimates:
$$|v_0(s')|\leq\frac{C}{s^2},\quad |v_1(s')|\leq\frac{C}{s^2},\quad |v_2(s')|\leq\frac{C\log s}{s^2},$$
$$\left\|\frac{v_-(\cdot,s')}{1+|\cdot|^3}\right\|_{L^\infty}\leq \frac{CA^2\log s}{s^{\frac52}},\quad \|{v_e(s')}\|_{L^\infty}\leq \frac{C}{s}.$$
Recall that $s\in[s',s'+\ell^*]$, which follows from $0\leq s-s'\leq s-\sigma\leq \ell\leq \ell^*$.
Applying Lemma \ref{l:BKNZ}, then for $s$ large enough (say, $s\geq e^{2\ell^*}((\ell^*)^2+1)^2$) we see that
$$\left\|\frac{(\cK(s,s')v(\cdot,s'))_-}{1+|\cdot|^3}\right\|_{L^\infty}\leq \frac{CA^2\log s}{s^{\frac52}},$$
$$\|{(\cK(s,s')v(s'))_e}\|_{L^\infty}\leq \frac{CA^2e^{s-\sigma}\log s}{s}.$$ 

Now, we observe (as in \cite{DuoGhoZaa21}) that
$$\Theta_-(s)=P_-\left(\int_\sigma^s(\cK(s,s')v(s'))_-ds'\right),$$
hence 
$$\left\|\frac{\Theta_-(\cdot,s)}{1+|\cdot|^3}\right\|_{L^\infty}\leq \frac{CA^2(s-\sigma)\log s}{s^{\frac52}}.$$
Similarly, we observe that
$$\Theta_e(s)=\int_\sigma^s(\cK(s,s')v(s'))_eds',$$
hence 
$$\|{\Theta_e}\|_{L^\infty}\leq \frac{CA^2(s-\sigma)e^{s-\sigma}\log s}{s}.$$

\bigskip

Collecting above arguments, we proved Proposition \ref{p:qdyn}.
\end{proof}

\bigskip

Now we continue to prove Proposition \ref{p:impr1}.
We argue only with the estimates on $q_-$ 
since the other parts of estimates on $q_2$ and $q_e$ can be proved in a similar manner (cf. \cite[Pages 1525-1527]{MerZaa97NL}).
Let us consider $K_0\geq K_5$, $A\geq A_5$, $\ell^*=\log A$ and $t_0\geq t_5$ 
where $t_5=t_5(K_0, A, \ell^*)< T$ is given by Proposition \ref{p:qdyn}.
Let us assume that the hypotheses stated in Proposition \ref{p:impr1} hold for some $t_*\geq t_0$.
Set $$s_0=-\log(T-t_0), \quad s_*=-\log(T-t_*)\quad \text{and}\quad s=-\log(T-t).$$

Two cases then arise:

$\bullet$ If $s-s_0\leq \log A$, applying \eqref{e:q-2} with $\ell=s-s_0\leq \log A$ and $\sigma=s_0$ we get
$$\begin{aligned}
\left\|\frac{q_-(\cdot,s)}{1+|\cdot|^3}\right\|_{L^\infty}&\leq\frac{C\log s}{s^{\frac52}}(1+(\log A)^3 e^{\log A})\\
&\leq \frac{A^{\underline{\alpha}}\log s}{2s^{\frac52}},
\end{aligned}$$
provided that $\underline{\alpha}>1$ and $A\geq A_{5,1}$ for some $A_{5,1}>1$.

$\bullet$ If $s-s_0> \log A$, applying \eqref{e:q-1} with $\sigma=s-\log A$ and $\ell=\log A (=\ell^*)$ we get 
$$\begin{aligned}
\left\|\frac{q_-(\cdot,s)}{1+|\cdot|^3}\right\|_{L^\infty}&\leq\frac{C\log s}{s^{\frac52}} 
\left({A^{\underline{\alpha}}}{e^{-\frac{\log A }{2}}}+{A^{\overline{\alpha}}}{e^{-\log^2 A }}+A^2e^{\log A }((\log A)^2  +1)+A^2\log A \right)\\  
&\leq \frac{A^{\underline{\alpha}}\log s}{2s^{\frac52}},
\end{aligned}$$
provided that $\underline{\alpha}>3$ and $A\geq A_{5,2}$ for some $A_{5,2}>1$.
 
We can conclude the proof by taking $A_2=\max\{A_5, A_{5,1},A_{5,2}\}$.

\begin{rem}
The current proof strategy streamlined via Lemma \ref{l:BKNZ} only leads to rough self-improving bounds but is considerably simpler.
See e.g. \cite[Lemma 3.13 and Proposition 3.1]{MerZaa97DMJ}, where the bound for $q_-$ has smaller size dependence in $A$.
These weaker bounds are however enough for the finite dimensional reduction.
\end{rem}

\begin{rem}
In proving the self-improving estimate on $\|q_e(\cdot,s)\|_{L^\infty}$, 
the condition $\overline{\alpha}\geq\underline{\alpha}+1$ is required when similar two-case analysis is applied to \eqref{e:qe1} and \eqref{e:qe2}.
\end{rem}

\section{A priori estimates in $\cR_2$} \label{s:impr2}

In this section we prove Proposition \ref{p:impr2}.
To simplify the presentation, 
we assume $\tau_1=0$, $\tau_2=1$ and $N(k)=\frac{1}{k^\beta}$.
The general case follows \textit{mutatis mutandis}.

\bigskip

The proof of the first estimate on $|k-\wh k|$ is same as the one for \cite[Lemma 4.2]{MerZaa97NL}. 
For the gradient estimate, we adapt an argument from \cite[Lemma 4.2]{Duo19} which uses parabolic regularity
(similar techniques will be used again later in the proof of Proposition \ref{p:profile}).
This enables us to get rid of the second order regularity assumption in the maximum principle approach of \cite{MerZaa97NL}.
This parabolic regularity argument also provides a proof of the estimate on $|k-\wh k|$ as in \cite{Duo19}. We thus omit it.

\bigskip

Recall that on $[0,1]$ we have
$$\wh k(\tau)=\left((\beta+1)(1-\tau)+\frac{(\beta+1)^2}{4\beta}\frac{K_0^2}{16}\right)^{\frac{1}{\beta+1}}
\geq\left(\frac{(\beta+1)^2}{4\beta}\frac{K_0^2}{16}\right)^{\frac{1}{\beta+1}}.$$
Thus by \eqref{e:R2-1}, $|1/k^{\beta+1}|\leq 1$ for $K_0$ large. 
Observe that $\nabla k$ satisfies the system
$$\frac{\pd}{\pd\tau}\nabla k=\Delta(\nabla k)+\frac{\beta}{k^{\beta+1}}\nabla k.$$
Let $$\wt k(\tau)=\wh k(\tau)\wt{\ep}$$ where $\wt{\ep}$ is a small vector in $\bRN$, and 
$$g=\chi_2\left(\nabla k-\wt k(\tau)\right).$$ 
Let $\chi_2\in C^\infty(\bRN)$ such that 
$$\chi_2\equiv1\text{\,\,for\,\,}|x|\leq\frac{3|\xi_0|}{2}, \quad \chi_2\equiv0\text{\,\,for\,\,}|x|\geq\frac{7|\xi_0|}{4},\quad 0\leq\chi_2\leq1,$$ 
$$|\nabla \chi_2|\leq \frac{C}{|\xi_0|}\quad\text{and}\quad|\Delta \chi_2|\leq\frac{C}{|\xi_0|^2}.$$
Then
$$\frac{\pd g}{\pd\tau}=\Delta g+\frac{\beta}{k^{\beta+1}}g+G,$$
where the remaining term $G$ can be rewritten as follows
$$\begin{aligned}
G&=\frac{\pd g}{\pd\tau}-\frac{\beta}{k^{\beta+1}}g-\Delta g\\
&=\chi_2\Delta(\nabla k)-\Delta\left(\chi_2\left(\nabla k-\wt k(\tau)\right)\right)\\
&=-(\Delta\chi_2)(\nabla k-\wt k(\tau))-2(\nabla \chi_2) (\nabla^2 k)\\
&=-(\Delta\chi_2)(\nabla k-\wt k(\tau))-2\left[\nabla\cdot((\nabla \chi_2)\otimes(\nabla k))-(\Delta\chi_2)(\nabla k)\right]\\
&=(\Delta\chi_2)(\nabla k+\wt k(\tau))-2\nabla\cdot\left((\nabla \chi_2)\otimes(\nabla k)\right).
\end{aligned}$$
By Duhamel formula, for any $\tau\in[0,1)$,
$$g(\tau)=e^{\tau\Delta}g(0)+\int_0^\tau e^{(\tau-\tau')\Delta}\left(\frac{\beta}{k^{\beta+1}}g+G\right)(\tau')d\tau',$$
and taking the $L^\infty$-norm, we have
$$\begin{aligned}
\|g(\tau)\|_{L^\infty}&\leq\|g(0)\|_{L^\infty}+C\int_0^\tau\|g(\tau')\|_{L^\infty}d\tau'\\
&\qquad+C\int_0^\tau\left(\|(\Delta\chi_2)(\nabla k+\wt k(\tau'))\|_{L^\infty}+\frac{1}{\sqrt{\tau-\tau'}}\|(\nabla \chi_2)\otimes(\nabla k)\|_{L^\infty}\right)d\tau'\\
&\leq\frac{C_0''^*}{|\xi_0|}+ C\left(|\wt{\ep}|+\frac{C_0'^*}{|\xi_0|^2}+\frac{C_0'^*}{|\xi_0|^3}\right)+C\int_0^\tau\|g(\tau')\|_{L^\infty}d\tau'.
\end{aligned}$$
By Gr\"onwall's lemma,
$$\|g(\tau)\|_{L^\infty}\leq C \left(\frac{C_0''^*}{|\xi_0|}+ C\left(|\wt{\ep}|+\frac{C_0'^*}{|\xi_0|^2}+\frac{C_0'^*}{|\xi_0|^3}\right)\right),$$
thereby, for $\xi\in B(0,\frac{3|\xi_0|}{2})$,
$$|\nabla k(\xi,\tau)|\leq C \left(\frac{C_0''^*}{|\xi_0|}+ C\left(|\wt{\ep}|+\frac{C_0'^*}{|\xi_0|^2}+\frac{C_0'^*}{|\xi_0|^3}\right)\right)\leq \frac{C_{3}C_0''^*}{|\xi_0|},$$
where $C_{3}>1$ is a universal constant, $|\xi_0|\geq |\xi_3|$ and $\xi_3=\xi_3(C_0'^*,C_0''^*)$ is large.

\section{A priori estimates in $\cR_3$} \label{s:impr3} 

In this section we prove Proposition \ref{p:impr3} by following \cite{MerZaa97NL}.
We choose to only give the proof of the gradient estimate in \eqref{e:hR3}.
It will be justified that the second order assumption on $|\nabla^2h|$ is in fact redundant (compared to \cite[(88)]{MerZaa97NL}).
It is also for this reason that we only have to impose $W^{1,\infty}$-regularity in defining $\fH$.

\bigskip

We argue by contradiction. Let us consider $t_\ep\in(t_0,t_*)$ such that
\begin{equation} \label{e:nablah1}
|\nabla h(x,t)-\nabla h(x,t_0)|_{L^\infty(\Omega^\perp_{\ep_0})}\leq\ep,\quad\forall\,t\in[t_0,t_\ep),
\end{equation}
and
\begin{equation} \label{e:nablah2}
|\nabla h(x,t_\ep)-\nabla h(x,t_0)|_{L^\infty(\Omega^\perp_{\ep_0})}=\ep,
\end{equation}
where $\Omega^\perp_{\ep_0}=\Omega\backslash B(0,\frac{\ep_0}{4})$.
We can assume $\ep\leq\ep^*$, where
$$\ep^*=\ep^*(H^*,\ep_0)=\frac14\min_{|x|\geq\frac{\ep_0}{6}}H^*(x)>0.$$
Also, by Assumption (ii) of the proposition, 
$$\|F'(h)(\cdot,t_0)\|_{L^\infty(\Omega^\perp_{\ep_0})}\leq C(H^*,\ep_0).$$

Next we follow the same techniques used in the proof of Proposition \ref{p:impr2}, but now we work directly on $h$.
Note that $\nabla h$ satisfies the equation 
\begin{equation} \label{e:nablah}
\frac{\pd}{\pd t}\nabla h=\Delta(\nabla h)-F'(h)\nabla h.
\end{equation}
Introduce $$h_1=\chi_3\nabla h,$$ 
where $\chi_3\in C^\infty(\bRN)$ satisfies
$$\chi_3\equiv1\text{\,\,for\,\,}|x|\geq\frac{\ep_0}{5}, \quad \chi_3\equiv0\text{\,\,for\,\,}|x|\leq\frac{\ep_0}{6},\quad 0\leq\chi_3\leq1,$$  
$$|\nabla \chi_3|\leq \frac{C}{\ep_0}\quad\text{and}\quad|\Delta \chi_3|\leq\frac{C}{\ep_0^2}.$$
Then
$$\frac{\pd h_1}{\pd t}=\Delta h_1-F'(h)h_1+G,$$
where the remaining term $G$ can be rewritten as follows
$$\begin{aligned}
G&=\frac{\pd h_1}{\pd t}+F'(h)h_1-\Delta h_1\\
&=-(\Delta\chi_3)(\nabla h)-2(\nabla \chi_3) (\nabla^2 h)\\ 
&=-(\Delta\chi_3)(\nabla h)-2\left[\nabla\cdot((\nabla \chi_3)\otimes(\nabla h))-(\Delta\chi_3)(\nabla h)\right]\\
&=(\Delta\chi_3)(\nabla h)-2\nabla\cdot\left((\nabla \chi_3)\otimes(\nabla h)\right).
\end{aligned}$$
In the second line above, we used the equation \eqref{e:nablah}.
We argue by two cases.

(i) $\Omega$ is a bounded domain. For all $t\in [t_0,t_*)$,
$$h_1(t)-e^{(t-t_0)\Delta}h_1(t_0)=\int_{t_0}^t e^{(t-t')\Delta}\left(G-F'(h)h_1\right)(t')dt'=:I(t).$$
Note that 
$$\int_{t_0}^t\left\|e^{(t-t')\Delta}\nabla\cdot\right\|_{L^\infty\rightarrow L^\infty}dt'\leq C\sqrt{t-t_0}.$$ 
Moreover, the size of $G$ is seen in $|\nabla h|$ which is bounded by \eqref{e:nablah1} and \eqref{e:R3-2}.
Hence, 
$$\begin{aligned}
&\|h_1(t)-h_1(t_0)\|_{L^\infty(\Omega^\perp_{\ep_0})}\\
&\qquad\leq\|e^{(t-t_0)\Delta}h_1(t_0)-h_1(t_0)\|_{L^\infty(\Omega^\perp_{\ep_0})}+\|I(t)\|_{L^\infty(\Omega^\perp_{\ep_0})}\\
&\qquad\leq\|e^{(t-t_0)\Delta}h_1(t_0)-h_1(t_0)\|_{L^\infty(\Omega^\perp_{\ep_0})} 
+C(\ep,\ep_0,\sigma_0,\sigma_1)\left(\sqrt{t-t_0}+(t-t_0)\right).
\end{aligned}$$
If $t_0\in [t_4,T)$ for some large $t_4=t_4(\ep,\ep_0,\sigma_0,\sigma_1)$, we have
$$\|\nabla h(t_\ep)-\nabla h(t_0)\|_{L^\infty(\Omega^\perp_{\ep_0})}=\|h_1(t_\ep)-h_1(t_0)\|_{L^\infty(\Omega^\perp_{\ep_0})}\leq\frac{\ep}{2},$$
which is a contradiction to the assumption \eqref{e:nablah2}. 

(ii) $\Omega=\bRN$. The proof uses \eqref{e:HRN} (see \cite{MerZaa97NL}). We omit the details.

\section{Finite dimensional reduction} \label{s:finredu}

In this section we prove Proposition \ref{p:finredu} on Finite Dimensional Reduction. 
This is the most important part of the constructive techniques.
In principal we follow \cite{MerZaa97NL} but we shall rewrite it in a manner which we hope is more easily accessible.

\bigskip

To be reader-friendly,
we point out that the parameter subscripts ($K_i$, $A_i$, $t_i$ and so on) are numbered first with respect to the order of propositions that appear in the roadmap to Theorem \ref{t:Quen},
then the order in the proof of Proposition \ref{p:impr1} (namely, Proposition \ref{p:qdyn}),
and at last the order for the estimates collected in the Appendice \ref{a:A}-\ref{a:B}.
We continue the counting of the parameter subscripts in this section.

\subsection{First determination of parameters for the shrinking set}

We now precise our \underline{Parameter Chronology} (or ``Parameter Tracking") for the shrinking set
$$S^*(t_0,t)=S^*(t_0, K_0, \ep_0, A,\alpha_0,\delta_0, C_0,\eta_0,t),$$ that is,
we shall give a chronological determination of the main involved parameters
$$\underline{t_0, \,\, K_0,  \,\, \ep_0,  \,\, A, \,\, \alpha_0, \,\, \delta_0, \,\,  C_0, \,\, \eta_0},$$
together with other preparatory ones. We begin with the first part of parameters.

\bigskip 

$\bullet$ \underline{$K_0$ and $\delta_0$}. Fix $K_0=4\max\{K_1,K_2,K_3\}$. Here, $(K_1,K_2,K_3)$ is taken accordingly from Propositions (\ref{p:init}, \ref{p:impr1}, \ref{p:impr2}).
Due to the dependence of $\wh k(1)$ in $K_0$, we fix thereby 
$$\delta_0=\delta_0(K_0)=\frac12\min\left\{\wh k(1),1\right\}.$$

\bigskip

$\bullet$ Preparatory step on $A$ and $C_0$. Let $A_{10}=A_{10}(K_0)$ be large enough so that $A_{10}\geq \max\{A_2,2\overline C\}$ 
where $A_2=A_2(K_0)>1$ is taken from Proposition \ref{p:impr1} with $K_0$ determined in last step, 
and $\overline C$ is taken from \eqref{e:01mod}. 
For $K_0$ and $A$, we define
\begin{equation} \label{e:C0K0A}
C_0(K_0,A):=2C_3\max\left\{C_{10}A^{\overline{\alpha}}+\|\nabla \wh \Phi\|_{L^\infty(B(0,2K_0))},\frac{10\wh k(1)}{(\beta+1)K_0}, C_0^*\right\},
\end{equation}
where $C_3>1$ is given in Proposition \ref{p:impr2},
$C_{10}$ is a large constant which shall appear in the proof of Lemma \ref{l:rigidk} below,
and $C_0^*=C_0^*(K_0)$ is taken from Proposition \ref{p:init} 

\bigskip

$\bullet$ \underline{$A$ and $C_0$}. Take $A\geq A_{10}$ and consider $C_0=C_0(K_0,A)$ as defined in \eqref{e:C0K0A}.

\bigskip

$\bullet$ (Tuning) Apply Proposition \ref{p:impr2} with $K_0$, $C_0'^*=2C_0$ and $C_0''^*=\frac{C_0}{4}$.
We obtain 
$$\xi_0^*=\xi_0^*(A)\text{\,\,with\,\,}|\xi_0^*|\geq1\quad \text{and}\quad \delta_3^*=\delta_3^*(A)\text{\,\,with\,\,}0<\delta_3^*\leq1$$ 
such that for all $\xi_0$ with $|\xi_0|\geq|\xi_0^*|$ and for all $\delta_3\leq \delta_3^*$,
the conclusion of Proposition \ref{p:impr2} holds with $\ep=\frac{\delta_0}{2}>0$. 
We point out that the relation $\ep=\frac{\delta_0}{2}$ is allowable by increasing $|\xi_0|$ and decreasing $\delta_3^*$.
Indeed, in Proposition \ref{p:impr2} one first arrives the conclusion with a bound $|k(\xi,\tau_1)-\wh k(\tau_1)|\leq \ep$, 
but $\ep\rightarrow0$ as $(\delta_3,|\xi_0|)\rightarrow(0,+\infty)$. 

\bigskip

$\bullet$ Let 
\begin{equation} \label{e:delta3A}
\delta_3=\delta_3(A)=\min\left\{\frac{\delta_3^*}{2},\delta_0\right\}.
\end{equation}

\subsection{A priori estimates in $\cR_2$ via the shrinking set}

In defining $S^*(t_0,t)$ we prescribed certain regularity estimates of $k_x(\xi,\tau(x,t))$ for $|\xi|\leq \alpha_0\sqrt{|\log \theta (x)|}$, $\alpha_0>0$, 
to guarantee $ h\in  {W^{1,\infty}}$ in $\cR_2(K_0,\ep_0,t)$.
With above choices of parameters ($K_0$, $A$, $C_0$ and $\delta_3$), 
we prove in next lemma that for small $\alpha_0$ these estimates are rigid (and improved for the initial time) 
in the annuli $\alpha_0\sqrt{|\log \theta (x)|}\leq |\xi|\leq \frac74\alpha_0\sqrt{|\log \theta (x)|}$.

\begin{lem}[\textit{A priori} estimates in $\cR_2$ via the shrinking set] \label{l:rigidk}
For any $A\geq A_{10}$, 
there exists $\alpha_{10}=\alpha_{10}(K_0,\delta_3(A))>0$ such that for all $\alpha_0\leq \alpha_{10}$,
there exists $\ep_{10}=\ep_{10}(\alpha_0,A)$ such that for all $\ep_0\leq \ep_{10}$, there exists $t_{10}=t_{10}(\ep_0,A)<T$ and $\eta_{10}=\eta_{10}(\ep_0,A)>0$
such that for all $\eta_0\leq \eta_{10}$ and $t_0\in[t_{10},T)$,
if for some $t_*\in[t_0,T)$
$$h(t)\in S^*(t_0, K_0, \ep_0, A,\alpha_0,\delta_0, C_0,\eta_0,t),\quad \forall\,t\in[t_0,t_*],$$
where $h$ is the solution to \eqref{e:Quen} with initial data as in \eqref{e:ht0}, 
then, by setting
$$\quad \tau_0=\max\left\{\frac{t_0-t(x)}{\theta(x)},0\right\}\quad\text{and}\quad \tau_*=\frac{t_*-t(x)}{\theta(x)},$$
where $t(x)$ is defined in \eqref{e:tx} and $\theta(x)=T-t(x)$, we have for all $x$ with
\begin{equation} \label{e:xmax*}
|x|\in\left[\frac{K_0}{4}\sqrt{(T-t_*)|\log(T-t_*)|},\ep_0\right] 
\end{equation}
that $k_x$ defined from $h$ via \eqref{e:kx} satisfies \eqref{e:knox} and

(i) For $|\xi|\leq \frac{7}{4}\alpha_0\sqrt{|\log \theta (x)|}$ and $\tau\in[\tau_0,\tau_*]$, 
$$k_x(\xi,\tau)\geq \frac12\wh k(\tau),\quad |\nabla_\xi k_x(\xi,\tau)|\leq \frac{2C_0}{\sqrt{|\log \theta (x)|}},$$

(ii) For $|\xi|\leq 2\alpha_0\sqrt{|\log \theta (x)|}$,
$$|k_x(\xi,\tau_0)-\wh k(\tau_0)|\leq \delta_3, \quad|\nabla_\xi k_x(\xi,\tau_0)|\leq \frac{C_0}{2C_3\sqrt{|\log \theta (x)|}}.$$
Here $C_0=C_0(K_0,A)$ and $\delta_3=\delta_3(A)$ are respectively given in \eqref{e:C0K0A} and \eqref{e:delta3A}.
\end{lem}

For the reader's convenience, we recall that $C_3>1$ is given in Proposition \ref{p:impr2}.
One designs $\frac{C_0}{2C_3}$ for the use of above lemma in conjunction with Proposition \ref{p:impr2}.

\begin{proof}
Since $x$ satisfies \eqref{e:xmax*}, 
and we shall select small $\ep_0>0$ and large $t_0$ which is close to quenching time $T$,
$H^*(x)$ is then given by \eqref{e:H*}.
Moreover, the $x\rightarrow0$ asymptotics formulae in Lemma \ref{l:Htheta} are then applicable.

We only prove the gradient estimate: for $|\xi|\leq \frac{7}{4}\alpha_0\sqrt{|\log \theta (x)|}$ and $\tau\in[\tau_0,\tau_*]$,
\begin{equation} \label{e:knabRig}
|\nabla_\xi k_x(\xi,\tau)|\leq \frac{2C_0}{\sqrt{|\log \theta (x)|}};
\end{equation}
for $|\xi|\leq 2\alpha_0\sqrt{|\log \theta (x)|}$,
\begin{equation} \label{e:knabRig'}
|\nabla_\xi k_x(\xi,\tau_0)|\leq \frac{C_0}{2C_3\sqrt{|\log \theta (x)|}}.
\end{equation}
The other estimates can be proved in a similar way.

Now we just proceed as in the proof of \cite[Lemma 2.6]{MerZaa97NL}. 

\bigskip

\textbf{Proof of \eqref{e:knabRig}.}
Denote 
$$X=X_x(\xi)=x+\xi\sqrt{T-t(x)}=x+\xi\sqrt{\theta(x)}.$$
Recall that
$$t=t_x(\tau)=(T-t(x))\tau+t(x)=\theta(x)\tau+t(x).$$
Therefore,
\begin{equation} \label{e:nabkh}
\nabla_\xi k_x(\xi,\tau)=\theta(x)^{\frac{1}{2}-\frac{1}{\beta+1}}(\nabla_x h)(X,t).
\end{equation}
Moreover, $\tau\in[\tau_0,\tau_*]$ corresponds to $t\in[\max\{t_0,t(x)\},t_*]$.

Let $\epsilon>0$ to be fixed later.
If $\alpha_0\leq \frac{\epsilon K_0}{8}$, 
then for $|\xi|\leq 2\alpha_0\sqrt{|\log \theta (x)|}$,
\begin{equation} \label{e:Xsimx}
(1-\epsilon)|x|\leq |X|\leq (1+\epsilon)|x|.
\end{equation}
Similar to construct the covering regions $\{\cR_i\}$, 
here we also partition $\Omega$ into three parts, 
according to the position of the new variable $X$ with respect to $\ep_0$ and 
$$r(t):=\frac{K_0}{4}\sqrt{(T-t)|\log(T-t)|}.$$ 
Therefore, each part in the resulted partition is a sub-region in some $\cR_i$, $i=1,2,3$.

\bigskip

\textbf{Case $|X|\leq r(t)$}. We can use the shrinking set definition in $\cR_1$ (with respect to the new spatial variable $X$), 
and we have for $u$ and $\nabla u$ that
$$\begin{aligned}
&\left|(T-t)^{\frac{1}{p-1}}u(X,t)-\Phi\left(\frac{X}{\sqrt{(T-t)|\log(T-t)|}}\right)\right|\\
&\qquad=\left|q\left(\frac{X}{\sqrt{(T-t)}},-\log(T-t)\right)+\frac{\kappa}{2(p-a)|\log(T-t)|}\right|\leq \frac{CA^{\underline{\alpha}} K_0^3}{\sqrt{|\log(T-t)|}}
\end{aligned}$$
and
$$\begin{aligned}
&\left|(T-t)^{\frac{1}{p-1}+\frac12}(\nabla_xu)(X,t)-\frac{1}{\sqrt{|\log(T-t)|}}(\nabla_z\Phi)\left(\frac{X}{\sqrt{(T-t)|\log(T-t)|}}\right)\right|\\
&\qquad=\left|(\nabla_yw)\left(\frac{X}{\sqrt{(T-t)}},-\log(T-t)\right)-\frac{1}{\sqrt{|\log(T-t)|}}(\nabla_z\Phi)\left(\frac{X/\sqrt{(T-t)}}{\sqrt{|\log(T-t)|}}\right)\right|\\
&\qquad=\left|(\nabla_yq)\left(\frac{X}{\sqrt{(T-t)}},-\log(T-t)\right)\right|\leq \frac{CA^{\overline{\alpha}}\log{|\log(T-t)|}}{{|\log(T-t)|}}\leq \frac{CA^{\overline{\alpha}}}{\sqrt{|\log(T-t)|}}.
\end{aligned}$$
In above estimates we used (a weaker form of) \eqref{e:qbEst} in Lemma \ref{l:qnablaq} and \eqref{e:parreg1} in Proposition \ref{p:Sobo}.
When transformed to $h$, together with the equality
$$\nabla h=\left(-\alpha^{-\frac{\beta}{\beta+1}}\right)\frac{\nabla u}{u^{1+\frac{1}{\alpha}}},$$ 
this gives by Taylor expansion that
$$\begin{aligned}
&\left|(T-t)^{-\frac{1}{\beta+1}+\frac12}(\nabla_xh)(X,t)-\frac{1}{\sqrt{|\log(T-t)|}}(\nabla_z\wh\Phi)\left(\frac{X}{\sqrt{(T-t)|\log(T-t)|}}\right)\right|\\
&\qquad\qquad\qquad\leq \frac{C_{10}A^{\overline{\alpha}}}{\sqrt{|\log(T-t)|}},
\end{aligned}$$
where $C_{10}>1$ is a universal constant.
By \eqref{e:nabkh} and \eqref{e:C0K0A}, this further implies
$$|\nabla_\xi k_x(\xi,\tau)|\leq\left(\frac{T-t}{\theta(x)}\right)^{\frac{1}{\beta+1}-\frac12}\frac{C_0}{\sqrt{|\log(T-t)|}}.$$

Now, $t\geq t(x)$ implies $T-t\leq \theta(x)$ thereby $|\log(T-t)|\geq|\log\theta(x)|$.
If $\beta\leq1$, then 
$$|\nabla_\xi k_x(\xi,\tau)|\leq\frac{C_0}{\sqrt{|\log(T-t)|}}\leq\frac{C_0}{\sqrt{|\log\theta(x)|}}.$$
If $\beta>1$, then
\begin{equation} \label{e:thetaTt}
\left(\frac{T-t}{\theta(x)}\right)^{\frac{1}{\beta+1}-\frac12}=\left(\frac{\theta(x)}{T-t}\right)^{\frac12-\frac{1}{\beta+1}}.
\end{equation}
Since $$(1-\epsilon)|x|\leq |X|\leq r(t),$$ we have
$$|x|\leq\frac{K_0}{4(1-\epsilon)}\sqrt{(T-t)|\log(T-t)|}.$$
As $|x|\rightarrow\theta(x)=T-t(x)$ is increasing and using \eqref{e:asymp}, we have
$$\theta(x)\leq \theta\left(\frac{K_0}{4(1-\epsilon)}\sqrt{(T-t)|\log(T-t)|}\right)\sim \frac{T-t}{(1-\epsilon)^2}.$$
If $\epsilon$ is small (depending on $\beta$), then \eqref{e:thetaTt} can be bounded by 2 and \eqref{e:knabRig} holds.

\bigskip

\textbf{Case $r(t)\leq|X|\leq\ep_0$}. Note that
$$k_X(\xi,\tau)=\frac{h\left(X+\sqrt{T-t(X)}\xi,t(X)+(T-t(X))\tau\right)}{(T-t(X))^{\frac{1}{\beta+1}}}.$$
Using \eqref{e:nabkh}, we have
$$\nabla_\xi k_x(\xi,\tau)=\left(\frac{\theta(X)}{\theta(x)}\right)^{\frac{1}{\beta+1}-\frac12}\nabla_\xi k_X\left(0,\frac{t-t(X)}{\theta(X)}\right).$$
Using the shrinking set definition in $\cR_2$ for $\nabla_\xi k_X$ at $\xi=0$, we have
$$\begin{aligned}
|\nabla_\xi k_x(\xi,\tau)|&\leq\left(\frac{\theta(X)}{\theta(x)}\right)^{\frac{1}{\beta+1}-\frac12}\frac{C_0}{\sqrt{|\log\theta(X)|}}\\
&=\frac{C_0}{\sqrt{|\log\theta(x)|}}\left(\frac{\theta(X)}{\theta(x)}\right)^{\frac{1}{\beta+1}-\frac12}\frac{\sqrt{|\log\theta(x)|}}{\sqrt{|\log\theta(X)|}}.
\end{aligned}$$
By \eqref{e:Xsimx} for $\epsilon$ small enough (thereby $X$ close to $x$), 
$$\left(\frac{\theta(X)}{\theta(x)}\right)^{\frac{1}{\beta+1}-\frac12}\frac{\sqrt{|\log\theta(x)|}}{\sqrt{|\log\theta(X)|}}\leq2,$$
hence \eqref{e:knabRig} holds.

\bigskip

\textbf{Case $|X|\geq\ep_0$}. Take $\eta_0\leq \epsilon \min_{|x'|\geq\ep_0}|\nabla h(x',t_0)|$ for some $\ep$.
Using the shrinking set definition in $\cR_3$ for $\nabla_\xi h$ at $X$ and the definition of initial data in \eqref{e:ht0}, we have
$$|\nabla h(X,t)|\leq (1+\epsilon)|\nabla h(X,t_0)|\leq(1+\epsilon)|\nabla h(\gamma X,t_0)|,$$
where $\gamma=1-\epsilon$ if $\beta>1$ and $\gamma=1+\epsilon$ if $\beta\leq1$, 
which follows simply from the monotonicity of $|\nabla H^*(x)|$ in $|x|$ (see \eqref{e:nablaH*}).
Also, note that 
$$\nabla h(\gamma X,t_0)=\nabla H^*(\gamma X,t_0).$$ 
By (ii) of Lemma \ref{l:Htheta}, we have
$$|\nabla h(X,t)|\leq (1+\epsilon)\frac{10\wh k(1)}{(\beta+1)K_0}\frac{\theta(\gamma x)^{\frac{1}{\beta+1}-\frac12}}{\sqrt{|\log\theta(\gamma x)|}}.$$
Hence, by \eqref{e:nabkh},
$$\begin{aligned}
|\nabla_\xi k_x(\xi,\tau)|&\leq (1+\epsilon)\frac{10\wh k(1)}{(\beta+1)K_0}
\frac{1}{\sqrt{|\log\theta(x)|}}\left(\frac{\theta(\gamma x)}{\theta(x)}\right)^{\frac{1}{\beta+1}-\frac12}\frac{\sqrt{|\log\theta(x)|}}{\sqrt{|\log\theta(\gamma x)|}}\\
&\leq \frac{C_0/2}{\sqrt{|\log\theta(x)|}}\left(\frac{\theta(\gamma x)}{\theta(x)}\right)^{\frac{1}{\beta+1}-\frac12}\frac{\sqrt{|\log\theta(x)|}}{\sqrt{|\log\theta(\gamma x)|}}.
\end{aligned}$$
By \eqref{e:Xsimx} for $\epsilon$ small enough (thereby $\gamma$ close to 1), as in last case \eqref{e:knabRig} holds.

\bigskip

In conclusion, we can choose small $\ep$, $\alpha\leq\alpha_{10}$ and $\eta\leq\eta_{10}$ so that \eqref{e:knabRig} holds.

\bigskip

\textbf{Proof of \eqref{e:knabRig'}.} If
$$|x|\geq\frac{K_0}{4}\sqrt{(T-t_0)|\log(T-t_0)|}$$
then $t(x)\leq t_0$ and $\tau_0(x)=\frac{t_0-t(x)}{\theta(x)}$. By Proposition \ref{p:init},
$$|\nabla_\xi k_x(\xi,\tau_0(x))|\leq \frac{C_0^*}{\sqrt{|\log\theta(x)|}}\leq\frac{C_0}{2C_3\sqrt{|\log\theta(x)|}}.$$
If
$$|x|\leq\frac{K_0}{4}\sqrt{(T-t_0)|\log(T-t_0)|},$$
then $t(x)\geq t_0$ and $\tau_0(x)=0$.
Then as in the proof of \eqref{e:knabRig},
$$|\nabla_\xi k_x(\xi,0)|\leq \frac{C_{10}A^{\overline{\alpha}}+\|\nabla \wh \Phi\|_{L^\infty(B(0,2K_0))}}{\sqrt{|\log\theta(x)|}}\leq\frac{C_0}{2C_3\sqrt{|\log\theta(x)|}}.$$

\bigskip

Collecting the estimates, this proves Lemma \ref{l:rigidk}.
\end{proof}

\bigskip

\textbf{Remark.} Since we used Proposition \ref{p:Sobo} in above arguments,
the determined parameter thresholds in Lemma \ref{l:rigidk} thus depend on those with subscripts $1.5$. 

\bigskip

\subsection{Second determination of parameters for the shrinking set}

Now, based on Lemma \ref{l:rigidk}, we continue the second part of our Parameter Chronology.

\bigskip

$\bullet$ \underline{$\alpha_0$ and $\ep_0$}. Fix 
$$\alpha_0=\min\left\{\frac{\alpha_1}{2},\alpha_{10},1\right\}.$$
Fix 
$$\ep_0=\ep_0(\alpha_0)\leq\min\{\ep_1,\ep_{10}\}$$ 
such that 
$$\alpha_0\sqrt{|\log\theta({\ep_0})|}\geq|\xi_0^*|.$$
Here, $(\alpha_1,\ep_1)$ and $(\alpha_{10},\ep_{10})$ are taken respectively from Proposition \ref{p:init}\footnotemark
\footnotetext{in this proposition we used $2\alpha_0$ rather than $\alpha_0$ in the shrinking set. 
Hence we set $\alpha_0\leq\frac{\alpha_1}{2}$.} and Lemma \ref{l:rigidk},
and $|\xi_0^*|$ is the parameter obtained in the tuning step. 

\bigskip

$\bullet$ \underline{$\eta_0$}. Choose 
$$\eta_0=\eta_0(\ep_0)\leq \min\{\eta_2,\eta_{10}\}.$$
Here $\eta_{10}$ is taken from Lemma \ref{l:rigidk},
and $\eta_2$ in Proposition \ref{p:impr1} depends on $\eta_5$ in Proposition \ref{p:qdyn},
which further depends on $\eta_6$ and $\eta_8$ from Lemmata \ref{l:qnablaq} and \ref{l:Tq}.

\bigskip

$\bullet$ Let 
\begin{equation} \label{e:sig0}
\sigma_0=\frac12\wh k(1)\theta\left(\frac{\ep_0}{6}\right)^{\frac{1}{\beta+1}}
\end{equation}
and
\begin{equation} \label{e:sig1}
\sigma_1=\max\left\{\frac32\wh k(0)\theta\left(\frac{\ep_0}{4}\right)^{\frac{1}{\beta+1}},
C_0\frac{\theta\left(\frac{\ep_0}{6}\right)^{\frac{1}{\beta+1}-\frac12}}{\sqrt{|\log\theta\left(\frac{\ep_0}{6}\right)|}},
C_0\frac{\theta\left(\frac{\ep_0}{4}\right)^{\frac{1}{\beta+1}-\frac12}}{\sqrt{|\log\theta\left(\frac{\ep_0}{4}\right)|}}\right\}. 
\end{equation}
These precise forms are designed for later purpose.

\bigskip

$\bullet$ \underline{$t_0$}. This is the final parameter. We choose $t_0<T$ so that
$$t_0\geq\max\{t_1,t_{1.5},t_2,t_4,t_{10},t_{\ep_0}\},$$
where $t_{\ep_0}$ is determined by the following relation 
$$\ep_0=12\sqrt{(T-t_{\ep_0})|\log(T-t_{\ep_0})|},$$
and $t_4$ depends in particular on $\sigma_0$ and $\sigma_1$, as given in Proposition \ref{p:impr3}. 

\bigskip

All the parameters are determined so far and we continue the proof.
Assume that $h(x,t_0;d_0,d_1)$ is given by \eqref{e:ht0} with $(d_0,d_1)\in\cD(t_0,K_0,A)$. 
By Proposition \ref{p:init},
$$h(t_0)\in S^*(t_0, K_0, \ep_0, A,\alpha_0,\delta_0, C_0,\eta_0,t_0).$$
This is Part (i) of Proposition \ref{p:finredu}.

\bigskip

Assume in addition that for some $t_*\in[t_0,T)$, we have for all $t\in[t_0,t_*]$
$$h(t)\in S^*(t_0, K_0, \ep_0, A,\alpha_0,\delta_0, C_0,\eta_0,t)$$
and
 $$h(t_*)\in \pd S^*(t_0, K_0, \ep_0, A,\alpha_0,\delta_0, C_0,\eta_0,t_*).$$
The remaining proof is analogous to \cite[Pages1514-1521]{MerZaa97NL}, 
and it is a consequence of the a priori estimates from Lemma \ref{l:rigidk} and Propositions \ref{p:impr1}, \ref{p:impr2} and \ref{p:impr3}.
According to the definition of $S^*(t)$, only three cases may occur.

\bigskip 

\textbf{Case 1}. $q(s_*)\in \pd V_{K_0,A}(s_*)$.

\bigskip

By Proposition \ref{p:impr1} on the self-improving properties of $q_2$, $q_-$ and $q_e$, we deduce 
\begin{equation} \label{e:pdQA}
(q_0(s_*),q_1(s_*))\in \pd\cQ_A(s_*).
\end{equation}

\bigskip 

\textbf{Case 2}. There exist $x$ and $\xi$ such that
$$|x|\in\left[\frac{K_0}{4}\sqrt{(T-t_*)|\log(T-t_*)|},\ep_0\right]\quad\text{and}\quad|\xi|\leq \alpha_0\sqrt{|\log\theta(x)|},$$
and either
$$|k_x(\xi,\tau_*)-\wh k(\tau_*)|=\delta_0,$$
or
$$|\nabla_\xi k_x(\xi,\tau_*)|=\frac{C_0}{\sqrt{\log\theta(x)}},$$
where $\tau_*=\frac{t_*-t(x)}{\theta(x)}<1$. 
Set $|\xi_0|=\alpha_0\sqrt{|\log\theta(x)|}$. 
By the monotonicity of $|x|\mapsto \theta(x)$,
$$|\xi_0|\geq\alpha_0\sqrt{|\log\theta({\ep_0})|}\geq|\xi_0^*|,$$
where $\xi_0^*$ is obtained in the tuning step.
Note that we take $\alpha_0\leq1$, thus $|\xi_0|\leq\sqrt{|\log\theta(x)|}$. 
Recall that $\tau_0=\max\left\{\frac{t_0-t(x)}{\theta(x)},0\right\}$. 
By Lemma \ref{l:rigidk} we have:

\indent\indent\indent (i) For $|\xi|\leq \frac{7}{4}\alpha_0\sqrt{|\log \theta (x)|}$ and $\tau\in[\tau_0,\tau_*]$, 
$$k_x(\xi,\tau)\geq \frac12\wh k(\tau),$$ 
$$|\nabla_\xi k_x(\xi,\tau)|\leq \frac{2C_0}{\sqrt{|\log \theta (x)|}}=\frac{2C_0\alpha_0}{|\xi_0|}\leq\frac{2C_0}{|\xi_0|},$$

\indent\indent\indent (ii) For $|\xi|\leq 2\alpha_0\sqrt{|\log \theta (x)|}$,
$$|k_x(\xi,\tau_0)-\wh k(\tau_0)|\leq \delta_3, $$
$$|\nabla_\xi k_x(\xi,\tau_0)|\leq \frac{C_0}{2C_3\sqrt{|\log \theta (x)|}}=\frac{C_0\alpha_0}{2C_3{|\xi_0|}}.$$
By Proposition \ref{p:impr2} (applied in the tuning step), for $|\xi|\leq \alpha_0\sqrt{|\log \theta (x)|}$, we have
$$|k_x(\xi,\tau_*)-\wh k(\tau_*)|\leq\ep=\frac{\delta_0}{2}$$
and
$$|\nabla_\xi k_x(\xi,\tau_*)|\leq C_3\frac{C_0\alpha_0}{2C_3{|\xi_0|}}= \frac{C_0/2}{\sqrt{|\log \theta (x)|}},$$
which contradict with the hypotheses within this case.

\bigskip 

\textbf{Case 3}. There exists $x\in\Omega$ such that $|x|\geq\frac{\ep_0}{4}$ and either
$$|h(x,t_*)-h(x,t_0)|=\eta_0,$$
or
$$|\nabla h(x,t_*)-\nabla h(x,t_0)|=\eta_0.$$
By \eqref{e:kcomp1} and \eqref{e:kcomp2},
$$\left|k_x(0,\tau)-\wh k(\tau)\right|\leq \delta_0,
\quad|\nabla_\xi k_x(0,\tau)|\leq \frac{C_0}{\sqrt{|\log \theta (x)|}}, \quad |x|\in \left[\frac{\ep_0}{6},\frac{\ep_0}{4}\right].$$
As in \cite{MerZaa97NL}, by \eqref{e:kx} and the fact that $\delta_0\leq\frac12\wh{k}(1)\leq\frac12\wh{k}(0)$ we obtain
$$\frac12\wh{k}(1)\theta(x)^{\frac{1}{\beta+1}}\leq h(x,t)\leq \frac12\wh{k}(0)\theta(x)^{\frac{1}{\beta+1}},$$
$$|\nabla h(x,t)|\leq C_0\frac{\theta(x)^{\frac{1}{\beta+1}-\frac12}}{\sqrt{|\log \theta(x)|}}.$$
It is elementary\footnotemark
\footnotetext{the latter two components in defining $\sigma_1$ in \eqref{e:sig1} take care of both $\beta\geq1$ and $\beta\leq1$.}
 to check that these two conditions are then transformed into \eqref{e:R3-1} for previous choices of $\sigma_0$ and $\sigma_1$ in \eqref{e:sig0}-\eqref{e:sig1}. 
Moreover, by \eqref{e:ht0} and $t_0\geq t_{\ep_0}$, \eqref{e:R3-2} is also satisfied. 
As a consequence, by $t_0\geq t_4$ and Proposition \ref{p:impr3} we have
$$|h(x,t_*)-h(x,t_0)|+|\nabla h(x,t_*)-\nabla h(x,t_0)|\leq\frac{\eta_0}{2},$$ 
which contradicts with the hypotheses within this case.

\bigskip

Since only \textbf{Case 1} happens, we have \eqref{e:pdQA}. 
This proves Part (ii).

\bigskip

Now, without loss of generality, we can assume 
$$q_0(s_*)=\mu \frac{A}{s_*^2}\quad\text{for}\,\, \mu=1\,\,\text{or}\,\, -1.$$
Recall that we take $A\geq A_{10}\geq2\overline{C}$. Thus by \eqref{e:01mod} for $m=0$, 
 $$\mu\frac{d}{ds}\left(q_0(s)-\mu\frac{A}{s^2}\right)\bigg|_{s=s_*}>0.$$ 
This shows that the flow of $q_0$ is transverse outgoing on the curve of $s\mapsto\mu\frac{A}{s^2}$ at $s=s_*$,
which proves Part (iii) hence completes the proof of Proposition \ref{p:finredu}.

\section{Derivation of the gradient profile via refined shrinking rates} \label{s:profile}

Here we prove Proposition \ref{p:profile} by combining the strategies of \cite{MerZaa97NL} and \cite{DuoGhoZaa21}.

\bigskip

\textbf{Step 1. Derivation of single-point extinction}.

\bigskip

To prove that the solution $h$ extinguishes at time $T$ only at the point $x=0$,
it suffices to prove \eqref{e:0extinc!}, \eqref{e:refined1} and \eqref{e:profsame}.
The latter two are proved in next steps, and the proof of \eqref{e:0extinc!} (using the conditions \eqref{e:del0eta0}) is same as \cite{MerZaa97NL}.

\bigskip

\textbf{Step 2. Derivation of intermediate profiles}.

\bigskip
 
This is analogous to \cite{DuoGhoZaa21} since by Remark \ref{r:bdd} and Proposition \ref{p:Sobo} we have
$$\forall\,\,s\geq s_0,\quad\|q(s)\|_{L^\infty}+\|\nabla q(s)\|_{L^\infty}\leq C(K_0,A)\frac{\log s}{s},$$
and we derive \eqref{e:refined1} and \eqref{e:refined2} via \eqref{e:u-int-p1} and \eqref{e:u-int-p2}.

\bigskip

\textbf{Step 3. Derivation of final profiles}.

\bigskip

For $|\mathbf{x}|\neq0$ small so that $$t(\mathbf{x})\geq t_0,$$ 
$k_\mathbf{x}(\xi,\tau)$ as defined in \eqref{e:kx} satisfies \eqref{e:keqn} for $\tau\in[0,1)$.

\bigskip

\textbf{Step 3.1. Final extinction profile}.

\bigskip

The proof is same as \cite{MerZaa97DMJ}. We include this piece of arguments with further details.
Note that we do not need precise estimates as commented in Remark \ref{r:Quen}.

By \eqref{e:refined1} applied to $x_0=0$, $t=t(\mathbf{x})$ and
$$z=\frac{\mathbf{x}+\xi\sqrt{\theta(\mathbf{x})}}{\sqrt{\theta(\mathbf{x})|\log\theta(\mathbf{x})|}},$$ 
we have 
\begin{equation} \label{e:k01}
\sup_{|\xi|\leq|\log\theta(\mathbf{x})|^{\frac14}}\left|\frac{1}{k_\mathbf{x}}(\xi,0)-\frac{1}{\wh\Phi}
\left(\frac{\mathbf{x}+\xi\sqrt{\theta(\mathbf{x})}}{\sqrt{\theta(\mathbf{x})|\log\theta(\mathbf{x})|}}\right)\right|\rightarrow0
\end{equation} 
as $\mathbf{x}\rightarrow0$ (in this case $t=t(\mathbf{x})\rightarrow T$).

Consider $|\mathbf{x}|\neq0$ small such that 
$$2\leq|\log\theta(\mathbf{x})|^{\frac14}\leq\alpha_0\sqrt{|\log\theta(\mathbf{x})|},$$
and take $2|\xi_0|=|\log\theta(\mathbf{x})|^{\frac14}$.
By \eqref{e:kcomp1} and \eqref{e:kcomp2} we have
\begin{equation} \label{e:k02}
|k_\mathbf{x}(\xi,\tau)-\wh k(\tau)|\leq\delta_0,\quad\forall\,\,|\xi|\leq|\log\theta(\mathbf{x})|^{\frac14},\quad\forall\,\,\tau\in[0,1),
\end{equation} 
\begin{equation} \label{e:k03}
|\nabla_\xi k_\mathbf{x}(\xi,0)|\leq\frac{C_0}{\sqrt{|\log\theta(\mathbf{x})|}}\leq\frac{\alpha_0C_0}{|\log\theta(\mathbf{x})|^{\frac14}},\quad\forall\,\,|\xi|\leq|\log\theta(\mathbf{x})|^{\frac14}.
\end{equation}
As $\delta_0\leq\frac{\wh k(1)}{2}<\frac{\wh k(\tau)}{2}$, \eqref{e:k02} implies
\begin{equation} \label{e:k04}
k_\mathbf{x}(\xi,\tau)\geq\frac{\wh k(\tau)}{2},\quad\forall\,\,|\xi|\leq|\log\theta(\mathbf{x})|^{\frac14},\quad\forall\,\,\tau\in[0,1).
\end{equation} 

Note that 
$$\lim_{\mathbf{x}\rightarrow0}\frac{\mathbf{x}+\xi\sqrt{\theta(\mathbf{x})}}{\sqrt{\theta(\mathbf{x})|\log\theta(\mathbf{x})|}}=\frac{K_0}{4},$$
thus \eqref{e:k01} yields indeed, via \eqref{e:khat0}, that as $\mathbf{x}\rightarrow0$,
\begin{equation} \label{e:k02'}
|k_\mathbf{x}(\xi,0)-\wh k(0)|\leq\delta(\mathbf{x})\rightarrow0,\quad\forall\,\,|\xi|\leq|\log\theta(\mathbf{x})|^{\frac14}.
\end{equation}
By \eqref{e:k02'}, \eqref{e:k03} and \eqref{e:k04}, we apply Proposition \ref{p:impr2} and obtain that
$$\forall\,\,|\xi|\leq\frac12|\log\theta(\mathbf{x})|^{\frac14},\quad\forall\,\,\tau\in[0,1),$$ 
$$|k_\mathbf{x}(\xi,\tau)-\wh k(\tau)|\leq\ep(\mathbf{x})\rightarrow0\quad\text{as}\quad \mathbf{x}\rightarrow0.$$ 
Letting $\xi=0$ and $\tau\rightarrow1$, by \eqref{e:Hsim} we get the final extinction profile \eqref{e:profsame}.

\bigskip

\textbf{Step 3.2. Final gradient profile}.

\bigskip

Let us consider $\mathbf{x}\neq0$ with $|\mathbf{x}|$ small enough.
We define
\begin{equation} \label{e:fkhatvec}
\wh{\fk_{\mathbf{x}}}(\tau)=\frac{(\beta+1)}{2\beta}\frac{K_0}{4}\frac{\mathbf{x}}{|\mathbf{x}|} \frac{1}{\sqrt{|\log\theta(\mathbf{x})|}}
\left((\beta+1)(1-\tau)+\frac{(\beta+1)^2}{4\beta}\left(\frac{K_0}{4}\right)^2\right)^{-\frac{\beta}{\beta+1}}
\end{equation}
which solves the vectorial ODE
\begin{equation} \label{e:fkhat}
\frac{d \wh{\fk_{\mathbf{x}}}}{d \tau}=\frac{\beta}{\wh k^{\beta+1}}\wh{\fk_{\mathbf{x}}},
\end{equation}
subject to initial value condition
$$\begin{aligned}\wh{\fk_{\mathbf{x}}}(0)&=(\nabla_\xi\wh\phi)(\xi=0,\tau=0)\\
&=(\nabla\wh\Phi)|_{x=\mathbf{x}}\left(\frac{K_0}{4}\right) \times \frac{1}{\sqrt{|\log\theta(\mathbf{x})|}}\\  
&=\frac{(\beta+1)}{2\beta}\frac{K_0}{4}\frac{\mathbf{x}}{|\mathbf{x}|} \frac{1}{\sqrt{|\log\theta(\mathbf{x})|}}\left(\beta+1+\frac{(\beta+1)^2}{4\beta}\left(\frac{K_0}{4}\right)^2\right)^{-\frac{\beta}{\beta+1}}.
\end{aligned}$$
Here, $\wh k$ is defined in \eqref{e:khatsol}, 
$\wh\phi=\wh\phi_{\mathbf{x}}(\xi,\tau)$ is the rescaled function of $\wh\Phi$ as given by \eqref{e:phix}, 
and $\nabla\wh\Phi$ is defined in \eqref{e:Phihatn}\footnotemark
\footnotetext{Recall that for $x_0=0$, we have the general substitutions $$z=\frac{z}{|z|}|z|\quad\text{and}\quad\frac{z}{|z|}=\frac{x}{|x|}.$$ 
In the situation $x=\mathbf{x}$ here, this enables us to write $$z|_{|z|=\frac{K_0}{4}}=\frac{K_0}{4}\frac{\mathbf{x}}{|\mathbf{x}|}.$$}.
Note that for $\tau\in[0,1)$, we have
\begin{equation} \label{e:fkhatmag}
|\wh{\fk_{\mathbf{x}}}(\tau)|\leq \frac{C}{\sqrt{|\log\theta(\mathbf{x})|}}.
\end{equation}

We claim the following estimates: for $|\mathbf{x}|\neq0$ small enough,
\begin{equation} \label{e:claim1}
\sup_{|\xi|\leq \frac14 |\log(T-t(\mathbf{x}))|^{\frac14}, \tau\in[0,1)}\left|k_{\mathbf{x}}(\xi,\tau)-\wh k(\tau)\right|\leq \frac{C}{|\log(T-t(\mathbf{x}))|^{\frac14}},
\end{equation}
\begin{equation} \label{e:claim2}
\sup_{|\xi|\leq \frac{1}{16} |\log(T-t(\mathbf{x}))|^{\frac14}, \tau\in[0,1)}\left|\nabla_\xi k_{\mathbf{x}}(\xi,\tau)-\wh{\fk_{\mathbf{x}}}(\tau)\right|\leq \frac{C}{ |\log(T-t(\mathbf{x}))|^{\frac34}}.
\end{equation}
We point out that \eqref{e:claim1} refines \eqref{e:k02'}, thus also leading to final extinction profile.
It is independent of \eqref{e:claim2} and can be proved by adapting \cite{Duo19}. 
We choose not to prove \eqref{e:claim1} since it is a simpler version of the proof strategy below for \eqref{e:claim2}.
For this reason, we can assume that \eqref{e:claim1} is proved and we use it to prove \eqref{e:claim2}.

First, by \eqref{e:kcomp2}, for $|\mathbf{x}|\neq0$ small enough so that 
$$\frac14|\log(T-t(\mathbf{x}))|^{\frac14}\leq\alpha_0\sqrt{|\log(T-t(\mathbf{x}))|},$$
we have
\begin{equation} \label{e:nabk12}
\sup_{|\xi|\leq \frac{1}{4} |\log(T-t(\mathbf{x}))|^{\frac14}, \tau\in[0,1)}\left|\nabla_\xi k_{\mathbf{x}}(\xi,\tau)\right|\leq \frac{C_0}{|\log(T-t(\mathbf{x}))|^{\frac12}}.
\end{equation}
Following \cite{DuoGhoZaa21}, we introduce the following vectorial functions
$$\cV(\mathbf{x},\xi,\tau)=\nabla_\xi k_{\mathbf{x}}(\xi,\tau)-\wh{\fk_{\mathbf{x}}}(\tau),$$ 
$$\wt\cV(\mathbf{x},\xi,\tau)=\chi_2(\xi)\cV(\mathbf{x},\xi,\tau),$$
where we use the following cut-off
$$\chi_2(\xi)=\chi_0\left(\frac{16|\xi|}{|\log(T-t(\mathbf{x}))|^{\frac14}}\right).$$
Recall that $\chi_0$ is defined in \eqref{e:chi0}, so
\begin{equation} \label{e:chi2}
|\nabla \chi_2|\leq \frac{C}{|\log(T-t(\mathbf{x}))|^{\frac14}}, \quad |\Delta \chi_2|\leq \frac{C}{|\log(T-t(\mathbf{x}))|^{\frac12}}.
\end{equation}
Thus we deduce the following equation for $\wt\cV$:
\begin{equation} \label{e:wtcV}
\frac{\pd\wt\cV}{\pd\tau}=\Delta\wt\cV+\frac{\beta}{\wh k^{\beta+1}}\wt\cV+\chi_2\left(\frac{\beta}{k^{\beta+1}}-\frac{\beta}{\wh k^{\beta+1}}\right)\wh{\fk_{\mathbf{x}}}-G,
\end{equation}
where 
$$\begin{aligned}G(\xi,\tau)&=(\Delta\chi_2)\cV+2\left(\nabla_\xi^2k\right)(\nabla\chi_2)\\
&=(\Delta\chi_2)(\cV-2\nabla_\xi k)+2\nabla\cdot((\nabla \chi_2)\otimes(\nabla_\xi k))\\
&=-(\Delta\chi_2)\left(\wh{\fk_{\mathbf{x}}}(\tau)+\nabla_\xi k\right)+2\nabla\cdot((\nabla \chi_2)\otimes(\nabla_\xi k)).
\end{aligned}$$ 
By \eqref{e:refined2} (transformed to $k_x$ and applied to $x_0=0$ and $t=t(\mathbf{x})$), we have
\begin{equation} \label{e:wtcV0}
|\wt\cV(\mathbf{x},\cdot,0)|\leq \frac{C}{|\log(T-t(\mathbf{x}))|^{\frac34}}.
\end{equation}
By \eqref{e:fkhatmag} and \eqref{e:claim1} (valid as we explained before), we have
\begin{equation} \label{e:integ1}
\left|\chi_2\left(\frac{\beta}{k^{\beta+1}}-\frac{\beta}{\wh k^{\beta+1}}\right)\wh{\fk_{\mathbf{x}}}\right|\leq \frac{C}{|\log(T-t(\mathbf{x}))|^{\frac34}}.
\end{equation}
By \eqref{e:fkhatmag}, \eqref{e:nabk12} and \eqref{e:chi2}, we have 
\begin{equation} \label{e:integ2}
\left|(\Delta\chi_2)\left(\wh{\fk_{\mathbf{x}}}(\tau)+\nabla_\xi k\right)\right|\leq \frac{C}{|\log(T-t(\mathbf{x}))|},
\end{equation}
\begin{equation} \label{e:integ3}
|(\nabla \chi_2)\otimes(\nabla_\xi k)|\leq \frac{C}{|\log(T-t(\mathbf{x}))|^{\frac34}}.
\end{equation}
Then as in \cite{DuoGhoZaa21}, we rewrite \eqref{e:wtcV} in integral form: 
$$\wt\cV(\tau)=e^{\tau\Delta}\wt\cV(0)+\int_0^\tau e^{(\tau-\tau')\Delta}\left(\frac{\beta}{\wh k^{\beta+1}}\wt\cV
+\chi_2\left(\frac{\beta}{k^{\beta+1}}-\frac{\beta}{\wh k^{\beta+1}}\right)\wh{\fk_{\mathbf{x}}}-G\right)(\tau')d\tau'.$$
Taking the $L^\infty$ norm, 
combining the mapping properties for $e^{(\tau-\tau')\Delta}$ and $e^{(\tau-\tau')\Delta}\divv$ with above estimates \eqref{e:wtcV0} on initial data
and \eqref{e:integ1}, \eqref{e:integ2} and \eqref{e:integ3} on the integrand, and using $\tau<1$,
we arrive at the following inequality
$$\|\wt\cV(\tau)\|_{L^\infty}\leq \frac{C}{|\log(T-t(\mathbf{x}))|^{\frac34}}+\int_0^\tau\|\wt\cV(\tau')\|_{L^\infty}d\tau'.$$ 
The claim \eqref{e:claim2} is then proved upon using Gr\"onwall's lemma.

In addition to that, we apply a parabolic regularity technique as in \cite{TayZaa19} to get
$$\forall\,\, \tau\in\left[\frac12,1\right)\text{\,\,and\,\,}|\xi|\leq\frac{1}{16} |\log(T-t(\mathbf{x}))|^{\frac14},\quad |\pd_\tau \nabla_\xi k_{\mathbf{x}}(\xi,\tau)|\leq C,$$
which ensures the existence of $\lim_{\tau\rightarrow1}\nabla_\xi k_{\mathbf{x}}(\xi,\tau)$.
Taking $\xi=0$ in \eqref{e:claim2} and letting $\tau\rightarrow1$, by \eqref{e:kx}, \eqref{e:claim2} and \eqref{e:nabHsim} we obtain that as $\mathbf{x}\rightarrow0$,
$$\begin{aligned}
(\nabla h)(\mathbf{x},T)&=(T-t(\mathbf{x}))^{\frac{1}{\beta+1}-\frac12}\lim_{\tau\rightarrow1}\nabla_\xi k(\mathbf{x},0,\tau)\\
&=\frac{(T-t(\mathbf{x}))^{\frac{1}{\beta+1}-\frac12}}{\sqrt{|\log(T-t(\mathbf{x}))|}}\lim_{\tau\rightarrow1}\sqrt{|\log(T-t(\mathbf{x}))|}\nabla_\xi k(\mathbf{x},0,\tau)\\
&=\frac{(T-t(\mathbf{x}))^{\frac{1}{\beta+1}-\frac12}}{\sqrt{|\log(T-t(\mathbf{x}))|}} \frac{(\beta+1)}{2\beta}\frac{K_0}{4}\frac{\mathbf{x}}{|\mathbf{x}|} 
\left(\frac{(\beta+1)^2}{4\beta}\left(\frac{K_0}{4}\right)^2\right)^{-\frac{\beta}{\beta+1}}\\
&\sim\nabla H^*(\mathbf{x})\frac{(\beta+1)K_0}{8\wh k(1)}\frac{(\beta+1)}{2\beta}\frac{K_0}{4}\left(\wh k(1)\right)^{-\beta}\\
&\sim\nabla H^*(\mathbf{x})\frac{(\beta+1)K_0}{8\left(\wh k(1)\right)^{\beta+1}}\frac{(\beta+1)}{2\beta}\frac{K_0}{4}=\nabla H^*(\mathbf{x}). 
\end{aligned}$$
We derive \eqref{e:profgrad} and complete the proof of Proposition \ref{p:profile} thereby Theorem \ref{t:Quen}.

\begin{rem}
As commented in Remark \ref{r:Quen},
for the derivation of $(\nabla h)(\mathbf{x},T)$
it is crucial that in \eqref{e:claim2} the decay in terms of $|\log(T-t(\mathbf{x}))|^{-1}$ is larger than $\frac12$.
\end{rem}

\appendix

\section{Direct regularity estimates via the shrinking set} \label{a:A}

We collect in this appendix some elementary estimates on $(q,\nabla q)$, $B(q)$ and $T(q)$ 
that are direct (not trivial however) consequences of the shrinking set conditions.

\begin{lem}[Elementary estimates on $(q,\nabla q)$] \label{l:qnablaq} 
For all $K_0\geq1$ and $\ep_0>0$, there exists $t_6=t_6(K_0,\ep_0)<T$ such that for all $t_0\in[t_6,T)$,
for all $A\geq1$, $\alpha_0>0$, $C_0>0$, $\delta_0\leq \frac{\wh k(1)}{2}$ and $\eta_0\leq\eta_6(\ep_0)$ for some $\eta_6(\ep_0)>0$,
we have the following property:

\bigskip

If $h$ is the solution to \eqref{e:Quen}, 
generated by the initial data $h(t_0;d_0,d_1)$ defined in \eqref{e:ht0}, that further satisfies
$$\forall \,\, t\in[t_0,T],\quad h(t)\in S^*(t_0, K_0, \ep_0, A,\alpha_0,\delta_0, C_0,\eta_0,t),$$
then for some $C=C(K_0,C_0)>1$
\begin{equation} \label{e:qest}
\|q(\cdot,s)\|_{L^\infty}\leq C \frac{A^{\overline{\alpha}}\log s}{s},
\end{equation}
\begin{equation} \label{e:qest'}
\left\|\frac{ q(\cdot,s)}{1+|\cdot|^3}\right\|_{L^\infty}\leq C\frac{A^{2}\log s}{ s^2},
\end{equation}
\begin{equation} \label{e:qnablae}
\|{(\nabla q)_e(\cdot,s)}\|_{L^\infty}\leq \frac{C}{\sqrt s},
\end{equation}
where $s=-\log(T-t)$ and $q$ is defined in \eqref{e:qdefn}.
\end{lem}

\begin{proof}
By the definition of $h(t)\in S^*(t_0,t)$, then $q(s)\in V_{K_0,A}(s)$, and we have
\begin{equation} \label{e:qbEst}
\begin{aligned}
|q_b(y,s)|&\leq \chi(y,s)\left(\sum_{m=0}^2|q_m(s)||h_m(y)|+|q_-(y,s)|\right)\\
&\leq C\chi(y,s)\left(\frac{A}{s^2}(1+|y|)+\frac{A^2\log s}{s^2}(1+|y|^2)+\frac{A^{\underline{\alpha}}\log s}{s^{\frac52}}(1+|y|^3)\right)\\
&\leq C\frac{A^2\log s}{s^2}\left(1+(2K_0\sqrt s)^2\right)+C\frac{A^{\underline{\alpha}}\log s}{s^{\frac52}}\left(1+(2K_0\sqrt s)^3\right)\\
&\leq C K_0^3\frac{A^{\underline{\alpha}}\log s}{s}, 
\end{aligned}
\end{equation}
which, combined with \eqref{e:rcomp} on $q_e$, proves \eqref{e:qest}. We used $K_0\geq1$, $A\geq1$ and $\underline{\alpha}\geq2$.

Meanwhile, we infer from the derivation in \eqref{e:qbEst} that for $s$ large enough,
\begin{equation} \label{e:qbEst'}
|q_b(y,s)|\leq C\chi(y,s)\frac{A^2\log s}{s^2}(1+|y|^3),
\end{equation}
whereas by \eqref{e:rcomp} on $q_e$
\begin{equation} \label{e:qeEst'}
|q_e(y,s)|\leq(1-\chi(y,s)) \frac{A^{\overline{\alpha}}\log s}{s}\leq C(1-\chi(y,s))\frac{A^{\overline{\alpha}}\log s}{s^{\frac52}}(1+|y|^3).
\end{equation}
From above two estimates, \eqref{e:qest'} follows immediately.

For the remaining estimates, first we observe that 
\begin{equation} \label{e:wnablae}
|(\nabla (\varphi+q))_e(y,s)|\leq \frac{C(K_0,C_0)}{\sqrt s}.
\end{equation}
This estimate is taken from the proof of \cite[Lemma B.1]{MerZaa97NL},
which uses only the Parts (ii)-(iii) information of the shrinking set (hence remains valid in our setting), 
together with the exterior part information of $h(\cdot,t_0;d_0,d_1)$.
Note that the modification in Part (iii) of the shrinking set is harmless.
Moreover, the constraints 
$$\delta_0\leq \frac{\wh k(1)}{2}\quad\text{and}\quad\eta_0\leq\eta_6(\ep_0)$$ 
are required in this step.
By the following crucial estimate 
$$|(\nabla\varphi)_e(y,s)|\leq \frac{C}{\sqrt s}$$ 
(which holds in fact in $\bRN$ but this will not be used), \eqref{e:qnablae} follows immediately. 
\end{proof}

\begin{lem}[Elementary estimates on $B(q)$] \label{l:Bq}
For all $K_0\geq1$, all $A\geq1$, there exists $s_7=s_7(K_0,A)$ such that for all $s\geq s_7$,
the condition $q(s)\in V_{K_0,A}(s)$ implies
$$|\chi(y,s)B(q)(y,s)|\leq C(K_0)|q|^2$$
and
$$|B(q)|\leq C|q|^{\bar p},$$
where $\bar p=\min\{p,2\}$. 
In particular, for some $C=C(K_0)>1$,
\begin{equation} \label{e:Bq'}
|B(q)(y,s)|\leq  \frac{C(1+|y|^3)}{s^{\frac52}}\quad\text{and}\quad |B(q)(y,s)|\leq \frac{C}{s}.
\end{equation} 
\end{lem}

\begin{proof}
The first part is \cite[Lemma 3.15]{MerZaa97DMJ}.
The second part can be found in \cite[(5.18)]{DuoGhoZaa21},
where $s\geq s_7(K_0,A)$ is used to have $A$-independent $C$ in \eqref{e:Bq'}.
\end{proof}

\begin{lem}[Elementary estimates on $T(q)$] \label{l:Tq}
For all $K_0\geq1$, $A\geq1$ and $\ep_0>0$, there exists $t_8=t_8(K_0,A,\ep_0)<T$ and $\eta_8=\eta_8(\ep_0)>0$ such that for all $t_0\in[t_8,T)$, 
$\alpha_0>0$, $C_0>0$, $\delta_0\leq \frac{\wh k(1)}{2}$ and $\eta_0\leq\eta_8$, we have the following property:

\bigskip

If $h$ is the solution to \eqref{e:Quen}, 
generated by the initial data $h(t_0;d_0,d_1)$ defined in \eqref{e:ht0}, that further satisfies
$$\forall \,\, t\in[t_0,T],\quad h(t)\in S^*(t_0, K_0, \ep_0, A,\alpha_0,\delta_0, C_0,\eta_0,t),$$
then for some $C=C(K_0,C_0)>1$
\begin{equation} \label{e:Tqb}
|\chi(y,s)T(q)(s)|\leq C\chi(y,s)\left(\frac{|q|}{s}+\frac{|\nabla q|}{\sqrt s}+|\nabla q|^2\right),
\end{equation}
\begin{equation} \label{e:Tqe}
|(1-\chi(y,s))T(q)(s)|\leq \frac Cs,
\end{equation}
where $s=-\log(T-t)$ and $q$ is defined in \eqref{e:qdefn}.
\end{lem}

In using \eqref{e:Tqe}, we note the following simple estimate: for $|y|\geq K_0\sqrt s$,
\begin{equation} \label{e:Tqe'}
\frac1s\leq C(K_0)\frac{|y|^3}{s^{\frac52}}\leq C(K_0)\frac{1+|y|^3}{s^{\frac52}}.
\end{equation}

\begin{proof}
We follow the proof of Lemma B.4 in \cite{MerZaa97NL}.
Consider
$$F(\theta)=-\frac{|\nabla(\varphi+ \theta q)|^2}{\varphi+ \theta q}+\frac{|\nabla\varphi|^2}{\varphi},\quad\theta\in[0,1].$$
Let us compute
$$F'(\theta)=q\frac{|\nabla(\varphi+ \theta q)|^2}{(\varphi+ \theta q)^2}-2\frac{(\nabla q)\cdot \nabla(\varphi+ \theta q)}{\varphi+ \theta q}.$$
The smallness of $q$ (which is guaranteed by Lemma \ref{l:qnablaq}) and positivity of $\varphi$ imply the boundedness of $(\varphi+ \theta q)^{-1}$.
From $F(0)=0$ and
$$F(1)=\int_0^1F'(\theta)d\theta,$$
we have the following rough estimate (where we also use $|q|\leq 1$)
$$\begin{aligned}
|\chi(y,s)T(q)|&\leq C\sup_{\theta\in[0,1]}|\chi(y,s)F'(\theta)|\\
&\leq C\chi(y,s)\left(|q||\nabla\varphi|^2+|\nabla\varphi||\nabla q|+|\nabla q|^2\right)\\
&\leq C\chi(y,s)\left(\frac{|q|}{s}+\frac{|\nabla q|}{\sqrt s}+|\nabla q|^2\right).
\end{aligned}$$
Moreover, $C$ in \eqref{e:Tqb} depends only on $K_0$.

For \eqref{e:Tqe}, one uses the crucial estimate
$$\left|(1-\chi(y,s))\frac{|\nabla \varphi|^2}{\varphi}(y,s)\right|\leq \frac{C}{s},$$
and the shrinking set conditions in $\cR_2$ and $\cR_3$. 
This is analogous to the proof of \eqref{e:qnablae},
and the reader is referred to the proof of \cite[Lemma B.1]{MerZaa97NL} for details.
Again, the constraints $\delta_0\leq \frac{\wh k(1)}{2}$ and $\eta_0\leq\eta_8(\ep_0)$ are required in this step.
\end{proof}

\section{Estimates which are independent of the shrinking set} \label{a:B}

We need estimates on $V$, $R$ and $L$ which are independent of the shrinking set.

\begin{lem} \label{l:V}
For $s\geq s_{9}$, the potential $V$ defined via \eqref{e:qeqn'} satisfies
\begin{equation} \label{e:V1}
-p\kappa^{p-1}\leq V(y,s)\leq \frac{C}{s},
\end{equation}
\begin{equation} \label{e:V2}
|V(y,s)|\leq \frac{C(1+|y|^2)}{s}, 
\end{equation}
\begin{equation} \label{e:V3} 
\left|V(y,s)+\frac{1}{4s}(|y|^2-2N)\right|\leq \frac{C(1+|y|^4)}{s^2}, 
\end{equation}
\begin{equation} \label{e:V4} 
|\nabla^i V(y,s)|\leq \frac{C}{s^\frac{i}{2}},\quad i=1,2.
\end{equation}
\end{lem}

\begin{proof}
See \cite[Lemma B.1]{NguZaa17}.
\end{proof}

\begin{lem} \label{l:R}
For any $s\geq1$, the remainder $R$ defined via \eqref{e:qeqn'} satisfies
\begin{equation} \label{e:R1}
\|R(\cdot,s)\|_{L^\infty}\leq\frac{C}{s},
\end{equation}
\begin{equation} \label{e:R2}
\left|R(y,s)-\frac{\alpha_0}{s^2}\right|\leq \frac{C(1+|y|^4)}{s^3}, 
\end{equation}
for some $\alpha_0=\alpha_0(a,p)\in\bR$. 
In particular,
$$|R_0(s)|\leq\frac{C}{s^2},\quad |R_1(s)|\leq\frac{C}{s^{3}},\quad |R_2(s)|\leq\frac{C}{s^{3}},$$
$$\left|R_b(y,s)-\frac{\alpha_0}{s^2}\right|\leq \frac{C(K_0)(1+|y|^3)}{s^{\frac52}}\quad 
\left(\text{{\rm{\,hence,\,\,}}} \left\|\frac{R_-(\cdot, s)}{1+|\cdot|^3}\right\|_{L^\infty}\leq \frac{C(K_0)}{s^{\frac52}}\right),$$
$$\|R_e(\cdot,s)\|_{L^\infty}\leq\frac{C}{s} \quad \left(\text{{\rm{\,hence,\,\,}}}  |R_e(y,s)|\leq\frac{C(1+|y|^3)}{s^{\frac52}}\right).$$
\end{lem}

\begin{proof}
The first part can be found in Lemma B.5 of \cite{MerZaa97NL} (naturally extended to $N\geq2$).
The second part follows from \eqref{e:R1} and \eqref{e:R2}, 
the decomposition \eqref{e:rdec} and the definitions involved in this decomposition (using in particular orthogonality).
\end{proof}

\begin{lem} \label{l:wtf} 
For $\wt f$ given in \eqref{e:Grad}, we have
$$\sup_{0\leq u<\infty}|\wt f(u)|+|\wt f'(u)|\leq C.$$
In particular, 
\begin{equation} \label{e:Lq}
|L(q)(s)|\leq \frac{C}{e^{\frac{ps}{p-1}}}.
\end{equation}
\end{lem}

\begin{proof}
Note that as $u\rightarrow\infty$, by \eqref{e:QuenF1} we have
$$\wt f(u)=\alpha^{\frac{\beta}{\beta+1}}u^{1+\frac{1}{\alpha}}F(\alpha^{\frac{1}{\beta+1}}u^{-\frac{1}{\alpha}})-u^p\rightarrow0.$$
The remaining arguments (using \eqref{e:QuenF2}) can be found in Lemma B.6 of \cite{MerZaa97NL}.
\end{proof}

\section{Local well-posedness in $\fH$ of the quenching problem} \label{a:C}

We show in this Appendix the local well-posedness in $\fH$ of \eqref{e:Quen}.

\bigskip

(1) Suppose that $\Omega$ is a bounded domain in $\bRN$ 
and $e^{t\Delta}$ is the heat semigroup associated to the Dirichlet Laplacian $\Delta=\Delta_\Omega$ with boundary data $\equiv1$.
Given $0<h_0\in \fH$ with $h_0\equiv1$ on $\pd\Omega$, we set 
$$\ep_0=\inf h_0\quad\text{and}\quad M_0=\|h_0\|_{W^{1,\infty}(\Omega)}.$$
Consider $$K_0=K_0(\Omega):=\left\{\fh\in L^\infty(\Omega): \inf \fh\geq\frac{\ep_0}{2}, \|\fh\|_{W^{1,\infty}(\Omega)}\leq 2CM_0\right\},$$
where $C\geq1$ is the operator norm of the (Dirichlet) heat semigroup $e^{t\Delta}$ on $W^{1,\infty}(\Omega)$.

First, we claim that there exists $T_0=T_0(\ep_0,M_0)>0$ such that
$$\phi: h\mapsto e^{t\Delta}h_0-\int_0^t e^{(t-s)\Delta}F(h(s))ds$$
is bounded and Lipschitz on $L^\infty(0,T_0;K_0)$, with a Lipschitz constant strictly less than $1$.
To see this, take an arbitrary $T>0$ and an arbitrary $h\in L^\infty(0,T;K_0)$. 
Since 
$$\frac{\ep_0}{2}\leq h(s)\quad\text{and}\quad\|h(s)\|_{W^{1,\infty}(\Omega)}\leq 2CM_0$$
$$\qquad\qquad\Longrightarrow\sup_{s\in(0,T)}\|F(h(s))\|_{W^{1,\infty}}\leq M$$ 
for some $M=M(\ep_0,M_0,\beta,C)$ where $\beta$ is given in \eqref{e:QuenF1}, we get
$$\left\|\int_0^t e^{(t-s)\Delta}F(h(s))ds\right\|_{W^{1,\infty}}\leq C TM.$$
Meanwhile, we have 
$$\ep_0\leq e^{t\Delta}h_0\quad\text{and}\quad\|e^{t\Delta}h_0\|_{W^{1,\infty}}\leq C\|h_0\|_{W^{1,\infty}}=C M_0.$$
Thus $$T\leq \frac{1}{CM}\min\left\{\frac{\ep_0}{2}, M_0\right\}\Longrightarrow\phi(h)(t)\in K_0.$$
Now take $h_1,h_2\in K_0$. Since
$$\|F(h_1(s))-F(h_2(s))\|_{W^{1,\infty}}\leq \gamma \|h_1(s)-h_2(s)\|_{W^{1,\infty}}$$
for some $\gamma=\gamma(\ep_0,M_0,\beta,C)$, we have
$$\|\phi(h_1)(t)-\phi(h_2)(t)\|_{W^{1,\infty}}\leq T\gamma \sup_{s\in(0,t)}\|h_1(s)-h_2(s)\|_{W^{1,\infty}}.$$
Choose $T_0$ so that $CT_0M\leq \min(\frac{\ep_0}{2}, M_0)$ and $T_0\gamma\leq\frac12$ proves the Lipschitz property.

By the claim just proved, for any $h_0\in\fH$ there exists a unique solution $h$ with trajectory in $K_0$, 
with maximal existence time denoted by $T$.
Moreover, either $T=\infty$ (the solution is global), or $T<\infty$, the solution escapes from $L^\infty(0,T;K_0)$ in the sense
$$\text{either}\quad\lim_{t\rightarrow T}\|h(t)\|_{W^{1,\infty}}=\infty,\quad\text{or}\quad \lim_{t\rightarrow T}\|1/h(t)\|_{L^\infty}=\infty.$$
By the arguments above, only the second scenario is allowed, i.e., $\lim_{t\rightarrow T}\inf h(t)=0$.

\bigskip

(2) Suppose that $\Omega=\bRN$ and $x_0\in\Omega$. 
Consider $$\wt h=h-\psi\quad \text{and} \quad\wt h_0=h_0-\psi,$$ 
where $\psi=\psi_{x_0}$ is given in Definition \ref{d:fH} for $\fH=\fH_{\psi_{x_0}}$.
Note that $$h_0\in\fH\Longrightarrow\wt h_0\in W^{1,\infty}(\bRN),$$
$$h(s)\in\fH\Longrightarrow \frac{1}{\wt h(s)+\psi}\in L^\infty(\bRN).$$
Moreover, \eqref{e:Quen} is equivalent to
$$\frac{\pd \wt h}{\pd t}=\Delta \wt h-F(\wt h+\psi)+\Delta \psi=:\Delta \wt h-\wt F(\wt h).$$
Here $\Delta \psi\in L^\infty(\bRN)$.
We then apply similar arguments as above to the new smooth nonlinear function $\wt F$ to solve $\wt h$ locally in $K_0(\bRN)$.
We conclude that the solution $h=\wt h+\psi$ solved this way is either global, or for some $T<\infty$, $\lim_{t\rightarrow T}\inf h(t)=0$.

\bigskip

Combing the two cases, we proved the local well-posedness of \eqref{e:Quen} in $\fH$.

\begin{rem}
In the first case of bounded domains, the local well-posedness arguments extend to other types of boundary conditions.
\end{rem}

\bibliographystyle{alpha}

\bibliography{Hua-Zaa-BlowupNLH} 
 
\end{document}